\documentclass{article}

\usepackage{amssymb}

\usepackage{proof}
\usepackage{latexsym}

\newcommand{\KP}{{\sf KP}}
\newcommand{\dg}{\mbox{{\rm dg}}}

\newcommand{\Natural}{\mathbb{N}}

\newcommand{\blem}{\begin{lemma}}
\newcommand{\elem}{\end{lemma}}
\newcommand{\bth}{\begin{theorem}}
\newcommand{\ethm}{\end{theorem}}
\newcommand{\benu}{\begin{enumerate}}
\newcommand{\eenu}{\end{enumerate}}
\newcommand{\bdes}{\begin{description}}
\newcommand{\edes}{\end{description}}
\newcommand{\bdf}{\begin{definition}}
\newcommand{\edf}{\end{definition}}
\newcommand{\bcor}{\begin{cor}}
\newcommand{\ecor}{\end{cor}}
\newcommand{\bprp}{\begin{proposition}}
\newcommand{\eprp}{\end{proposition}}
\newcommand{\bclm}{\begin{claim}}
\newcommand{\eclm}{\end{claim}}
\newcommand{\brem}{\begin{remark}}
\newcommand{\erem}{\end{remark}}
\newcommand{\bprf}{{\bf Proof}.\hspace{2mm}}
\newcommand{\eprf}{\hspace*{\fill} $\Box$}
\newcommand{\ovl}{\overline}
\newcommand{\beqn}{\begin{equation}}
\newcommand{\eeqn}{\end{equation}}
\newcommand{\beqnarr}{\begin{eqnarray}}
\newcommand{\eeqnarr}{\end{eqnarray}}
\newcommand{\beqnarrs}{\begin{eqnarray*}}
\newcommand{\eeqnarrs}{\end{eqnarray*}}
\newcommand{\spand}{\,\&\,}

\newcommand{\restrict}{\!\upharpoonright\!}

\newtheorem{theorem}{Theorem}[section]
\newtheorem{definition}[theorem]{Definition}
\newtheorem{proposition}[theorem]{Proposition}
\newtheorem{lemma}[theorem]{Lemma}
\newtheorem{cor}[theorem]{Corollary}
\newtheorem{claim}[theorem]{Claim}
\newtheorem{remark}[theorem]{Remark}

\newcommand{\alp}{\alpha}

\newcommand{\veps}{\varepsilon}

\newcommand{\Del}{\Delta}
\newcommand{\ome}{\omega}

\newcommand{\bet}{\beta}
\newcommand{\gam}{\gamma}
\newcommand{\Gam}{\Gamma}
\newcommand{\kap}{\kappa}
\newcommand{\sig}{\sigma}
\newcommand{\Sig}{\Sigma}
\newcommand{\tht}{\theta}

\newcommand{\lam}{\lambda}
\newcommand{\Lam}{\Lambda}
\newcommand{\vphi}{\varphi}

\newcommand{\fal}{\forall}
\newcommand{\exi}{\exists}

\newcommand{\Rarw }{\Rightarrow}

\newcommand{\lrarw}{\leftrightarrow}
\newcommand{\Lrarw}{\Leftrightarrow}

\newcommand{\cala}{{\cal A}}

\newcommand{\cald}{{\cal D}}

\newcommand{\calf}{{\cal F}}

\newcommand{\calh}{{\cal H}}

\newcommand{\calL}{{\cal L}}

\newcommand{\calP}{{\cal P}}

\newcommand{\la}{\langle}
\newcommand{\ra}{\rangle}

\newcommand{\msten}{\mbox{\hspace{10mm}}}
\newcommand{\msfiv}{\mbox{\hspace{5mm}}}

\title{
Hydras for $\ome_{1}$
}

\author{Toshiyasu Arai
\\
Graduate School of Science,
Chiba University
\\
1-33, Yayoi-cho, Inage-ku,
Chiba, 263-8522, JAPAN
\\
tosarai@faculty.chiba-u.jp
}

\begin{document}
\maketitle
\date{}
%\date{Jan 28-Feb. 19, 2017}

\begin{abstract}
%% Text of abstract
In this paper we introduce a hydra battle.
Each hydra will eventually die out, but the fact is not provable in a set theory with urelements of 
natural numbers and the assumption that
`there exists an uncountable regular ordinal'.
\end{abstract}

%% main text

\section{Introduction}\label{sect:intro}

L. Kirby and J. Paris \cite{Kirby-Paris82} introduces hydra battles for the first-order arithmetic PA, and
W. Buchholz \cite{BuchholzAPAL87} extends it for the theory $\mbox{ID}_{n}$ 
of $n$-th fold iterated positive elementary inductive definitions over natural numbers $(n\leq\ome)$.
The termination of hydra battles is independent from $\mbox{ID}_{n}$
since the number of moves in the battles dominates every provably recursive function in the theories.

In this paper let us modify the hydra battle of Buchholz.
Our modification is not elegant, but close to finite proof figures in set theories.

Each hydra $a\in H_{0}(\calf_{0})$ defined in Definition \ref{df:hydradf} is a finitary object, i.e., 
a term over a fixed finite set of function symbols and a finite subset $\calf_{0}$ of a set $\calf_{\mu}$ 
of function symbols defined by a $\mu$-operator.
The set $H_{0}(\calf_{0})$ is a computable subset of $\Natural$ under a suitable encoding, i.e.,
a hydra is identified with its code, a natural number.
Each hydra $a$ denotes an ordinal $v(a)$, the \textit{value} of the ordinal term $a$.
The relation $v(a)=\alp$ is defined recursively on $a$.
To codify the relation $v(a)=\alp$, we need 
a finite set of pairs $\la b,\bet\ra$ of natural numbers (hydras) $b$ and ordinals $\bet$.

For hydras $a,b$ and a natural number $n$, a stepping-down relation $b\in a[n]$
is defined in Definition \ref{df:dom},
where $a[n]$ denotes a finite set of the possible responses of the hydra $a$
under the Hercules' chop of the right-most head, relative to $\calf_{0}$ and $n$. 
The stepping-down relation is defined by referring the values of hydras.
The ternary relation $\{(a,n,b)\in\Natural^{3}: b\in a[n]\}$ on integers is far from being computable.
To define the relation, objects in three types are utilized: natural numbers $\Natural$, ordinals, and hereditarily finite sets
of natural numbers and ordinals.

Each hydra will eventually die out, but the fact is not provable in a set theory
$T_{1}(\Natural)$.
The intended model of $T_{1}(\Natural)$ is the set of hereditarily finite sets ${\rm HF}_{\Natural\cup L_{\rho_{0}}}$
over urelements $\Natural\cup L_{\rho_{0}}$, a disjoint union of the set $\Natural$ of natural numbers
and the $\rho_{0}$-th level $L_{\rho_{0}}$ of constructible sets, where $\rho_{0}>\ome_{1}$ is an ordinal
such that $L_{\rho_{0}}\models(\Pi_{1}\mbox{{\rm -Collection}})$.
Thus the set-theoretic part $L_{\rho_{0}}$ of the urelements is a model
of a set theory $T_{1}:=\KP\ome+(V=L)+(\Pi_{1}\mbox{{\rm -Collection}})+(\ome_{1})$ analyzed in \cite{CEomega1},
the Kripke-Platek set theory with the axiom of Infinity,
the axiom of constructibility $V=L$,
the $\Pi_{1}$-Collection
and the axiom $(\ome_{1})$ stating that `there exists an uncountable regular ordinal'.

For hydras $a\in H_{0}(\calf_{0})$ let $h_{a}^{\calf_{0}}(n)=k$ for the least number $k\geq n$ 
such that $a[n][n+1]\cdots[k-1]=\{0\}$
with the zero hydra $0$ if such a $k$ exists, where
$a[n][m]:=\bigcup\{b[m]:b\in a[n]\}$.

Now our theorem is stated as follows, where the class of initial hydras, a subclass of $H_{0}(\calf_{0})$,
is defined in Definition \ref{df:WB}, and
$T_{1}^{+}(\Natural)=T_{1}(\Natural)+TI(\veps_{\rho_{0}+1})$ denotes a theory obtained from $T_{1}(\Natural)$ 
by adding the transfinite induction schema
along a $\Del$-well ordering $<^{\veps}$ of type $\varepsilon_{\rho_{0}+1}$ 
for the order type $\rho_{0}$ of the class $Ord$ in a transitive and wellfounded model $V$ 
of the Kripke-Platek set theory with the axiom of the infinity.
For natural numbers $n$,
$\ome_{n}(\rho_{0}+1)$ is defined recursively by
$\ome_{0}(\rho_{0}+1)=\rho_{0}+1$ and $\ome_{n+1}(\rho_{0}+1)=\ome^{\ome_{n}(\rho_{0}+1)}$.

\bth\label{th:hydra}
\benu
\item\label{th:hydra-1}
$T_{1}^{+}(\Natural)$ proves 
the statement $\mbox{{\rm (H)}}_{\ome_{1}}$, which says 
that
the number-theoretic function $h_{a}^{\calf}$ is totally defined for {\rm any} initial hydra $a$ and any finite set 
$\calf\subset\calf_{\mu}$ of function symbols.

\item\label{th:hydra0}
$T_{1}(\Natural)$ proves the statement $\mbox{{\rm (H)}}_{\ome_{1}}$ for {\rm each} initial hydra $a$.
%, i.e.,
%$T_{1}$ proves that for {\rm each} initial hydra $a_{0}$,
%$h_{a_{0}}^{\calf_{0}}$ is totally defined
%with respect to any finite set $\calf_{0}\subset\calf_{\mu}$ of function symbols.

\item\label{th:hydra1}
$T_{1}(\Natural)$ does not prove the full statement $\mbox{{\rm (H)}}_{\ome_{1}}$.

\eenu

\end{theorem}

Let us give a sketch of our prof of the unprovability result.
Theorem \ref{th:hydra}.\ref{th:hydra1} follows from the fact that
each $T_{1}(\Natural)$-provably total $\Sig_{2}$-functions on integers
is dominated by a function $1+h_{a}^{\calf_{0}}$ for an initial hydra $a$, cf.\,Lemma \ref{lem:consis}.
As contrasted with the proofs in \cite{BuchholzAPAL87,slowgrow},
our proof of the fact is in the scheme
of the consistency proofs in G. Gentzen\cite{Gentzen38} and in G. Takeuti\cite{Takeuti67},
in which ordinals $o(\calP)$ are associated with finite proof figures $\calP$ in such a way that
$o(\calP)>o(r(\calP))$ for a proof figure $r(\calP)$ (of a contradiction).
Similarly finitary objects such as finite proof figures and hydras (ordinal terms) are manipulated
to yield a rewriting step $r$ for finite proof figures $\calP$ of a sentence $\fal n\exi m\, R(n,m)$ 
with a $\Sig_{2}$-formula $R$.
Given a proof figure $\calP_{n}$ of a sentence $\exi m\, R(\bar{n},m)$ 
with a numeral $\bar{n}$,
proof figures $\calP_{k}$ are defined recursively by $\calP_{k+1}=r(\calP_{k})$. 
Assuming that $\exi m\, R(\bar{n},m)$ does not hold,
the series $\{\alp_{k}\}_{k\geq n}$ of ordinals $\alp_{k}=o(\calP_{k})$ would be an infinite descending chain, and
hence $\exi m\, R(\bar{n},m)$ has to be true.

Although each ordinal term can be regarded as a hydra,
the rewriting step $r$ in \cite{Gentzen38, Takeuti67} does not yield a stepping-down
on ordinals such as $\alp_{k+1}=o(r(\calP_{k}))=(o(\calP_{k}))[k]=\alp_{k}[k]$.
We need to modify the rewriting step $r(\calP)$ and ordinal assignment $o(\calP)$
in such a way that 
the response of hydras to Hercules' chop corresponds exactly to
a rewriting step on proof figures.
Thus our proof here is similar to one in \cite{pntind}.

Assuming that such a rewriting step $r(\calP)$ has been defined, 
we see that $1+h_{\alp}^{\calf_{0}}(n)$ is a bound on the witnesses of the sentence $\exi m\, R(\bar{n},m)$ as follows,
We can assume that any numeral $\bar{m}$
occurring in the $k$-th proof figure $\calP_{k}$ is less than $\max\{n_{0},1+n\}+k-n$ for a constant $n_{0}$
independent from $n$.
Thus for $n\geq n_{0}$,
we can find
a witnessing numeral $\bar{m}$ of the true sentence $\exi m\, R(\bar{n},m)$
such that
$m\leq 1+h_{\alp}^{\calf_{0}}(n)$ with $\alp=o(\calP)$ and a finite set $\calf_{0}$.

In \cite{Gentzen38, Takeuti67} both of rewriting step $r$ and ordinal assignment $o$
are primitive (or even elementary) recursive.
Our rewriting step $r$ is again far from being computable.
\\

\noindent
Let us mention the contents of the paper.
In section \ref{sect:Clps}
let us recall $\Sig_{1}$-Skolem hulls, a paraphrase of the regularity of ordinals, and 
ordinals for regular ordinals.
All of these come from \cite{liftupZF} with modifications for urelements.
In section \ref{sect:hydra}, the sets of hydras $a$ and their responses $a[z]$ to Hercules' chop are introduced.
From these a hydra battle is defined.
%and the main theorem \ref{th:hydra} is stated.
Theorems \ref{th:hydra}.\ref{th:hydra-1} and \ref{th:hydra}.\ref{th:hydra0} are readily seen.
In section \ref{sect:finiteproof}, 
permissible ordinal assignments (hydras) $o(\Gam)$ to sequents $\Gam$ occurring in proofs are defined, and
each proof is shown to have a permissible ordinal assignment.
Finally we define a rewriting step $\mathcal{P}\mapsto\mathcal{P}^{\prime}$ on (finite) proofs for which
there exists a permissible ordinal assignment $o^{\prime}$ such that
$o^{\prime}(\mathcal{P}^{\prime})=(o(\mathcal{P}))[n]$ for each ordinal assignment $o$ for $\mathcal{P}$, and a proof of Theorem \ref{th:hydra}.\ref{th:hydra1} is concluded in section \ref{sect:consisprf}.

\section{The theory $T_{1}(\Natural)$}\label{sect:Clps}
In this section
the theory $T_{1}(\Natural)$ is defined.

Let $Ord$ denote the class of all ordinals.
$\rho_{0}$ denotes the least ordinal above $\ome_{1}$ such that 
$L_{\rho_{0}}\models (\Pi_{1}\mbox{{\rm -Collection}})$.
For $X\subset L_{\rho_{0}}$,
$\mbox{{\rm Hull}}(X)$ denotes the $\Sig_{1}$-Skolem hull of $X$ in $L_{\rho_{0}}$.
The Mostowski collapsing function
\[
F_{X}:\mbox{{\rm Hull}}(X)\lrarw  L_{\gam}
\]
for an ordinal $\gam\leq\rho_{0}$ such that $F_{X}\restrict Y=id\restrict Y$ 
for any transitive
$Y\subset \mbox{{\rm Hull}}(X)$.
Let us denote, though $\rho_{0}\not\in dom(F)=\mbox{{\rm Hull}}(X)$
\[
F_{X}(\rho_{0}):=\gam
.\]
The following theory $T(\ome_{1})$ is a conservative extension of the theory 
$T_{1}:=\KP\ome+(V=L)+(\Pi_{1}\mbox{{\rm -Collection}})+(\ome_{1})$, cf.\,\cite{CEomega1}.

\bdf\label{df:regext}
$T(\ome_{1})$ {\rm denotes the set theory defined as follows.}
{\rm Its language is} $\{\in, P,P_{\rho_{0}},\ome_{1}\}$ {\rm for a binary predicate} $P$, {\rm a unary predicate} $P_{\rho_{0}}$
 {\rm  and an individual constant} $\ome_{1}$.

{\rm Its axioms are obtained from those of} $\KP\ome+(\Pi_{1}\mbox{{\rm -Collection}})$ 
{\rm
 in the expanded language,
 %\footnote{
%This means that the predicates $P,P_{\rho_{0}}$ do not occur in 
%$\Del_{0}$-formulas
%for $\Del_{0}$-Separation and 
%$\Pi_{1}$-formulas
%for $\Pi_{1}$-Collection.
%},
the axiom of constructibility}
$V=L$
{\rm together with the following axiom schemata.}
{\rm For a formula} $\vphi$ {\rm and an ordinal} $\alp$,
$\vphi^{\alp}$ {\rm denotes the result of restricting every unbounded quantifier}
$\exi z,\fal z$ {\rm in} $\vphi$ {\rm to} $\exi z\in L_{\alp}, \fal z\in L_{\alp}$.

 $x\in Ord$ {\rm is a} $\Del_{0}${\rm -formula saying that `}$x$ {\rm is an ordinal'.}
 \\
$ (\ome<\ome_{1}\in Ord)$, $(P(x,y) \to \{x,y\}\subset Ord \land  x<y<\ome_{1})$
{\rm and}
$(P_{\rho_{0}}(x) \to x\in Ord)$.

\beqn\label{eq:Z1}
P(x,y) \to a\in L_{x}  \to \vphi[\ome_{1},a] \to \vphi^{y}[x,a]
\eeqn
{\rm for any} $\Sig_{1}${\rm -formula} $\vphi$ {\rm in the language} $\{\in\}$.

.\beqn\label{eq:Z2}
a\in Ord\cap\ome_{1} \to \exi x, y\in Ord\cap\ome_{1}[a<x\land P(x,y)]
\eeqn

\beqn\label{eq:Z4}
P_{\rho_{0}}(x) \to a\in L_{x} \to \vphi[a] \to \vphi^{x}[a]
\eeqn
{\rm for any} $\Sig_{1}${\rm -formula} $\vphi$ {\rm in the language} $\{\in\}$.

\beqn\label{eq:Z5}
a\in Ord \to \exi x\in Ord[a<x\land P_{\rho_{0}}(x)]
\eeqn

\edf

\bprp\label{prp:Pi1}
For a $\Del_{0}$-formula $\tht(u,v,w)$ in the language $\{\in\}$,
\[
T(\ome_{1})\vdash \fal w\,\tht(u,v,w)\lrarw\exi x\in P_{\rho_{0}}\tau(x,u,v)
\]
where
$\tau(x,u,v)\equiv[u,v\in L_{x} \land \fal w\in L_{x} \tht(u,v,w)]$.
\eprp
\bprf
This is seen from $(V=L)$, (\ref{eq:Z4}) and (\ref{eq:Z5}).
\eprf
\\
 
Let $tran(c):\equiv(\fal x\in c(x\subset c))$.
$\Pi_{1}$-Collection
\[
\fal u\in a\exi v \fal w\,\tht \to\exi c[tran(c)\land a\in c\land \fal u\in a\exi v\in c\fal w\,\tht]
\]
follows from
\beqn\label{eq:Z6}
\fal u\in a\, A(u) \to \exi c[tran(c)\land a\in c\land \fal u\in a\, A^{(c)}(u)]
\eeqn
where $A(u)\equiv(\exi x\in P_{\rho_{0}}\exi v\, \tau(x,u,v))$ for  
$A^{(c)}(u)\equiv(\exi x\in P_{\rho_{0}}\cap c\exi v\in c\,\tau)$.
\\

Next let us interpret the set theory $T(\ome_{1})$ in a theory $T^{ord}(\ome_{1})$ of ordinals
as in \cite{ptMahlo}.
The base language is $\mathcal{L}_{0}=\{<,0,+,\cdot,\lam x.\ome^{x}\}$.
Each of functions $1,\max$ and the G\"odel pairing function $j$ is $\Del_{0}$-definable in $\mathcal{L}_{0}$, 
cf.\,Appendix B
of \cite{ptMahlo}.
For each bounded formula $\cala(X,a,b)$ in the base language $\mathcal{L}_{0}$, introduce a binary predicate symbol $R^{\cala}$
with its defining axiom $b\in R^{\cala}_{a}:\equiv R^{\cala}(a,b) \lrarw \cala(R^{\cala}_{<a},a,b)$
where $c\in R^{\cala}_{<a}:\Lrarw \exi d<a(c\in R^{\cala}_{d})$.
$\mathcal{L}_{1}$ denotes the resulting language with these predicates $R^{\cala}$.
$\KP\ome+(V=L)$ is interpretable in a theory $T_{2}$ with the axiom for $\Pi_{2}$-reflection, cf.\,Appendix A
of \cite{ptMahlo}.
Each epsilon number $\alp$ is identified with the $\mathcal{L}_{1}$-structure
$\la\alp; <,0,+,\cdot,\lam x.\ome^{x}, R^{\cala}\ra$.
A G\"odel's surjective map $F:Ord\to L$ maps each epsilon number (or even a multiplicative principal number) $\alp$ onto $L_{\alp}$,
and $a\epsilon b\Lrarw F(a)\in F(b)\,(a,b\in L_{\alp})$ is a $\Del_{0}$-relation in the language $\mathcal{L}_{1}$.

For $\Pi_{2}$-formula $A$ in the language $\mathcal{L}_{1}$,
$A(t) \to \exi y[t<y\land A^{(y)}(t)]$
is an instance of $\Pi_{2}$-reflection, which follows from $(V=L)$ and $\Del_{0}$-Collection,
where $A^{(y)}$ denotes the result of restricting unbounded quantifiers $Q x\,(Q\in\{\exi,\fal\})$ to $Qx<y$.

The language of the theory $T^{ord}(\ome_{1})$ is defined to be 
$\mathcal{L}_{2}=\mathcal{L}_{1}\cup\{\ome_{1},P,P_{\rho_{0}}\}$.
The axiom (\ref{eq:Z1}) is translated to
\beqn\label{eq:Z1ord}
P(x,y) \to a<x \to \vphi[\ome_{1},a] \to \vphi^{y}[x,a]
\eeqn
 for $\Sig_{1}$-formulas $\vphi$ in $\mathcal{L}_{1}$.
The axiom (\ref{eq:Z2}) becomes
\beqn\label{eq:Z2ord}
a<\ome_{1} \to \exi x, y<\ome_{1}[a<x\land P(x,y)]
\eeqn
The axiom (\ref{eq:Z4}) turns to
\beqn\label{eq:Z4ord}
P_{\rho_{0}}(x) \to a<x \to \vphi[a] \to \vphi^{x}[a]
\eeqn
The axiom (\ref{eq:Z5}) is formulated in
\beqn\label{eq:Z5ord}
\exi x[a<x\land P_{\rho_{0}}(x)]
\eeqn

Finally consider $\Pi_{1}$-Collection.
For a $\Del_{0}$-formula $\tht(u,v,w)$ in the language $\{\in\}$,
let $\tau(x,u,v)\equiv[u,v\in L_{x}\land \fal w\in L_{x}\tht(u,v,w)]$.
Then we see $\fal w\,\tht(u,v,w)\lrarw\exi x\in P_{\rho_{0}}\tau(x,u,v)$ from $(V=L)$, (\ref{eq:Z4}) and (\ref{eq:Z5}).
Hence $\Pi_{1}$-Collection
\[
\fal u\in a\exi v \fal w\,\tht \to\exi c[a\in c\land \fal u\in a\exi v\in c\fal w\,\tht]\:
(c \mbox{ is transitive})
\]
follows from
\[
\fal u\in a\, A(u) \to \exi c[a\in c\land \fal u\in a\, A^{(c)}(u)]
\]
where $A(u)\equiv(\exi x\in P_{\rho_{0}}\exi v\, \tau(x,u,v))$ for  
$A^{(c)}(u)\equiv(\exi x\in P_{\rho_{0}}\cap c\exi v\in c\,\tau)$.
The latter is translated in the language $\mathcal{L}_{2}$ to
\beqn\label{eq:Z6ord}
\fal u<a\, A(u) \to \exi c>a\fal u<a\, A^{(c)}(u)
\eeqn
where $A(u)\equiv(\exi x\in P_{\rho_{0}}\exi v\, \tau(x,u,v))$ with a $\Del_{0}$-formula $\tau$ in $\mathcal{L}_{1}$.

Let $T^{ord}(\ome_{1})$ denote the resulting extension of the theory $T_{2}$ of ordinals 
with axioms (\ref{eq:Z1ord}), (\ref{eq:Z2ord}), (\ref{eq:Z4ord}), (\ref{eq:Z5ord}) and (\ref{eq:Z6ord}),
in which $T(\ome_{1})$ is interpreted.
\\

%\subsection{$T_{1}(\Natural)$}

$\mathcal{L}(\sf{PA})$ denotes a language for the first-order arithmetic with an individual constant 
$0^{N}$,
a unary function symbol $S$ for the successor, and relation symbols for primitive recursive relations.
$<^{N}$ denotes the less than relation on integers.
Let $N, ON$ be unary relation symbols, $\emptyset$ an individual constant,
 and $J$ a binary function symbol.
$\mathcal{L}(\Natural,\in)=\mathcal{L}({\sf PA})\cup\mathcal{L}_{2}\cup\{N,ON,Set\}\cup\{\in, =,\emptyset,J\}$ 
denotes the language for the set theory $T_{1}(\Natural)$ with urelements in $\Natural\cup \rho_{0}$,
where $J(a,x)=a\cup\{x\}$ for sets $a$, and urelements or sets $x$.
For a collection $A$ of sets over the urelements,
$A_{\Natural}=\la\Natural\cup\rho_{0}; A, \in_{A}\ra$ is a standard structure for the language, where
$\Natural$ is the standard model of the first order arithmetic,
$L_{\rho_{0}}\models(\Pi_{1}\mbox{-Collection})$ with $\rho_{0}>\ome_{1}$, and
$\in_{A}=\{(x,y)\in(\Natural\cup\rho_{0}\cup A)\times A: x\in y\}$.
The relation symbol $N$ denotes the collection $\Natural$, and $ON$ the collection $\rho_{0}$ of 
of urelemets in the structure.

The axioms in $T_{1}(\Natural)$ are classified into four groups.
\benu
\item(Ontological axioms)
Equality axioms, and
$\fal x(N(x)\veebar ON(x)\veebar Set(x))$ with the exclusive disjunction $\veebar$.

Variables $n,m,a,b,c,\ldots$ range over urelements in $N$(natural numbers),
variables $\alp,\bet,\ldots$ over urelements in $ON$(ordinals), and
variables $x,y,z$ over urelements in either sort and sets.

 \benu
 \item
 $N(0^{N})$, $\fal x(N(x)\lrarw N(S(x)))$ and 
 for each primitive recursive relation $R$
 $\fal x_{1},\ldots,x_{n}(R(x_{1},\ldots,x_{n})\to \bigwedge_{i}N(x_{i}))$.
 \item
 $ON(0^{ON})$, $ON(\ome_{1})$,
 $\fal\alp,\bet(ON(\alp+\bet)\land ON(\alp\cdot\bet)\land ON(\ome^{\alp}))$,
 and 
 $\fal x_{1},\ldots,x_{n}(R(x_{1},\ldots,x_{n})\to \bigwedge_{i}ON(x_{i}))$
for each relation $R$ in $\mathcal{L}_{2}$.
 \item
 $\fal x,y(x\in y\to Set(y))$,
 $Set(\emptyset)$, $\fal x,y(Set(x)\to Set(J(x,y))$ and $\fal x,y(Set(J(x,y))\to Set(x))$.
 \eenu

\item(Arithmetic axioms)
Axioms in {\sf PA} for $0^{N},S$ and primitive recursive relations and the complete induction schema
\[
F(0^{N})\land\fal n(F(n)\to F(S(n)))\to\fal n\, F(n)
\]
for each formula $F$ in the language $\mathcal{L}(\Natural,\in)$.

\item(Ordinal-theoretic axioms)
Axioms in $T^{ord}(\ome_{1})$ for $<,0^{ON},+,\cdot,\lam x.\ome^{x}, R^{\cala}$,
 (\ref{eq:Z1ord}), (\ref{eq:Z2ord}), (\ref{eq:Z4ord}), (\ref{eq:Z5ord}) and (\ref{eq:Z6ord}) for 
 $P,P_{\rho_{0}}$, and
 the transfinite induction schema
 \[
 \fal \alp(\fal\bet<\alp\,F(\bet)\to F(\alp))\to\fal\alp F(\alp)
 \]
 for each formula $F$ in the language $\mathcal{L}(\Natural,\in)$.
 
\item(Set-theoretic axioms)
Extensionality $\fal x,y(Set(x)\land Set(y)\land \fal z(z\in x\lrarw z\in y)\to x=y)$,
the defining axiom for $J$, $\fal x,y,z(Set(x)\to (z\in J(x,y)\lrarw (z\in x\lor z=y)))$,
and 
\[
F(\emptyset) \land \fal x,y(Set(x)\land F(x) \to F(J(x,y)))\to \fal x(Set(x)\to F(x))
\]
 for each formula $F$ in the language $\mathcal{L}(\Natural,\in)$.
 
\eenu

\section{Ordinals for $\ome_{1}$}\label{sect:ordinals}

Let $Ord^{\varepsilon}$ and $<^{\varepsilon}$ be $\Delta$-predicates on the universe of sets such that
for any transitive and wellfounded model $V$ of the Kripke-Platek set theory with the axiom of the infinity,
$<^{\varepsilon}$ is a well ordering of type $\varepsilon_{\rho_{0}+1}$ on $Ord^{\varepsilon}$
for the order type $\rho_{0}$ of the class $Ord$ in $V$.
For natural numbers $n$,
$\ome_{n}(\rho_{0}+1)\in Ord^{\veps}$ is defined recursively by
$\ome_{0}(\rho_{0}+1)=\rho_{0}+1$ and $\ome_{n+1}(\rho_{0}+1)=\ome^{\ome_{n}(\rho_{0}+1)}$.
The $\Del$-ordering $<^{\veps}$ is seen to be a canonical ordering as stated in
the following Proposition \ref{prp:canonical}.

\bprp\label{prp:canonical}
\benu
\item\label{prp:canonical0}
$\KP\ome$ proves the fact that $<^{\veps}$ is a linear ordering.

\item\label{prp:canonical2}
For any formula $\vphi$
and each $n<\ome$,
\beqn\label{eq:trindveps}
\KP\ome\vdash\fal x\in Ord^{\veps}(\fal y<^{\veps}x\,\vphi(y)\to\vphi(x)) \to 
\fal x<^{\veps}\ome_{n}(\rho_{0}+1)\vphi(x)
\eeqn
\eenu
\eprp

$T_{1}^{+}=T_{1}+TI(\veps_{\rho_{0}+1})$ denotes the theory obtained from $T_{1}$ by adding the transfinite induction schema
along the ordering $<^{\veps}$.
In this section we work in the stronger theory $T_{1}^{+}$ otherwise stated.

For simplicity let us identify the code $x\in Ord^{\veps}$ with
the `ordinal' coded by $x$,
and $<^{\veps}$ is denoted by $<$ when no confusion likely occurs.
Note that the ordinal $\rho_{0}$ 
is the order type of the class of ordinals in the intended model $L_{\rho_{0}}$ of $T_{1}$.
Define simultaneously 
 the classes $\calh_{\alp}(X)\subset \veps_{\rho_{0}+1}$
and the ordinals $\Psi_{\ome_{1}} (\alp)$ and $\Psi_{\rho_{0}}(\alp)$ 
for $\alp<^{\veps}\veps_{\rho_{0}+1}$ and sets $X\subset \veps_{\rho_{0}+1}$ as follows.
We see that $\calh_{\alp}(X)$ and $\Psi_{\kap} (\alp)\, (\kap\in\{\ome_{1},\rho_{0}\})$ are (first-order) definable as a fixed point in $T_{1}$.

Recall that $\mbox{Hull}(X)\subset L_{\rho_{0}}$ and
$F_{X}:\mbox{Hull}(X)\lrarw L_{\gam}$ 
for $X\subset L_{\rho_{0}}$ and a $\gam=F_{X}(\rho_{0})\leq\rho_{0}$.

\bdf\label{df:Cpsiregularsm}

$\calh_{\alp}(X)$ {\rm is defined recursively as follows.}

\benu
\item
$\{0,\ome_{1},\rho_{0}\}\cup X\subset\calh_{\alp}(X)$.

\item
 $x, y \in \calh_{\alp}(X)\Rarw x+ y,\ome^{x}\in \calh_{\alp}(X)$.

\item
$\gam\in \calh_{\alp}(X)\cap\alp
\Rarw 
\Psi_{\rho_{0}}(\gam)\in\calh_{\alp}(X)
$.

\item
$\gam\in \calh_{\alp}(X)\cap\alp
\Rarw 
x=\Psi_{\ome_{1}}(\gam)\in\calh_{\alp}(X) \spand F_{x\cup\{\ome_{1}\}}(\rho_{0})\in\calh_{\alp}(X)$.

\item

{\rm Let $A(x;y_{1},\ldots,y_{n})$ be a $\Del_{0}$-formula in the language $\{\in\}$.
For $\{\alp_{1},\ldots,\alp_{n}\}\subset\calh_{\alp}(X)$,
$\mu x.\, A(x;\alp_{1},\ldots,\alp_{n})\in\calh_{\alp}(X)$,
where $\mu x.A(x;\alp_{1},\ldots,\alp_{n})=\bet$ for 
 the least ordinal $\bet$ such that $A(\bet;\alp_{1},\ldots,\alp_{n})$
if such an ordinal exists.
Otherwise $\mu x. A(x;\alp_{1},\ldots,\alp_{n})=0$.
}

\eenu

{\rm For} $\kap\in\{\ome_{1},\rho_{0}\}$ {\rm and} $\alp<\veps_{\rho_{0}+1}$
\[
\Psi_{\kap}(\alp):=
\min\{\bet\leq\kap :  \calh_{\alp}(\bet)\cap \kap \subset \bet \}
.\]
\edf

The ordinal $\Psi_{\kap}(\alp)$ is well defined and $\Psi_{\kap}(\alp)\leq \kap$
for $\kap\in\{\ome_{1},\rho_{0}\}$.

\bprp\label{prp:definability}
Both of 
$x=\calh_{\alp}(X)$ and
$y=\Psi_{\kap}(\alp)\,(\kap\in \{\ome_{1},\rho_{0}\})$
are $\Sig_{2}$-predicates. 
\eprp

\blem\label{lem:lowerbndreg}
\benu
\item
For each $n<\ome$,
$T_{1}\vdash \fal\alp<\ome_{n+1}(\rho_{0}+1)
\fal\kap\in\{\ome_{1},\rho_{0}\}\exi x<\kap[x=\Psi_{\kap}(\alp)]$.
\item
$T_{1}^{+}\vdash \fal\alp<\veps_{\rho_{0}+1}
\fal\kap\in\{\ome_{1},\rho_{0}\}\exi x<\kap[x=\Psi_{\kap}(\alp)]$.
\eenu
\elem

\section{Hydras}\label{sect:hydra}

In this section we work in the stronger theory $T_{1}^{+}(\Natural)=T_{1}(\Natural)+TI(\veps_{\rho_{0}+1})$ 
otherwise stated.
The sets of hydras $a$ and their responses $a[z]$ to Hercules' chop are introduced.
From these a hydra battle is defined
and the main theorem \ref{th:hydra} is stated.
It turns out that the battle is well-defined for each hydra in the theory $T_{1}$, cf.\,Proposition \ref{prp:H0}.\ref{prp:H0.3}.
%A \textit{move} in the hydra battles yields finitely many successors.

Each \textit{hydra} is a term over
symbols 
\[
\{0,+,\cdot,\ome,\oplus,D_{0},D_{1},D_{2}, F\}\cup
\{\times,\otimes\}\cup\calf_{\mu}
\]
where $0$
is a constant, 
each of $\ome,D_{0},D_{1},D_{2},F$ is a unary function symbol, 
$+$ a function symbol for
branching and $\cdot,\times,\otimes$ binary function symbols.
$\oplus$ is a punctuation mark.
In a hydra $c\oplus b$, $c$ is a `stock' of hydras.
The response of hydras to Hercules' chop may depend on the current stock.
The stock is kept until the hydra $b$ becomes $0$, cf.\, {\bf (sd.3)} in Definition \ref{df:dom}, 
and even enlarged when the battle goes, cf.\,{\bf (sd.5.3)}.
$f_{A}(x_{1},\ldots,x_{n})$ 
in $\calf_{\mu}$ is an $n$-ary function symbol for $\Del_{0}$-formula 
$A(x;x_{1},\ldots,x_{n})$ in the language $\mathcal{L}_{1}$.

For $a\neq 0$,
$D_{0}(a),D_{1}(a)$ denote collapsing functions $\Psi_{\ome_{1}}(a),\Psi_{\rho_{0}}(a)$, resp. 
defined in subsection \ref{sect:ordinals},
while $1:=D_{0}(0)$, $\ome_{1}:=D_{1}(0)$, $\rho_{0}:=D_{2}(0)$ and $D_{2}(a)$ denotes $\ome^{a}$ when $a\neq 0$.
$\calf_{\mu}$ is the set of $\mu$-operators $f_{A}$ for $\Del_{0}$-formulas $A$ on $\mathcal{L}_{1}$-structure $\rho_{0}$:
\beqnarrs
&& f_{A}(x_{1},\ldots,x_{n})  =  \mu x. A(x;x_{1},\ldots,x_{n}) 
\\
& = & \left\{
\begin{array}{ll}
\min\{d<\rho_{0} : A(d;x_{1},\ldots,x_{n})\}  & \mbox{{\rm if }} \rho_{0}\models\exi x\, A(x; x_{1},\ldots,x_{n})
\\
0 & \mbox{{\rm otherwise}}
\end{array}
\right.
\eeqnarrs

Let $\calf_{0}=\{f_{A}\}_{A}\subset\calf_{\mu}$ be a finite set of function symbols.
In the following Definition \ref{df:hydradf}, 
the set $H(\calf_{0})$ of hydras over $\calf_{0}$ and the set $Tm(\calf_{0})$ of terms over 
function symbols in $\{+,\cdot,\lam x. \ome^{x},F\}\cup\calf_{0}$ are defined simultaneously.
Each hydra and term is a finitary object, and can be identified with an integer.
It is clear that both of these sets are computable subsets of integers.

\bdf\label{df:hydradf}
{\rm (Simultaneous inductive definition of $H(\calf_{0})$ and $Tm(\calf_{0})$.)}
\benu

\item
{\rm 
$\{0\}\cup\{D_{i}(0) : i=0,1\}\subset H(\calf_{0})\cap Tm(\calf_{0})$ and $D_{2}(0)\in H(\calf_{0})$.}

\item
{\rm
$0\not\in \{a_{0},\ldots,a_{n}\}\subset H(\calf_{0})[\subset Tm(\calf_{0})]\, (n>0)\Rarw
 (a_{0}+\cdots+a_{n})\in H(\calf_{0})[\in Tm(\calf_{0})]$, resp.}

\item
{\rm $0\neq n<\ome, 0\neq t\in Tm(\calf_{0})\cup\{D_{2}(0)\}  \Rarw 
n\times t\in H(\calf_{0})$, where $n=\underbrace{1+\cdots+1}_{n'\footnotesize{s}\, 1}$ with $1:=D_{0}(0)$.}
%\in H_{0}(\calf_{0})$.}

\item
{\rm $0\neq n<\ome \Rarw 
n\otimes \ome\in H(\calf_{0})$.
%, where $\ome=\ome^{1}$.
}

\item
$0\neq a\in H(\calf_{0}) \Rarw D_{2}(a)\in H(\calf_{0})$.

\item
{\rm Let $c$ be a finite list $(c_{1},\ldots,c_{n})\, (n\geq 0)$ of hydras $c_{k}$ of the form
$D_{i_{k}}(d_{k}\oplus e_{k})$ in $H(\calf_{0})$, and $a\in H(\calf_{0})$ with $a\neq 0$ when 
$c=\emptyset$(empty list).
Then 
$D_{i}(c\oplus a)\in H(\calf_{0})\cap Tm(\calf_{0})$ for $i=0,1$, and
 $F(c\oplus a)\in Tm(\calf_{0})$.
When $c=\emptyset$(empty list), $\emptyset\oplus a$ denotes $a$.}

\item
$\{s,t\}\subset Tm(\calf_{0}) \Rarw \{s\cdot t,\ome^{t}\}\subset Tm(\calf_{0})$.

\item
{\rm For $f_{A}\in\calf_{0}$, if
$\{t_{1},\ldots,t_{n}\}\subset Tm(\calf_{0})$, then $f_{A}(t_{1},\ldots,t_{n})\in Tm(\calf_{0})$.}
\eenu
 \edf 
 
Terms are generated from `constants' $0, D_{i}(c\oplus a)\, (i=0,1)$ 
by function symbols $+,F,\cdot, \lam x.\ome^{x}$ and
$f_{A}$,
while hydras are generated from `constants' $0$ and $n\times t$ by function symbols $+$ and $D_{i}\,(i=0,1,2)$.
 
For hydras $D_{v}(c\oplus a)$, its \textit{local stock} is defined to be $stk(D_{v}(c\oplus a)):=c$.

The \textit{value} $v(a)<\veps_{\rho_{0}+1}$
of hydras and terms $a\in H(\calf_{0})\cup Tm(\calf_{0})$ is defined.

\bdf\label{df:hydradfvalue}
\benu

\item
$v(0)=0$, $v(D_{0}(0))=1$, $v(D_{1}(0))=\ome_{1}$, $v(D_{2}(0))=\rho_{0}$.

\item
{\rm
$v(a_{0}+\cdots+a_{n})=v(a_{0})\#\cdots\# v(a_{n})$ for the natural sum $\#$ on ordinals.}

\item
{\rm
$v(n\times t)=n\cdot v(t)$, $v(n\otimes \ome)=\ome$, 
$v(s\cdot t)=v(s)\cdot v(t)$, and
$v(\ome^{t})=\ome^{v(t)}$.}

\item
$v(D_{2}(a))=\ome^{v(a)}$.

\item

{\rm 
$v(D_{0}(c\oplus a))=\Psi_{\ome_{1}}(v(c\oplus a))$, $v(D_{1}(c\oplus a))=\Psi_{\rho_{0}}(v(c\oplus a))$ with
 $v(c\oplus a)=v(c_{1})\#\cdots \# v(c_{n})\# \ome^{v(a)}$ for the list $c=(c_{1},\ldots,c_{n})$, and
  $v(F(c\oplus a))=F_{x\cup\{\ome_{1}\}}(\rho_{0})$ for $x=v(D_{0}(c\oplus a))$.
}

\item
$v(f_{A}(t_{1},\ldots,t_{n})) = \mu x. A(x;v(t_{1}),\ldots,v(t_{n}))$.

\eenu
 \edf

Subsets $H_{i}(\calf_{0}), \cald_{i}(\calf_{0})\,(i=0,1)$ of $H(\calf_{0})\cap Tm(\calf_{0})$ are defined by
$H_{0}(\calf_{0})=\{a\in H(\calf_{0})\cap Tm(\calf_{0}): v(a)<\ome_{1}\}$ {\rm and}
$H_{1}(\calf_{0})=\{a\in H(\calf_{0})\cap Tm(\calf_{0}): v(a)<\rho_{0}\}$.
%$\cald_{i}(\calf_{0})\subset H_{i}(\calf_{0})$
%{\rm for} $i=0,1$.
Note that
$v(a)<\rho_{0}$ for any $a\in Tm(\calf_{0})$.

We see that $v(a)=\alp\,(a\in H(\calf_{0})\cap Tm(\calf_{0}), \alp\in Ord^{\veps})$
is a $\Del_{2}$-predicate from Proposition \ref{prp:definability}.

Let us identify the hydras and terms $a$ with the ordinals $v(a)$, and let
\beqnarrs
a<b:\Lrarw v(a)<v(b) & \Lrarw & \exi\alp,\bet\in Ord^{\veps}[v(a)=\alp<^{\veps}\bet=v(b)]
\\
& \Lrarw & \fal\alp,\bet\in Ord^{\veps}[v(a)=\alp\to v(b)=\bet\to\alp<^{\veps}\bet]
\eeqnarrs
$a<b$ is again a $\Del_{2}$-predicate on integers $a,b$, where
$\Del_{2}$ denotes a class in the Levy hierarchy.
Let $v(D_{1}(c\oplus))=\Psi_{\rho_{0}}(v(c))$ and $v(D_{0}(c\oplus))=\Psi_{\ome_{1}}(v(c))$
for lists $c=(c_{1},\ldots,c_{n})$ of hydras and $v(c)=v(c_{1})\#\cdots \# v(c_{n})$.

\bprp\label{prp:valuelinear}
$T_{1}(\Natural)$ 
proves the following facts for each $n\in\Natural$ and hydras $a,b$ with $v(a),v(b)<\ome_{n}(\rho_{0}+1)$:
$a<b:\Lrarw v(a)<v(b)$ is a linear ordering on quotient sets of hydras by the equivalence relation 
$a\simeq b:\Lrarw v(a)=v(b)$.
\eprp
\bprf
This is seen from Lemma \ref{lem:lowerbndreg}.
\eprf

\bdf\label{df:hydradfsize}
{\rm
The \textit{size} $|t|\in\Natural$ of terms $t\in Tm(\calf_{0})$ is defined to be the total number of occurrences
of symbols $0,+,\cdot,\ome,\oplus,D_{0},D_{1},D_{2}, F,\times,\otimes$ and $f_{A}\in\calf_{\mu}$ in $t$.
}
\edf

\bprp\label{prp:sizefinite}
For each finite set $\calf_{0}$ of function symbols $f_{A}$,
there exists a constant $c$ such that for any $k$,
the number of terms in size$\leq k$ is bounded by $c^{k}$,
$\#\{t\in Tm(\calf_{0}): |t|\leq k\}\leq c^{k}$.
\eprp

\bdf
{\rm For terms $t,s\in Tm(\calf_{0})\cup\{D_{2}(0)\}$, $i=0,1$, and lists $c\subset H(\calf_{0})$ let}
\beqnarrs
s<_{i}c & :\Lrarw & v(s)\in H_{v(c)}(v(D_{i}(c)))
%defined through G_{a}?
%v(s)\in\bigcap\{H_{v(c)}(v(D_{i}(c\oplus d))): d\in H(\calf_{0})\}
\\
multi_{t,2}(\calf_{0})  & := & \{s: s\in Tm(\calf_{0}), s<t\}
\\
multi_{t,1}(c;\calf_{0}) & := & \{s\in multi_{t,2}(\calf_{0}):  s<_{1}c\}
%v(s)\in\calh_{v(c)}(v(D_{1}(c\oplus)))\}
\eeqnarrs

\edf

\bprp\label{prp:increase}
Assume $s<_{i}c$ and $v(c)\in H_{v(c)}(v(D_{i}(c)))$.
Then
$v(s)\in\bigcap\{H_{v(c\oplus d)}(v(D_{i}(c\oplus d))): d\in H(\calf_{0})\}$.
\eprp
\bprf
Let $\kap_{0}=\ome_{1}$ and $\kap_{1}=\rho_{0}$.
Suppose $v(s),v(c)\in H_{v(c)}(v(D_{i}(c)))$.
Then $H_{v(c)}(v(D_{i}(c\oplus d)))\cap\kap_{i}\subset H_{v(c\oplus d)}(v(D_{i}(c\oplus d)))\cap\kap_{i}\subset (D_{i}(c\oplus d))$.
Hence $v(D_{i}(c))\leq v(D_{i}(c\oplus d))$, and we obtain $H_{v(c)}(v(D_{i}(c)))\subset H_{v(c\oplus d)}(v(D_{i}(c\oplus d)))$.
\eprf

\bdf\label{df:dom}{\rm (Stepping-down or Hydra's response)}\\
{\rm Let $\calf_{0}\subset\calf_{\mu}$ be a finite set of function symbols, and $a\in H(\calf_{0})$ a hydra.
Its \textit{domain}
$dom(a)$ and hydras $a[z]\subset H(\calf_{0})$ are defined for $z\in dom(a)$.

$dom(a)$ is one of sets $\emptyset(=0),1(=\{0\}),\Natural, H_{i}(\calf_{0})\, (i=0,1)$
or
one of sets
\\
$multi_{t,2}(\calf_{0}),multi_{t,1}(c_{1};\calf_{0})$
for a term $t\in Tm(\calf_{0})\cup\{D_{2}(0)\}$ with $v(t)\neq 0$, and some lists $c_{1}\subset H(\calf_{0})$.}

 \bdes
  \item[{(sd.0)}]
 $dom(0)=\emptyset$. 
 
 \item [{(sd.1)}]
 $dom(1)=1$,
 $1[0]=\{0\}$ {\rm where} $1=D_{0}(0)$.
 
 \item [{(sd.2)}]
 $dom(D_{i+1}(c\oplus 0))= H_{i}(\calf_{0})$;
 $(D_{i+1}(c\oplus 0))[z]=\{z\}$ {\rm for} $i=0,1$.

 \item [{(sd.3)}] 
 {\rm Let} $m>0$ {\rm and} $0\neq t\in Tm(\calf_{0})\cup\{D_{2}(0)\}$.
$dom(m\times t) = multi_{t,2}(\calf_{0})$; 
 $(m\times t)[s] = \{(m\times s)+(m-1)\}$,
 {\rm where}
$(m\times 0)+b:=b$.

 \item [{(sd.4)}] 
 {\rm Let} $m>0$.
$dom(m\otimes \ome) = \Natural$; 
 $(m\otimes \ome)[n] = \{m\cdot (n+1))\}$,
 {\rm where}
$m\cdot (n+1)=1+\cdots+1$ {\rm with $m(n+1)$ times $1$'s.}

 \item [{(sd.5)}] 
  {\rm Let $a=D_{i}(c\oplus b)$, where $b$ is a non-zero hydra.}

   \item[{(sd.5.1)}] 
    {\rm If} $b=b_{0}+1$, {\rm then} $dom(a)=\Natural$; 
    $a[n]=\{(D_{i}(c\oplus b_0))\cdot 2\}$ {\rm for} $n\in\Natural$, where
   $(D_{i}(c\oplus b_0))\cdot 2:=D_{i}(c\oplus b_0)+D_{i}(c\oplus b_0)$.
   
    \item[{(sd.5.2)}] 
     {\rm If either
     $dom(b)\in\{1,\Natural, H_{j}(\calf_{0}): j< i\}$, or $dom(b)\in\{multi_{t,1}(c;\calf_{0}):t\in Tm(\calf_{0}), c\subset H(\calf_{0})\}$ and $i=1$, 
     then
 $dom(a)=dom(b)$; $a[z]=\{D_{i}(c\oplus d):d\in b[z]\}$ for $z\in dom(b)$.}

     \item[{(sd.5.3)}] 
      {\rm If} $dom(b)\in\{ H_{j}(\calf_{0}): j\geq i\}$,
      {\rm then} $dom(a)=1$.
       {\rm Let} $\ell=D_{i}(c\oplus (b[1]))$
    {\rm and} $r=D_{i}((c+D_{2}(b[1])+1)\oplus (b[1]))$ {\rm with} $stk(r)=c+D_{2}(b[1])+1$,
     $a[0]:=\{\ell+r\}$ {\rm if} $i=1$. $a[0]:=\{r\}$ {\rm if} $i=0$,
     {\rm where
     $c+D_{2}(b[1])+1:=c*\la D_{2}(b[1]),1\ra$, a concatenated list.}
           
     \item[{(sd.5.4)}] 
     {\rm If $dom(b)=multi_{t,2}(\calf_{0})$ for a term $t$ and $i=1$, then} 
     $(D_{1}(c\oplus b))[s]=\{D_{1}(c\oplus (b[s]))\}$ {\rm for
     $s\in dom(D_{1}(c\oplus b)):=multi_{t,1}(c;\calf_{0})$.}
     
       \item[{(sd.5.5)}] 
        {\rm If $dom(b)\in\{multi_{t,2}(\calf_{0}), multi_{t,1}(c_{1};\calf_{0}):t\in Tm(\calf_{0}), c_{1}\subset H(\calf_{0})\}$
        and $i=0$,
        then $dom(D_{0}(c\oplus b))=\Natural$. 
        Let $multi_{n}$ denote the set
        \[
         \{s\in multi_{t,2}(\calf_{0}) : |s|\leq 2^{2^{n}}, s<_{1}c_{1}, s<_{0}c\}
         %v(s)\in \calh_{v(c_{1})}(v(D_{1}(c_{1}\oplus )))\cap\calh_{v(c)}(v(D_{0}(c\oplus)))\}
         \]
         {\rm if }
         $dom(b)=multi_{t,1}(c_{1};\calf_{0})$.
         {\rm Otherwise}
        \[
        multi_{n}=
         \{s\in multi_{t,2}(\calf_{0}): |s|\leq 2^{2^{n}}, s<_{0}c\}
         %v(s)\in \calh_{v(c)}(v(D_{0}(c\oplus)))\}.
        \]
        Then
    $a[n]:=\{D_{0}(c\oplus (b[s])): s\in multi_{n}\}$.}

 \item [{(sd.6)}]
 {\rm Let $a= (a_{0}+\cdots+a_{k}) \, (k>0)$, where each
 $a_{i}$ is a non-zero hydra. 
 Then
 $dom(a)=dom(a_{k})$;
  $a[z]=\{a_{0}+\cdots+a_{k-1}+b: b\in (a_{k}[z])\}$.}

 \edes

\edf
When $a[z]=\{b\}$ is a singleton, we write $b=a[z]$.

Note that the case $dom(a)=\Natural$ occurs essentially only in the cases {\bf (sd.4)} and {\bf (sd.5.5)}.
The latter case is close to the definition of the fundamental sequences in \cite{BCW} based on norm bounding.
Moreover $a[n]$ is not a singleton only in this case.
%We have $multi_{s,2}(\calf_{0})=multi_{s',2}(\calf_{0})$ when $v(s)=v(s')$.
%This means that it does not matter which maximal term $s_{n}$ is chosen from the finite set 
%$multi_{n}$ in the Hydra battles in the case {\bf (sd.5.5)}.

The term $t$ of $m\times t$ in the case {\bf (sd.3)} is regarded as the rightmost head of the hydra $m\times t$.
When Hercules chops off the head, the hydra chooses a term $s$ from an
 \textit{infinite} set $multi_{t,2}(\calf_{0})$,
and $m\times t$ turns to $m\times s+(m-1)$.
On the other side the $m\times t$ is the rightmost head in a hydra 
$a=D_{0}(c_{0}\oplus b_{0}[\cdots D_{1}(c_{1}\oplus b_{1}[\cdots m\times t])])$,
the hydra $a$ builds a term $s$ from constants $0,1,\ome_{1}$
and a finite number of function symbols in
$\{+,\cdot, \lam x.\ome^{x},D_{0},D_{1},D_{2},F\}\cup\{\times\}\cup\calf_{0}$.
Moreover the hydra obeys the restrictions 
$s<_{1}c_{1}$, $s<_{0}c$,
%$v(s)\in \calh_{v(c_{1})}(v(D_{1}(c_{1}\oplus )))\cap\calh_{v(c)}(v(D_{0}(c\oplus)))$ 
and $|s|\leq 2^{2^{n}}$.
In particular the hydra has to choose a term $s$ from the 
\textit{finite} set 
$multi_{n}$, cf.\,Proposition \ref{prp:sizefinite}, and the hydra turns to
$D_{0}(c_{0}\oplus b_{0}[\cdots D_{1}(c_{1}\oplus b_{1}[\cdots (m\times s+(m-1))])])\in a[n]$, cf.\,{\bf (sd.5.5)}.
This is the only case when the response of hydras may depend on its stock.

\bprp\label{prp:H0}
$T_{1}(\Natural)$ proves the following facts for {\rm each} $a\in H(\calf_{0})$.
\benu
\item\label{prp:H0.1}
If $a\in H_{0}(\calf_{0})$ and $b\in a[n]$, then
$dom(a)\in\{0,1,\Natural\}$, 
$b\in H_{0}(\calf_{0})$ and $|b|\leq \max\{2|a|+2^{2^{n}},3|a|, |a|(n+1)\}$ for $n\in dom(a)$.

\item\label{prp:H0.3}
For any $z\in dom(a)$,
$\exi ! x\subset H(\calf_{0})(x=a[z])$.

\eenu
\eprp
\bprf
\ref{prp:H0}.\ref{prp:H0.1}.
$dom(a)\in\{0,1,\Natural\}$ is seen easily for $a\in H_{0}(\calf_{0})$.
Let $b=m\times t$ and $b[s]=m\times s+(m-1)$ for an $s$ such that $|s|\leq 2^{2^{n}}$.
Then $|b[s]|\leq 4m-2+2^{2^{n}}\leq 2|b|+2^{2^{n}}$ with $|b|=2m+|t|$.
From this we see that $|b|\leq 2|a|+2^{2^{n}}$ in the case {\bf (sd.5.5)}.

Next consider the case {\bf (sd.4)}.
$a=m\otimes \ome$ and $b=m\cdot(n+1)$.
Then $|b|=2m(n+1)-1\leq(2m+1)(n+1)=|a|(n+1)$.

Finally consider the case {\bf (sd.5.3)}, $a=D_{1}(c_{1}\oplus b)$ with $dom(b)= H_{1}(\calf_{0})$.
Then $a[n]=\{\ell+r\}$ for $\ell=D_{1}(c_{1}\oplus b[1])$ and $r=D_{1}((c_{1}+D_{2}(b[1])+1)\oplus b[1])$.
Hence
$|\ell+r|=|\ell|+|r|+1\leq |a|+2|a|$.
\\
\ref{prp:H0}.\ref{prp:H0.3}.
This is seen from Proposition \ref{prp:valuelinear}.
\eprf
\\

When $dom(a)=1$, let $a[n]:=a[0]$ for any $n\in\Natural$, and $0[n]:=0$, where $dom(a)=0$ iff $a=0$.
Also $a[n][m]:=\bigcup\{b[m]:b\in a[n]\}\cup\{0\}$ for $dom(a)\in\{0,1,\Natural\}$.
\bdf\label{df:Hardy}
{\rm
For a finite set $\calf_{0}$ of function symbols and a hydra $a\in H_{0}(\calf_{0})$,
$h_{a}^{\calf_{0}}:\mathbb{N}\to\mathbb{N}$ denotes a (possibly partial) 
number-theoretic function defined as follows.
Let $h_{0}^{\calf_{0}}(n):=n$, and for $a\neq 0$,
$h_{a}^{\calf_{0}}(n):\simeq h_{a[n]}^{\calf_{0}}(n+1)$, i.e.,
\[
h_{a}^{\calf_{0}}(n):\simeq \min\{k>n: a[n][n+1]\cdots[k-1]=\{0\}\}
\]
where the stepping-down $a[n]$ is determined from $\calf_{0}$.
}
\edf

\brem
{\rm
Let us consider some restricted hydras.
First consider $(\times,\otimes)$-free hydras $a$.
Then the cases {\bf (sd.3)}, {\bf (sd.5.4)} and {\bf (sd.5.5)} nor $dom(a)=\Natural$ never occur, and 
neither the stock nor the set $\calf_{0}$ plays a r\^{o}le
in the stepping-down.
Hence $a[n]$ does not depend on $n$, and
$h_{a}$ is seen to be an $\mathcal{E}^{4}$-function.

Next consider hydras $H(\emptyset)$ over $\calf_{0}=\emptyset$ without the Mostowski collapsing function $F$.
We see then that the relation $s<t\Lrarw v(s)<v(t)$ on terms $s,t$ is computable, and
${\sf ID}_{2}$ proves that $h_{a}^{\emptyset}$ is defined for \textit{each} hydra $a$ with $dom(a)\in\{0,1,\Natural\}$,
while the fact that `$h_{a}^{\emptyset}$ is defined for \textit{every} hydra $a$ with $dom(a)\in\{0,1,\Natural\}$'
is independent over ${\sf ID}_{2}$, cf.
%\,Remark \ref{rem:rec} and 
section \ref{sect:consisprf}.

In general $dom(a)$ and $a[z]$ for $z\in dom(a)$ are far from being computable in the case {\bf (sd.5.5)}.
}
\erem

Initial hydras $a$ in the following definition are assigned to given proofs in $T_{1}(\Natural)$, 
cf.\,section \ref{sect:finiteproof}, and
are seen to enjoy $a[n]\cdots[n+m-1][n+m]<a[n]\cdots[n+m-1]$ for any $n,m\in\mathbb{N}$, 
which means that $v(c)<c(b)$ for any $b\in a[n]\cdots[n+m-1]$ and any $c\in b[n+m]$.
cf.\,Proposition \ref{prp:WB} and Lemma \ref{prp:H0.2}.
Hence the function $h_{a}^{\calf_{0}}$ is seen to be total for initial hydras $a$.

\bdf\label{df:WB}{\rm (Initial hydras)}
\\
{\rm $\mathcal{I}$ denotes the set of hydras generated from $D_{i}(0)\,(i=0,1,2)$, 
$n\times D_{2}(0), n\otimes\ome\,(n>0)$ by applying $+$.

Then each hydra 
$a=D_{0}(D_{2}^{(k+2)}(D_{2}(0)+1)\oplus D_{2}^{(k)}(D_{1}(\emptyset\oplus D_{2}^{(k)}(b))))$
%+k$ 
for a $b\in\mathcal{I}$ and a $k<\ome$ is said to be an \textit{initial hydra},
where 
$D_{2}^{(k)}(b)=D_{2}(\cdots D_{2}(D_{2}(b))\cdots)$ with $k$'s $D_{2}$.
}
\edf
%$k=0$ $stk(a)=D_{2}(0)+1=(D_{2}(0),D_{0}(0))$

%\bdf\label{df:hydrabattle}{\rm (Hydra battles)}\\
%{\rm Let $\calf_{0}\subset\calf_{\mu}$ be a finite set of function symbols.
%Let $a_{0}$ be a given initial hydra.
%Define hydras $a_{n}\in H_{0}(\calf_{0})$ over $\calf_{0}$ recursively as follows.}
%\benu
%\item\label{df:hydrabattle0}
%{\rm If $a_{n}=0$, then $a_{n+1}=0$.

%In what follows assume $a_{n}\neq 0$.}

%\item\label{df:hydrabattle1}
%$a_{n+1}=a_{n}[n]$ {\rm where the RHS is determined from $\calf_{0}$.}

%\eenu
%\edf

Let $\mbox{(H)}_{\ome_{1}}$ denote the statement saying that
%\begin{quote}
for any initial hydra $a_{0}$, any finite set $\calf_{0}\subset\calf_{\mu}$ of function symbols and any $n\in\Natural$, 
there exists an $m\in\Natural$ such that $h_{a_{0}}^{\calf_{0}}(n)\simeq m$.
%there exists an $n<\ome$ such that
%$a_{n}=0$.
%\end{quote}
This means that an initial hydra $a_{0}$ first chooses a finite set $\calf_{0}$ of function symbols
and an $n\in\Natural$ arbitrarily.
The hydra responds to Hercules' chop to its right-most head using $\calf_{0}$ and the number $n+m$
in the $m$-th round according to Definition \ref{df:dom}.
Then the hydra eventually die out in the battle, no matter which term $s$ is chosen in the case {\bf (sd.5.5)}.

For Theorems \ref{th:hydra}.\ref{th:hydra-1} and \ref{th:hydra}.\ref{th:hydra0},
we show that if hydras (ordinal terms) $a$ and $z$ enjoy a condition in the following
Definition \ref{df:wellbehaved},
then so does $b$ and $v(b)<v(a)$ for any $b\in a[z]$.

%The set $E_{i}(a)$ is defined as in \cite{BuchholzAPAL86}.
To prove Lemma \ref{prp:H0.2} below, it is convenient for us to split the set $G_{i}(a)$ in \cite{BuchholzAPAL86}.
In the following definition $g_{i}(a)$ denotes a subset of $G_{i}(a)$.
For a multiplicative hydra $n\times t$,
$g_{i}(n\times t)=\emptyset$ and $E_{i}(n\times t)=\{n\times t\}$.
%&, $M_{i}(n\times t)=\{t\}$, and $G_{i}(n\times t)=G_{i}(t)$.

\bdf\label{df:gM}{\rm (Finite sets $g_{i}(a)$ and $E_{i}(a)$)}
\benu

\item
{\rm 
$g_{i}(0)=E_{i}(0)=g_{i}(D_{v}(0))=E_{i}(D_{v}(0))=\emptyset$ for $i=0,1$ and $v=0,1,2$.}

\item
{\rm
$g_{i}(a_{0}+\cdots+a_{n})=\bigcup\{g_{i}(a_{k}): k\leq n\}$
and $E_{i}(a_{0}+\cdots+a_{n})=\bigcup\{E_{i}(a_{k}): k\leq n\}$.}

\item
$g_{i}(n\times t)=\emptyset$.
$E_{i}(n\times t)=\{n\times t\}$
$g_{i}(D_{2}(a))=g_{i}(a)$.
$E_{i}(D_{2}(a))=E_{i}(a)$.

\item
$g_{i}(n\otimes \ome)=E_{i}(n\otimes \ome)=\emptyset$.

\item
{\rm For $v=0,1$ and $b=(b_{1},\ldots,b_{n})$, let}
  \[
  g_{i}(D_{v}(b\oplus a))=\left\{
  \begin{array}{ll}
  \{b\oplus a\}\cup g_{i}(b_{1})\cup\cdots\cup g_{i}(b_{n})\cup g_{i}(a) & \mbox{{\rm if }} i\leq v
  \\
  \emptyset &  \mbox{{\rm if }} i> v
  \end{array}
  \right.
  \]
   \[
  E_{i}(D_{v}(b\oplus a))=\left\{
  \begin{array}{ll}
    \{D_{v}(b\oplus a)\}\cup E_{i}(b_{1})\cup\cdots\cup E_{i}(b_{n})\cup E_{i}(a) &  \mbox{{\rm if }}  i\leq v
  \\
  \emptyset &  \mbox{{\rm if }}  i> v
  \end{array}
  \right.
  \]

\eenu
 \edf

%For a set $G$ of hydras and a list of hydras $c$,
%let $G<c:\Lrarw\fal a\in G(v(a)<v(c))$.

It is easy to see by induction on the size $|a|$ of hydras $a$ that 
\beqn\label{eq:DH}
E_{i}(a)\subset \calh_{c}(D_{i}(c\oplus)) \Rarw a\in\calh_{c}(D_{i}(c\oplus))
\eeqn
where
$a\in\calh_{c}(D_{i}(c\oplus)):\Lrarw v(a)\in\calh_{v(c)}(v(D_{i}(c\oplus)))$.

%\beqn\label{eq:MDH}?
%M_{i}(a)\subset\calh_{c}(D_{i}(c\oplus)) \spand G_{i}(a)<c \Rarw a\in\calh_{c}(D_{i}(c\oplus))
%\eeqn
For a set $G$ of hydras and a list of hydras $c$,
let $G<_{i}c:\Lrarw\fal a\in G(a<_{i}c)$.

\bdf\label{df:wellbehaved}
{\rm A hydra} $b$ {\rm is said to be} well-behaved {\rm if every $D$-subhydra $D_{i}(c_{i}\oplus a)\, (i=0,1)$ of $b$
enjoys the following condition:}
\beqn\label{eq:wellbehaved}
\{c_{i}\}\cup E_{i}(a)<_{i}c_{i} \spand g_{i}(a)<c_{i}
%c_{i}\oplus a\in\calh_{c_{i}}(D_{i}(c_{i}\oplus))
%G_{i}(c_{i}\oplus a)<c_{i}
\eeqn
\edf
This means that $D_{i}(c_{i}\oplus a)\, (i=0,1)$ is well-behaved if (\ref{eq:wellbehaved}) holds and
both $c_{i}$ and $a$ are well-behaved.
Each of $0$, $D_{k}(0)\,(k=0,1,2)$ and $n\times t$ is well-behaved.
$a_{0}+\cdots+a_{n}$ is well-behaved iff each $a_{i}$ is well-behaved.
If $a$ is well-behaved, then so is $D_{2}(a)$.

%Note that any multiplicative hydra $m\times t$ is well-behaved since there is no $D$-subhydra of $m\times t$.

\bprp\label{prp:WB}
Each initial hydra is well-behaved.
\eprp
\bprf
%For $b\in\mathcal{I}$, we see $G_{1}(b)=G_{0}(b)=\emptyset$.
For $c=D_{2}^{(k+2)}(D_{2}(0)+1)$ and $a=D_{2}^{(k)}(D_{1}(\emptyset\oplus D_{2}^{(k)}(b)))$
we see $\{c\}\cup E_{0}(a)<_{0}c$ from
%\calh_{c}(D_{0}(c\oplus))$ from
$v(\emptyset\oplus D_{2}^{(k)}(b))=\ome_{k+1}(v(b))<\ome_{k+2}(\rho_{0}+1)=v(D_{2}^{(k+2)}(D_{2}(0)+1))=v(c)$
since $v(b)<\rho_{0}\cdot\ome=\ome^{\rho_{0}+1}$.
Also $g_{0}(a)=\{\emptyset\oplus D_{2}^{(k)}(b)\}$.
\eprf

\blem\label{prp:H0.2}
$T_{1}(\Natural)$ proves the following facts for {\rm each} $a\in H(\calf_{0})$:
Let $a\in H(\calf_{0})$ be a well-behaved hydra.
Then
$b$ is well-behaved and $b<a$ if $b\in a[z]$ and $z\in dom(a)$ is either well-behaved or 
$dom(a)\in \{multi_{t,2}(\calf_{0}), multi_{t,1}(c;\calf_{0})\}_{c}$.
\elem
\bprf
By induction on the sizes $|a|$ of hydras $a$.
Let $a$ be a well-behaved hydra.
Consider the case when $a=D_{v}(c\oplus b)$ with $v=0,1$ and $b\neq 0$.
We have
\beqn\label{eq:H0.2}
\{c\}\cup E_{v}(b)<_{v}c \spand g_{1}(b)<c
\eeqn
{\bf (sd.5.1)}.
If $b=b_{0}+1$, then $a[n]=\{D_{v}(c\oplus b_{0})\cdot 2\}$.
We have $g_{v}(b_{0})=g_{v}(b)$, $E_{v}(b_{0})=E_{v}(b)$, and $D_{v}(c\oplus b_{0})$ is well-behaved by (\ref{eq:H0.2}).
Also $c\oplus b_{0}<c\oplus b$ and $\{c,b_{0}\}\subset\calh_{c}(a)$ by (\ref{eq:wellbehaved}).
Hence $D_{v}(c\oplus b_{0})\in\calh_{c}(a)\cap D_{v+1}(0)\subset a$.
%since $G_{v}(b_{0})=G_{v}(b_{0}+1)$.
\\
{\bf (sd.5.2)}.
If $dom(a)=dom(b)$, then $a[z]=\{D_{v}(c\oplus d):d\in b[z])\}$.
By IH $d$ is well-behaved and $d<b$.
It suffices to show that $E_{v}(d)<_{v}c$, and
$g_{1}(d)<c$.
We have $g_{1}(d)\subset g_{1}(b)$.

If either $dom(b)=\Natural$ or $dom(b)= H_{0}(\calf_{0})\spand v=1$, then $E_{v}(d)\subset E_{v}(b)$
with $E_{v}(z)=\emptyset$,
and $D_{v}(c\oplus d)$ is well-behaved.
Let $dom(b)=multi_{t,1}(c_{1};\calf_{0})\ni z$ and $v=1$. 
Then $b=b[D_{1}(c_{1}\oplus b_{0}[m\times t])]$ for some $b_{0}$ with $dom(b_{0})=multi_{t,2}(\calf_{0})$
and $m<\ome$.
We have $z<_{1}c_{1}$, and
$E_{v}(b[z])\subset E_{v}(b)\cup \{m\times z)$.
On the other hand we have $c_{1}<c_{1}\oplus b_{0}\in g_{1}(b)<c$.
% and $D_{1}(c_{1}\oplus)<D_{1}(c\oplus)$ by $c_{1}\in\calh_{c}(D_{1}(c\oplus))$.
Hence $z<_{1}c$,
%$z\in\calh_{c}(D_{v}(c\oplus))$, 
and $a[z]$ is well-behaved.

From $\{c,b[z]\}\subset\calh_{c}(D_{v}(c\oplus b))$,
we see that $c\oplus b[z]\in \calh_{c}(D_{v}(c\oplus b))$.
Hence we obtain $a[z]<a$ by $b[z]<b$. 
\\
{\bf (sd.5.3)}.
If $dom(b)= H_{i}(\calf_{0})$ with $i\geq v$, then 
$a=D_{v}(c\oplus b[D_{i+1}(0)])$ and $a[n]\in\{\ell+r,r\}$
for $n=z\in dom(a)=\ome$, $\ell=D_{v}(c\oplus b[1])$ and $r=D_{v}((c+D_{2}(b[1])+1)\oplus b[1])$.
By IH $b[1]$ is well-behaved.
We see that $g_{v}((c+D_{2}(b[1])+1)\oplus b[1])\subset g_{v}(c\oplus b)<c<1\# c$,
and
$E_{v}(b[1])\subset E_{v}(b)$.
%$M_{v}(b[1])\subset M_{v}(b)$ and $M_{v}(c+D_{2}(b[1])+1)=M_{v}(c\oplus b[1])\subset\calh_{c}(D_{v}(1+c+d))$ for any $d$.
Hence $a[n]$ is well-behaved.
It is clear that $v(c)\#\ome^{v(b[1])}\cdot 2+1=v((c+D_{2}(b[1])+1)\oplus b[1])<v(c\oplus b)=v(c)\#\ome^{v(b)}$ by $b[1]<b$, and
$a[n]<a$.
\\
{\bf (sd.5.4)}.
If $dom(b)=multi_{t,2}(\calf_{0})$ and $v=1$, then
$a[z]=D_{1}(c\oplus b[z])$ for $z\in multi_{t,1}(c;\calf_{0})$.
By IH $b[z]$ is well-behaved and $b[z]<b$.
We have $g_{1}(b[z])\subset g_{1}(b)$ and $E_{1}(b[z])\subset E_{1}(b)\cup\{n\times z:n\in\ome\}$.
From $z\in multi_{t,1}(c;\calf_{0})$ we see $z<_{1}c$,
and hence $a[z]$ is well-behaved.
Also $a[z]<a$ follows from $b[z]<b$ as in the case $dom(b)=\ome, H_{0}(\calf_{0})$.
\\
{\bf (sd.5.5)}.
Finally let $dom(b)=multi_{t,2}(\calf_{0}),multi_{t,1}(c_{1};\calf_{0})$ and $v=0$.
This means that $b=b[m\times t]$ or $b=b[D_{1}(c_{1}\oplus b_{0}[m\times t])]$ for some $b_{0}$ and $m$ such that
$dom(b_{0})=multi_{t,2}(\calf_{0})$.
Let $multi_{n}$ denote the set in Definition \ref{df:dom}.
Then $n=z\in dom(a)=\Natural$ and $D_{0}(c\oplus b[s])\in a[n]$ for a term $s\in multi_{n}$,
where $b[s]=b[(m\times t)[s]]$ or
$b[s]=b[D_{1}(c_{1}\oplus b_{0}[(m\times t)[s]])]$
with $(m\times t)[s]=m\times s+(m-1)$.

From IH we see that $b[s]$ is well-behaved,
$b_{0}[(m\times t)[s]]<b_{0}[m\times t]$ and $b[s]<b$.
We have $E_{0}(b[s])\subset E_{0}(b)\cup\{m\times s\}<_{0}c$.
by $s\in multi_{n}$.
If $dom(b)=multi_{t,2}(\calf_{0})$, then $g_{0}(c\oplus b[s])\subset g_{0}(c\oplus b)$.
If $dom(b)=multi_{t,1}(c_{1};\calf_{0})$, then
$g_{0}(c\oplus b[s])\subset g_{0}(c\oplus b)\cup\{c_{1}\oplus b_{0}[(m\times t)[s]]\}
\leq g_{0}(c\oplus b)<c$ by $c_{1}\oplus b_{0}[(m\times t)[s]]<c_{1}\oplus b_{0}[m\times t]\in g_{0}(b)<c$.
Hence $D_{0}(c\oplus b[s])$ is well-behaved.

%We have $G_{0}(c\oplus b[s_{n}])<c$.
% $M_{0}(c\oplus b[s_{n}])\subset\calh_{c}(D_{0}(c))$
%and $c\in\calh_{c}(D_{0}(c\oplus b))$.
We obtain  $c\oplus b[s]\in\calh_{c}(D_{0}(c\oplus b))$.
% by (\ref{eq:DH}).
Hence $D_{0}(c\oplus b[s])<a$ by $b[s]<b$.
\eprf
\\

Lemma \ref{prp:H0.2} yields Theorem \ref{th:hydra}.\ref{th:hydra-1}, and
Theorem \ref{th:hydra}.\ref{th:hydra0} by Lemma \ref{lem:lowerbndreg}.

\subsection{Urelements}

Theorem \ref{th:hydra}.\ref{th:hydra1} follows from the following Lemma \ref{lem:consis}.

\bdf\label{df:Sigma2*}
{\rm
A formula in the language $\mathcal{L}_{2}$ is a \textit{$\Del_{0}$-formula} if every quantifier in it is bounded.
%$\Sig_{n}$-formulas and their dual $\Pi_{n}$-formulas in $\mathcal{L}_{2}$ are defined as usual.

The set of \textit{$\Sig_{2}^{*}$-formulas} in the language $\mathcal{L}(\Natural,\in)$ is defined recursively as follows.
\benu
\item
For each $\Pi_{1}$-formula $\fal x\,A(x)\,(A\in\Del_{0})$ in $\mathcal{L}_{2}$,  
$\fal\alp\, A(\alp):\equiv(\fal \alp(ON(\alp)\to A(\alp)))$ is a $\Sig_{2}^{*}$-formula.
\item
Each arithmetic literal $R(t_{1},\ldots,t_{n})$ is a $\Sig_{2}^{*}$-formula with relation symbols $R$ for primitive recursive relation.
\item
If $A_{0}$ and $A_{1}$ are $\Sig_{2}^{*}$-formulas, then so are $A_{0}\lor A_{1}$ and $A_{0}\land A_{1}$.
\item
If $A$ is a $\Sig_{2}^{*}$-formula, then so are $\exi \alp\, A$,
$\exi n\, A$, $\exi x(Set(x)\land A)$ and $\fal n<^{N}t\, A, \fal x\in t\,A$,
where $\exi\alp\, A:\equiv(\exi \alp(ON(\alp)\land A))$,
$\exi n\, A:\equiv(\exi n(N(n)\land A))$ and similarly for $\fal n<^{N}t\, A$ and $\fal x\in t\,A$.
%\item
%If $A$ is a $\Sig_{2}^{*}$-formula, then so is $\exi x(Set(x)\land A)$.
\eenu
If $A$ is a $\Sig_{2}^{*}$-formula, then $\lnot A$ is a \textit{$\Pi_{2}^{*}$-formula}.
}
\edf

\bprp
The relation $\{(a,b,n)\in\Natural^{3}: b\in a[n], a\in H_{0}(\calf_{0})\}$ is written in a $\Pi^{*}_{2}$-predicate $R$ where
hydras are coded by natural numbers in $\Natural$.
\eprp
\bprf
Definition \ref{df:dom} of $b=a[n]$ is done by definition by cases.
Consider the case {\bf (sd.5.5)}.
$b\in a[n]$ is defined from a term $s$ in the finite set $multi_{n}$.
It suffices to show that the relation $\{(a,\alp)\in\Natural\times \veps_{\rho_{0}+1}: v(a)=\alp\}$ is a $\Sig_{2}^{*}$-relation.
Then $s\in multi_{n}$ is seen to be a $\Pi_{2}^{*}$-relation
since $v(s)<v(t)$ iff $\fal\alp,\bet(v(s)=\alp\land v(t)=\bet\to\alp<\bet)$.

Now $v(a)=\alp$ iff there exists a function on the set of subterms of $a$ to ordinals$<\veps_{\rho_{0}+1}$
enjoying the inductive clauses in Definition \ref{df:hydradfvalue}.
We see from the equivalence that $v(a)=\alp$ is $\Sig_{2}^{*}$ from
Proposition \ref{prp:definability}.
\eprf

\blem\label{lem:consis}
Each $T_{1}(\Natural)$-provably total $\Sig^{*}_{2}$-functions on integers
is dominated by a function $1+h_{a}^{\calf_{0}}$ for an initial hydra $a$:
Let $R(n,m)$ be a $\Sig^{*}_{2}$-formula for which
$T_{1}(\Natural)\vdash\fal n\exi m R(n,m)$.
%$T_{1}(\Natural)\vdash\fal n_{1},\ldots,n_{k}\exi m_{1},\ldots,m_{\ell}\, R(n_{1},\ldots,n_{k},m_{1},\ldots,m_{\ell})$.

Then there exists an $n_{0}\in\mathbb{N}$, an initial hydra $a_{0}$ and a finite set $\calf_{0}$ such that
$\fal n\exi m\leq 1+h_{a_{0}}^{\calf_{0}}(n)[n_{0}\leq n \Rarw R(n, m)]$ holds.
\elem

Lemma \ref{lem:consis} is shown in the next sections \ref{sect:finiteproof} and \ref{sect:consisprf}.
Assuming Lemma \ref{lem:consis}, let us show
Theorem \ref{th:hydra}.\ref{th:hydra1}.
Suppose $T_{1}(\Natural)$ proves the statement 
$\mbox{{\rm (H)}}_{\ome_{1}}\Lrarw\left(\fal a,\calf, n\exi m(h_{a}^{\calf}(n)\simeq m)\right)$,
and hence $\fal a,\calf, n\exi k\exi m(h_{a}^{\calf}(n)\leq m\land k=m+2)$,
where $a$ ranges over
initial hydras, $\calf$ over finite subsets of $\calf_{\mu}$.
$h_{a}^{\calf}(n)\leq m$ denotes a formula $R(\lceil a\rceil,\lceil\calf\rceil, n,m)$ with 
a $\Sig^{*}_{2}$-formula $R$ saying that
`for any sequence $\sig=(\sig_{0},\ldots,\sig_{m-n-1})$ of hydras $\sig_{i}$, if $\sig_{0}=a$,
and $\fal i<m-n-1\left(\sig_{i+1}\in\sig_{i}[n+i]\cup\{0^{N}\}\right)$
%`$T\subset{}^{<\ome}\ome$ is a finite tree of hydras $T(\sig)\,(\sig\in T)$ of depth $dp(T)=m-n$ such that 
%$T(\emptyset)=a$ for the root, 
%and for any $\sig\in T$,
%the set $\{T(\sig*(k))\}_{k}$ of hydras on the sons $\sig*(k)$ of $\sig$ in $T$ is equal to the set $(T(\sig))[n+lh(\sig)]$ 
with respect to $\calf$,
then $\sig_{m-n-1}=0^{N}$'.
with the code $0^{N}$ of the zero hydra $0$.
By Proposition \ref{prp:H0}.\ref{prp:H0.1}, `for any sequence $\sig$ of hydras' is a bounded quantifier.

By Lemma \ref{lem:consis} pick an $n_{0}\in\mathbb{N}$, an initial hydra $a_{0}$ and a finite set $\calf_{0}$
such that $2+h_{a}^{\calf}(n)\leq 1+h_{a_{0}}^{\calf_{0}}(\max\{\lceil a\rceil,\lceil\calf\rceil,n\})$ holds
for any $a,\calf,n$ such that
 $\max\{\lceil a\rceil,\lceil\calf\rceil,n\}\geq n_{0}$.
Let $a=a_{0}$, $\calf=\calf_{0}$.
Then $h_{a_{0}}^{\calf_{0}}(n)<h_{a_{0}}^{\calf_{0}}(n)$ for any 
$n\geq\max\{\lceil a_{0}\rceil,\lceil\calf_{0}\rceil,n_{0}\}$.
This is a contradiction.

\bcor\label{th:hydra2}
$T_{1}$ does not prove the full statement $\mbox{{\rm (H)}}^{set}_{\ome_{1}}$ in the set-theoretic language.
\ecor
\bprf
Corollary \ref{th:hydra2} follows from Theorem \ref{th:hydra}.\ref{th:hydra1} as follows.
Let $n\simeq x$ denote the relation between natural numbers $n$ and ordinals $x$ such that
$n\simeq x$ iff there exists a bijection between $\{0,\ldots,n\dot{-}1\}$ and $\{y\in Ord: y<x\}$.
It is clear that $T_{1}(\Natural)$ proves that $0^{N}\simeq 0^{ON}$, $n\simeq x \Lrarw n+1\simeq x\cup\{x\}$
and $\fal n\in\Natural\exi ! x<\ome(n\simeq x)\land \fal x<\ome\exi n!\in\Natural(n\simeq x)$.
Moreover for each primitive recursive relation R, we have in $T_{1}(\Natural)$ that 
$\bigwedge_{i}(n_{i}\simeq x_{i}) \to (R(n_{1}, \ldots,n_{k}))\lrarw R^{set}(x_{1},\ldots,x_{k}))$
for the set-theoretic counter part $R^{set}$ of $R$. For example
$n_{1}\simeq x_{1}\land n_{2}\simeq x_{2}\to(n_{1}<^{N}n_{2}\lrarw x_{1}\in x_{2})$.

From this we see that $T_{1}(\Natural)$ proves the equivalence
$\mbox{{\rm (H)}}_{\ome_{1}}\lrarw \mbox{{\rm (H)}}^{set}_{\ome_{1}}$.
\eprf

\section{Finite proof figures}\label{sect:finiteproof}

In this section \ref{sect:finiteproof} and the next section \ref{sect:consisprf}
we work in the theory $T_{1}^{+}(\Natural)=T_{1}(\Natural)+TI(\veps_{\rho_{0}+1})$.

In this section  an extension $T_{c}(\Natural)$ of the theory $T_{1}(\Natural)$ with individual constants and function constants
is formulated in one-sided sequent calculus, and
permissible ordinal assignments to sequents occurring in proofs are defined in subsection \ref{subsec:ordinalassignment}.
Each proof in $T_{1}(\Natural)$ is shown to have a permissible ordinal assignment in subsection \ref{subsec:embedding}.

The language $\calL_{c}$ of $T_{c}(\Natural)$ is obtained from the language $\mathcal{L}(\Natural,\in)$ of $T_{1}(\Natural)$
by adding names (individual constants) $c_{\alp}$ 
of each $a=D_{i}(c\oplus b)\in Tm(\calf_{0})\, (i=0,1, b\neq 0)$, 
names $c_{a}$ of each $a=F(c\oplus b)\in Tm(\calf_{0})$, 
and
`function symbols' $f_{A}(y_{1},\ldots,y_{n})\in\calf_{\mu}$ for each $\Del_{0}$-formula $A(x;y_{1},\ldots,y_{n})$ in
the set-theoretic language $\{\in\}$.
The constant $c_{a}$ is identified with $a\in Tm(\calf_{0})$.
Formulas are assumed to be in negation normal form.

%$\mbox{Hull}(x)$ denotes the $\Sig_{1}$-Skolem hull of sets $x$ of sets on $L_{\rho_{0}}$.

\bdf
\benu

\item
\textit{Terms} {\rm in $\calL_{c}$ are generated as follows.}
\benu
\item
{\rm Each variable and each constant $0^{N},0^{ON},\emptyset,D_{i}(c\oplus b), F(c\oplus b)$ is a term.}

\item
{\rm If $t,s$ are terms, then so are $S(t), \ome^{t}, J(t,s)$.}

\item
{\rm If $t_{1},\ldots,t_{n}\, (n>1)$ are terms, then so are $t_{1}+\cdots+t_{n}$ and $t_{1}\cdot t_{2}$.}

\item
{\rm For $\Del_{0}$-formula $A(x;y_{1},\ldots,y_{n})$ with $f_{A}\in\calf_{\mu}$ and} closed {\rm terms $t_{1},\ldots,t_{n}$,
\\
$\mu x. A(x;t_{1},\ldots,t_{n})\equiv f_{A}(t_{1},\ldots,t_{n})$ is a} closed {\rm term.}
\eenu

\item
{\rm A term in $\calL_{c}$ is a \textit{well formed term} if it is one of number terms, set terms
or ordinal terms defined below.
}
\item
{\rm
A term $t$ is a \textit{number term} iff $t\equiv(S(\cdots(S(u))\cdots))$, where $u$ is either a variable or 
$u\equiv 0^{N}$. A closed number term $S(\cdots(S(0^{N}))\cdots)$ with $k$ times successor function $S$ is 
a \textit{numeral} denoted by $\bar{k}$.
}

\item
{\rm
\textit{Set terms} in $\calL_{c}$ are generated as follows.
 \benu
 \item
 {\rm Each variable and the constant $\emptyset$ is a set term.
 }
 \item
 {\rm If $t$ is a set term and $s$ a well formed term, then $J(s,t)$ is a set term. 
 }
  \eenu
  
 \item
 {\rm
 \textit{Ordinal terms} in $\calL_{c}$ are generated as follows.}
 \benu
 \item
 {\rm
  Each variable and each constant $0^{ON},D_{i}(c\oplus b), F(c\oplus b)$ is an ordinal term.}
  \item
   {\rm If $t$ is an ordinal term, then so is $\ome^{t}$.}
  \item
  {\rm If $t_{1},\ldots,t_{n}\, (n>1)$ are ordinal terms, then so are $t_{1}+\cdots+t_{n}$ and $t_{1}\cdot t_{2}$.}
  \item
  {\rm For $\Del_{0}$-formula $A(x;y_{1},\ldots,y_{n})$ with $f_{A}\in\calf_{\mu}$ and \textit{closed} terms $t_{1},\ldots,t_{n}$,
$\mu x. A(x;t_{1},\ldots,t_{n})\equiv f_{A}(t_{1},\ldots,t_{n})$ is a \textit{closed} ordinal term.}
 \eenu
}
%$F_{x\cup\{\ome_{1}\}}(y)$ {\rm denotes the Mostowski collapse} $F_{x\cup\{\ome_{1}\}}:\mbox{{\rm Hull}}(x\cup\{\ome_{1}\})\lrarw \gam$
%{\rm with} $x=F_{x\cup\{\ome_{1}\}}(\ome_{1})$ {\rm and} $\gam=F_{x\cup\{\ome_{1}\}}(\rho_{0})$.
\eenu
\edf

The \textit{value} $v(t)\in\Natural\cup\rho_{0}\cup \mbox{HF}_{\Natural\cup\rho_{0}}$ of closed terms $t$ is defined as follows.
$v(t)=0<\rho_{0}$ when $t$ is not a well formed term.
$v(S(t))=v(t)+1$ for number terms $t$, i.e., $v(\bar{k})=k\in\Natural$.
$v(J(t,s))=v(t)\cup\{v(s)\}$ for set terms $t$ and well formed $s$.

%Each term in $\mathcal{L}_{c}$ is a term in $Tm(\calf_{\mu})$.?

\bdf
\benu
\item
{\rm A \textit{literal} is one of atomic formulas $N(t), R(t_{0},\ldots,t_{n-1})$,
$ON(t), t_{0}<t_{1}, R^{\cala}(t_{0}t_{1}), P(t_{0},t_{1}), P_{\rho_{0}}(t)$, $Set(t),s\in t, s=t$
or their negations, where $R$ is a relation symbol for an $n$-ary primitive recursive relation on integers.}

\item
{\rm The \textit{truth}  of closed literals is defined as follows.}
 \benu
 \item
 {\rm
 $ON(t)$ is \textit{true} if $v(t)$ is an ordinal.
$N(t)$ is \textit{true} if $t$ is a numeral.
$Set(t)$ is \textit{true} if $t$ is a closed set term.
}

 \item
 {\rm
  Let $R$ be a relation symbol for a primitive recursive relation on $\Natural$.
  Then $R(t_{1},\ldots,t_{n})$ is \textit{true} if all of $t_{1},\ldots,t_{n}$ are numerals and
  $\Natural\models R(t_{1},\ldots,t_{n})$ holds.
  }
  
  \item
  {\rm
$s< t$ is \textit{true} if $v(s),v(t)<\rho_{0}$ and $v(s)< v(t)$. 
$R^{\cala}(s,t)$ is \textit{true} if $v(s),v(t)<\rho_{0}$ and $R^{\cala}(v(s),v(t))$ holds.
$P(s,t)$ is \textit{true} if $v(s),v(t)<\rho_{0}$ and $v(s)=x=\Psi_{\ome_{1}}(\bet)$ and
$v(t)=F_{x\cup\{\ome_{1}\}}(\rho_{0})$ for some $\bet$.
$P_{\rho_{0}}(t)$ is \textit{true} if $v(t)<\rho_{0}$ and $v(t)=\Psi_{\rho_{0}}(\bet)$ for some $\bet$.
}

\item
{\rm
$s\in t$ is \textit{true} if $t$ is a closed set term and $v(s)\in v(t)$ holds.
$s=t$ is \textit{true} if $s$ and $t$ are closed term in the same sort, and $v(s)=v(t)$ holds.
}

\item
{\rm
 A closed literal $\lnot L$ is \textit{true} if $L$ is not true.
 }
  \eenu

\item
{\rm An} $E$\textit{-formula} {\rm is either a literal or a formula of one of the shapes $A_{0}\lor A_{1}, \exi x\, A(x)$.}
\eenu
\edf
If a formula is obtained from a $\Del_{0}$-formula in the language $\mathcal{L}_{2}$
by substituting $\mathcal{L}_{c}$-terms for variables, then the formula is a \textit{$\Del_{0}$-formula}.
By the definition the predicates $P,P_{\rho_{0}}$, $N,ON,\in$ do not occur in $\Del_{0}$-formulas.

%By \textit{$\Del_{0}$-formula} we mean a bounded formula in the language $\mathcal{L}_{c}$ in which predicates $P,P_{\rho_{0}}$
%do not occur.

The \textit{truth} of $\Del_{0}$-sentences is defined from one of literals.
A \textit{$\Sig_{1}$-formula} or a \textit{$\Pi_{1}$-formula} is defined similarly.
%These formulas are obtained from formulas in $\mathcal{L}_{1}\cup\{\ome_{1}\}$ by substituting $\mathcal{L}_{c}$-terms for variables.

%$N(t)$ is denoted by $t\in\Natural$.
%$s\in t\to t\not\in\Natural$, $R(t_{0},\ldots,t_{n-1})\to \bigwedge_{i<n}(t_{i}\in\Natural)$.

The following are axioms and inference rules in $T_{c}(\Natural)$.
\textit{Proof figures} are constructed from these axioms and inference rules.

Relations between occurrences $A,B$ of formulas in a proof such as `$A$ is a \textit{descendant} of $B$'
or equivalently `$B$ is an \textit{ancestor} of $A$',
and `an occurrence of inference rule is \textit{implicit} or \textit{explicit}'
 are defined as in 
 \cite{ptMahlo}.
\\

\noindent
[{\bf Axioms}]
\[
\infer[(ax)]{\Gam,A}{}
\]
where $A$ is either a true closed literal or a true closed $\Del_{0}$-formula or
an arithmetic axiom whose universal closure holds in the standard model $\Natural$ or
an ontological axiom or the defining axiom for $J$ or the axiom of Extensionality.

\[
\infer[(taut)]{\Gam,\lnot A, A}{}
\msfiv\mbox{   for literals and $\Del_{0}$-formulas $A$.}
\]
This means $ \dg(A)=1$ in Definition \ref{df:depthfml} below.
If there occurs no fee variable in an axiom $(ax), (taut)$,
then it contains either a true literal or a true $\Del_{0}$-sentence.

Cf.\,(\ref{eq:Z2}).
\[
\infer[(P\exi)]{\Gam, (s\not<\ome_{1},) \exi x, y<\ome_{1}[s<x\land P(x,y)]}{}
\]
When $s<\ome_{1}$ is a true literal, $s\not<\ome_{1}$ may be absent.
%, and $x<\ome_{1}:\equiv x\in\ome_{1}$.

Cf.\,(\ref{eq:Z5}).
\[
\infer[(P_{\rho_{0}}\exi)]{\Gam,  (\lnot ON(s),)\exi x[s<x\land P_{\rho_{0}}(x)]}{}
\]
When $ON(s)$ is a true literal, $\lnot ON(s)$ may be absent.

\bdf
{\rm
A term $t$ is said to be an \textit{$N$-simple term} iff
if
$t\equiv S(t_{0})$ for a term $t_{0}$, then either $t_{0}$ is a numeral
or a variable.

A term $t$ is said to be an \textit{$S$-simple term} iff
if $t\equiv J(s_{0},s_{1})$, then either $s_{0}$ is a closed set term or a variable.
}
\edf
[{\bf Inference rules}]
In each case the main (principal) formula is assumed to be in the lower sequent $\Gam$.
Namely
$(A_{0}\lor A_{1})\in\Gam$ in $(\lor)$, $(A_{0}\land A_{1})\in\Gam$ in $(\land)$,
$(\exi x\, A(x))\in\Gam$ in $(\exi)$, $(\exi x< t\, A(x))\in\Gam$ in $(b\exi)$, $(\fal x\, A(x))\in\Gam$ in $(\fal)$, $(\fal x< t\,A(x))\in\Gam$ in $(b\fal)$.

The variable $x$ in $(\fal),(b\fal)$ is an eigenvariable.

\[
\infer[(\lor)]{\Gam}{\Gam,A_{i}}
\:
\infer[(\land)]{\Gam}{\Gam, A_{0} & \Gam,A_{1}}
\]
\[
\infer[(b\exi)]{(s\not< t,)\Gam}
{
\Gam,A(s)
}
\:
\infer[(\fal)]{\Gam}{\Gam, A(x)}
\:
\infer[(b\fal)]{\Gam}{\Gam,x\not< t, A(x)}
\]
where in $(b\exi)$, the formula $s\not< t$ may be absent when $s< t$ is a closed true literal.
\[
\infer[(\exi)]{\Gam}{\Gam, A(s)}
\]
where $A(x)$ is not of the form $N(x)\land A_{0}(x)$ nor $Set(x)\land A_{0}(x)$.
\[
\infer[(\exi)^{N}]{\Gam}
{
\Gam, N(s)
&
\Gam, A(s)
}
\]
where $\exi x(N(x)\land A(x))$ is in $\Gam$, and the \textit{instance term} $s$ is $N$-simple.

Let $s$ be a non-simple term, and $t$ be a term such that
$s\equiv S(t)$.
Assume that we have proofs of $\Gam,N(s)$ and of $\Gam,A(s)$.
Then $\Gam,\exi x(N(x)\land A(x))$ is derivable using the restricted inference $(\exi)^{N}$ as follows.
For simplicity assume that $t$ is $N$-simple.
{\small
\[
\infer{\Gam,\exi x(N(x)\land A(x))}
{
 \infer[(\exi)^{N}]{\Gam,\exi x(N(x)\land x=t)}
 {
  \infer{\Gam,N(t)}
  {
   \infer*{\Gam,N(St)}{}
   &
   \lnot N(St),N(t)
   }
   &
   \hskip-1.0cm
   t=t
   }
  &
  \hskip-0.2cm
  \infer{\lnot \exi x(N(x)\land x=t),\Gam,\exi x(N(x)\land A(x))}
  {
  \infer[(\exi)^{N}]{\lnot N(y), y\neq t,\Gam,\exi x(N(x)\land A(x))}
   {
    \lnot N(y),N(Sy)
    &
    \hskip-1.3cm
    \infer{y\neq t,\Gam,A(Sy)}
    {
     \infer*{\Gam,A(St)}{}
     &
     \hskip-0.5cm
     \infer{y\neq t,\lnot A(St),A(Sy)}
     {
      y\neq t, Sy=St
      &
      \infer*{Sy\neq St,\lnot A(St),A(Sy)}{}
     }
    }
   }
  }
}
\]
}
where $y$ is a fresh variable,
both $t$ and $Sy$ are $N$-simple, both $\lnot N(St),N(t)$ and $\lnot N(y),N(Sy)$ are ontological axioms,
and
both $t=t$ and $y\neq t,Sy=St$ are equality axioms. 

\[
\infer[(\exi)^{S}]{\Gam}
{
\Gam, Set(s)
&
\Gam, A(s)
}
\]
where $\exi x(Set(x)\land A(x))$ is in $\Gam$, and the \textit{instance term} $s$ is $S$-simple.
As in the case for $(\exi)^{N}$, we can restrict inferences for introducing existential quantifiers on sets to ones 
with $S$-simple instance terms.

\[
\infer[(ind)_{<}]{(s\not< t,)\Gam}{\Gam,\lnot ON(\alp),\lnot\fal \bet<\alp A(\bet), A(\alp) & \Gam,\lnot A(s)}
\]
where 
$s\not< t$ may be absent in the lower sequent when
$s<t$ is a true closed literal.
The formula $A(x)$ is the \textit{induction formula}, 
the term $t$ is the \textit{induction term}, and
the left upper sequent $\Gam,\lnot ON(\alp),\lnot\fal \bet<\alp A(\bet), A(\alp)$ is the \textit{induction sequent}
  of the $(ind)_{<}$.
The variable $\alp$ is the \textit{eigenvariable} of the rule $(ind)_{<}$. 
The \textit{degree} of the $(ind)_{<}$ is defined to be $\dg(\fal \bet< s A(\bet))$.
\[
\infer[(ind)_{\Natural}]{(\lnot N(s),)\Gam}
{
\Gam, A(\bar{0})
&
\Gam,\lnot N(n), \lnot A(n), A(Sn) 
& \Gam,\lnot A(s)}
\]
where 
$\lnot N(s)$ may be absent in the lower sequent when
$N(s)$ is a true closed literal, i.e., $s$ is a numeral.
The formula $A(x)$ is the \textit{induction formula}, 
\textit{induction term} $s$ is $N$-simple
%either a variable or a closed term,
%numeral or closed with $\lnot N(s)$, 
and
the middle upper sequent $\Gam,\lnot N(n), \lnot A(n), A(Sn) $ is the \textit{induction sequent}
 of the $(ind)_{\Natural}$.
The variable $n$ is the \textit{eigenvariable} of the rule $(ind)_{\Natural}$. 
The \textit{degree} of the $(ind)_{\Natural}$ is defined to be $\dg(A(s))$.
\[
\infer[(ind)_{\in}]{(\lnot Set(s),)\Gam}
{
\Gam,A(\emptyset)
&
\Gam,\lnot Set(x), \lnot A(x),A(J(x,y))
& 
\Gam,\lnot A(s)}
\]
where $\lnot Set(s)$ may be absent in the lower sequent when $Set(s)$ is a true closed literal.
The formula $A(x)$ is the \textit{induction formula},
the \textit{induction term} $s$ is $S$-simple,
%closed set term or closed with $\lnot Set(s)$, and
the middle upper sequent $\Gam,\lnot Set(x), \lnot A(x),A(J(x,y))$ is the \textit{induction sequent}
 of the $(ind)_{\in}$.
The variables $x,y$ are the \textit{eigenvariable} of the rule $(ind)_{\in}$. 
%F(\emptyset)\land \fal x(\lnot Set(x)\to F(x)) \land \fal x,y(Set(x)\land F(x) \to F(J(x,y)))\to \fal x\, F(x)
The \textit{degree} of the $(ind)_{\in}$ is defined to be $\dg(A(s))$.

$(ind)$ denotes one of these three induction schemata $(ind)_{<},(ind)_{\Natural},(ind)_{\in}$.

\[
\infer[(cut)]{\Gam,\Lam}{\Gam,\lnot A & A,\Lam}
\]
$A$ is an $E$-formula called the \textit{cut formula} of the $(cut)$.
\[
\infer[(Rfl)]{(\lnot ON(t),)\Gam}
{
\Gam, \fal x< t\, A(x)
&
\lnot ON(y), t\not< y,\exi x< t\, \lnot A^{(y)}(x),\Gam
}
\]
where $\lnot ON(t)$ may be absent in the lower sequent when $ON(t)$ is a true closed literal.
$t$ is a term, $y$ is an eigenvariable,  and 
$A(x)\equiv(\exi z\exi w[P_{\rho_{0}}(z)\land B(x,w)])\, (B\in\Del_{0})$, 
$A^{(y)}(x):\equiv(\exi z< y\exi w< y[P_{\rho_{0}}(z)\land B(x,w)])$,
cf.\,(\ref{eq:Z6}).
%$\tau(x,u,v)\equiv[u,v\in L_{x}(\Natural) \land \fal w\in L_{x}(\Natural) \tht(u,v,w)]$.
%\beqn\label{eq:Z6}
%\fal u\in a\, A(u) \to \exi c[tran(c)\land a\in c\land \fal u\in a\, A^{(c)}(u)]
%where $A(u)\equiv(\exi x\in P_{\rho_{0}}\exi v\, \tau(x,u,v))$ for  
%$A^{(c)}(u)\equiv(\exi x\in P_{\rho_{0}}\cap c\exi v\in c\,\tau)$.

\[
\infer[(P\Sig_{1})]{\Gam, (\lnot P(t_{0},t_{1}), s\not<t_{0},)  \vphi^{t_{1}}[t_{0},s]}{\Gam,\vphi[\ome_{1},s]}
\]
$\vphi$ is an arbitrary $\Sig_{1}$-formula in the set-theoretic language $\{\in\}$,
cf.\,(\ref{eq:Z1}).
When $P(t_{0},t_{1})$ or $s<t_{0}$ is a true literal, these may be absent.
%\beqn\label{eq:Z1}
%P(x,y) \to a\in L_{x}  \to \vphi[\ome_{1},a] \to \vphi^{y}[x,a]

\[
\infer[(P_{\rho_{0}}\Sig_{1})]{\Gam,(\lnot P_{\rho_{0}}(t), s\not<t,) \vphi^{t}[s]}{\Gam,\vphi[s]}
\]
$\vphi$ is an arbitrary $\Sig_{1}$-formula in the language $\{\in\}$, cf.\,(\ref{eq:Z4}).
When $P_{\rho_{0}}(t)$ or $s<t$ is a true literal, these may be absent.
%\beqn\label{eq:Z4}
%P_{\rho_{0}}(x) \to a\in L_{x} \to \vphi[a] \to \vphi^{x}[a]

\[
\infer[(h)]{\Gam,\Del}{\Gam}
\]

\[
\infer[(D_{1})_{\alp}]{\Lam,\Gam^{(\alp)}}{\Lam, \Gam}
\]
{\rm where} 
$\alp=D_{1}(c_{1}\oplus\alp_{0})$ for some $c_{1}\oplus\alp_{0}\neq 0$ with $c_{1}=stk(\alp)$.
Each formula in $\Gam$ is one of the closed formulas 
$\fal x<t\, A(x)$, $A(s_{0})$, and $\exi w[P_{\rho_{0}}(s_{1})\land B(s_{0},s_{1},w)]$,
where $B$ is a $\Del_{0}$-formula, 
$A(x)\equiv(\exi z\exi w[P_{\rho_{0}}(z)\land B(x,z,w)])$
$t,s_{0},s_{1}$ are closed ordinal-terms.
Each \textit{implicit} formula in $\Lam$ is a bounded sentence.
Note that there occurs no unbounded universal quantifier in implicit formulas in $\Lam\cup\Gam$.

\[
\infer[(D_{0})_{\alp}]{\Lam}{\Lam}
\]
where each formula in $\Lam$ is either a false closed $\Del_{0}$-formula or
a closed subformula of a $\Sig_{2}^{*}$-sentence.
$\alp=D_{0}(c_{0}\oplus \alp_{0})$ for some $c_{0}\oplus \alp_{0}$.
$c_{0}=stk(\alp)$ is the \textit{stock} of the $(D_{0})_{\alp}$.

\[
\infer[(pad)_{b}]{\Gam,\Del}{\Gam}
\msfiv
\infer[{}_{c}(pad)]{\Gam,\Del}{\Gam}
\]
for $b,c\in H(\calf_{\mu})$.

\subsection{Ordinal assignment}\label{subsec:ordinalassignment}
In this subsection let us define permissible ordinal assignments.

\bdf\label{df:height}
{\rm The \textit{height} $h(\Gam)=h(\Gam;\calP)<\ome\cdot 2$ of sequents $\Gam$ in a proof figure $\calP$.}
\benu
\item
$h(\Gam)=0$ {\rm if $\Gam$ is the end-sequent of $\calP$.}

\item
$h(\Gam)=\ome\cdot i$ {\rm if $\Gam$ is the upper sequent of a $(D_{i})$.}

\item
$h(\Gam)=h(\Del)+1$ {\rm if $\Gam$ is the upper sequent of an $(h)$ with its lower sequent $\Del$.}

\item
$h(\Gam)=h(\Del)$ {\rm if $\Gam$ is an upper sequent of a rule other than $(h)$ and $(D_{i})$ with its lower sequent $\Del$.}
\eenu
{\rm Let
$h_{0}(\Gam)=h(\Gam)$ if $h(\Gam)<\ome$.
$h_{0}(\Gam)=h(\Gam)-\ome$ if $h(\Gam)\geq \ome$.}
\edf

\bdf\label{df:depthfml}{\rm The \textit{degree} $ \dg(A)<\ome$ of formulas $A$.}
\benu
\item
$ \dg(A)=1$ {\rm if $A$ is either a literal or a $\Del_{0}$-formula.}

{\rm In what follows} $A$ {\rm is neither a literal nor a $\Del_{0}$-formula.}

\item
$ \dg(A)= \dg(A_{0})+ \dg(A_{1})+2$ {\rm if $A\equiv(A_{0}\lor A_{1}), (A_{0}\land A_{1})$.}

\item
$ \dg(A)= \dg(B)+1$ {\rm if $A\equiv(\exi x\, B(x)),(\fal x\, B(x))$.}

\item
$ \dg(A)= \dg(B)+1$ {\rm if $A\equiv(\exi x< t\, B(x)),(\fal x< t\, B(x))$.}
\eenu
\edf

\bdf
{\rm A proof figure is said to be \textit{height regulated} if it enjoys the following conditions:}
\bdes
%\item[(h0)] 
%{\rm The end-sequent of $\calP$ consists solely of closed false $\Del_{0}$-formulas.}

\item[(h1)] 
{\rm There occurs no free variable in any sequent $\Gam$ if $h(\Gam)<\ome$.}

\item[(h2)] 

{\rm Let $\Gam,\exi x[s<x\land P_{\rho_{0}}(x)]$ be an axiom $(P_{\rho_{0}}\exi)$ in $\mathcal{P}$,
and $J$ be a $(cut)$ whose cut formula is a descendant $C\equiv(\exi x[s<x\land P_{\rho_{0}}(x)])$
of $C$ in the axiom.
Then $h(\Del)\geq\ome$ for the upper sequent $\Del$ of the $(cut)\, J$.}

\item[(h3)] 
{\rm For any $(cut)$ in $\calP$,
$ \dg(C)\leq h_{0}(\Gam,\Del)$ for its cut formula $C$ and the lower sequent $\Gam,\Del$.}

\item[(h4)]
{\rm For any $(ind)\, J$ in $\calP$ with its lower sequent $\Gam$,
$\ome+ \dg(J)\leq h(\Gam)$ holds, and
there are no nested $(ind)$ rules, i.e., there occurs no $(ind)$ above the rule $(ind)$.}

\item[(h5)]
{\rm There exists a rule $(D_{1})$ below a $(Rfl)$. 
Let $J$ be the lowest such rule $(D_{1})$ with the lower sequent $\Del$.
Then $h(\Del)\geq \dg(\exi x< t \lnot A^{(y)}(x))$.}
\[
\infer[(D_{1})\, J]{\Del}
{
 \infer*{\cdots}
   {
    \infer[(Rfl)]{\Gam}
     {
     \Gam, \fal x< t\, A(x)
     &
    t\not< y, \exi x< t \lnot A^{(y)}(x),\Gam
     }
 }
}
\]

\item[(h6)]
{\rm If a rule $(D_{1})\, J_{0}$ is above another $(D_{1})\, J_{1}$, then the only rules between $J_{0}$ and $J_{1}$ are
$(D_{1})$'s.}

%\item[(h7)]
%{\rm $\calP$ ends with an inference rule $(D_{0})$.}

\edes
\edf

\bdf\label{df:ordinalassignment}
{\rm Let $\calf_{0}\subset\calf_{\mu}$ be a finite set of function symbols.}
{\rm An} $\calf_{0}$-ordinal assignment {\rm for a proof figure $\calP$ attaches a hydra (an ordinal) $o(\Gam)\in H(\calf_{0})$
 to each occurrence of a sequent $\Gam$ in $\calP$
which enjoys the following conditions. Let us write 
\[
\infer{\Gam; b}{\cdots\Gam_{i}; a_{i}\cdots}
\]
when the lower sequent $\Gam$ receives an ordinal $b$, i.e., $o(\Gam)=b$, and
$o(\Gam_{i})=a_{i}$ for upper sequents $\Gam_{i}$.}
\\
{\bf Axioms} $\Gam$.
\benu

\item
{\rm If} $\Gam$ {\rm is one of axioms $(ax),(taut)$, then $o(\Gam)=1=D_{0}(0)$.}

\item
{\rm For a $(P_{\rho_{0}}\exi)\, \Gam$, $o(\Gam)=D_{2}(0)$.}

\item
{\rm For a} $(P\exi)\, \Gam$, $o(\Gam)=D_{1}(0)$.
\eenu

\noindent
{\bf Rules}.
{\rm Let $\Gam$ be the lower sequent of a rule $J$ with its upper sequents $\Gam_{i}$:}
\[
\infer[J]{\Gam}{\cdots\Gam_{i}\cdots}
\]

\benu

\item
{\rm $J$ is one of the rules $(\fal),(b\fal),(P\Sig_{1})$ or $(P_{\rho_{0}}\Sig_{1})$:}
 $o(\Gam)=o(\Gam_{0})$.
 
%\item\label{df:ordinalassignmentbfal}
%{\rm $J$ is a $(b\fal)$:
%$o(\Gam)=o(\Gam_{0})+b$ for some $0\neq b\in  H(\calf_{0})$ if $J$ is explicit with a $\Del_{0}$-main formula.
% In this case the rule is denoted $(b\fal)_{b}$.
%Otherwise}  $o(\Gam)=o(\Gam_{0})$.

\item
{\rm $J$ is either a $(\land)$ or a $(cut)$:
$o(\Gam)=o(\Gam_{0})+o(\Gam_{1})$.}

\item
{\rm $J$ is one of rules $(\lor), (\exi), (b\exi)$:
 $o(\Gam)=o(\Gam_{0})+b$ for some $0\neq b\in  H(\calf_{0})$.
 In this case the rule is denoted $(\lor)_{b}, (\exi)_{b}$, etc.}

\item
{\rm $J$ is an $(\exi)^{N}$:
 $o(\Gam)=o(\Gam_{0})+o(\Gam_{1})+b$ for some $0\neq b\in  H(\calf_{0})$.
 In this case the rule is denoted $(\exi)^{N}_{b}$.}
 
 \item
{\rm $J$ is an $(\exi)^{S}$:
 $o(\Gam)=o(\Gam_{0})+o(\Gam_{1})+b$ for some $0\neq b\in  H(\calf_{0})$.
 In this case the rule is denoted $(\exi)^{S}_{b}$.}
 
\item 
{\rm $J$ is a $(pad)_{b}$:
 $o(\Gam)=o(\Gam_{0})+b$ for $b\in  H(\calf_{0})$.}
 
\item 
{\rm $J$ is a ${}_{b}(pad)$:
 $o(\Gam)=b+o(\Gam_{0})$ for $b\in  H(\calf_{0})$.}

\item
{\rm $J$ is an $(h)$:}
 $o(\Gam)=D_{2}(o(\Gam_{0}))$.

\item
{\rm $J$ is a $(Rfl)$:}
 $o(\Gam)=o(\Gam_{0})+o(\Gam_{1})+D_{2}(0)$.

\item
{\rm $J$ is an $(ind)_{<}$:
\[
\infer[(ind)_{<}]{(s\not< t,) \Gam;b}{\Gam,\lnot ON(\alp),\lnot\fal \bet< \alp A(\bet), A(\alp);a_{1} & \Gam,\lnot A(s);a_{2}}
\]
Let $mj(t)=D_{2}(0)$ if $t$ is not closed.
Otherwise $mj(t)=t'$ for some $t'\in Tm(\calf_{0})\cup\{D_{2}(0)\}$ such that 
$v(t)\leq v(t')$.

Then
$b= (a_{1}+a_{2}+1)\times mj(t)$, cf.\,{\bf (p1)} below.}

\item
{\rm $J$ is an $(ind)_{\Natural}$:
\[
\infer[(ind)_{\Natural}]{(\lnot N(s),)\Gam;b}
{
\Gam, A(\bar{0});a_{0}
&
\Gam,\lnot N(n), \lnot A(n), A(Sn) ;a_{1}
& \Gam,\lnot A(s);a_{2}
}
\]
$b=a_{0}+a_{2}+a_{1}\otimes \ome$.
}

\item
{\rm $J$ is an $(ind)_{\in}$:
\[
\infer[(ind)_{\in}]{(\lnot Set(s),)\Gam;b}
{
\Gam,A(\emptyset);a_{0}
&
\Gam,\lnot Set(x), \lnot A(x),A(J(x,y));a_{1}
& 
\Gam,\lnot A(s);a_{2}
}
\]
$b=a_{0}+a_{2}+a_{1}\otimes \ome$.
}

\item
{\rm $J$ is a rule $(D_{1})$:}
\[
o(\Gam)=
\left\{
\begin{array}{ll}
D_{1}(c_{1}\oplus o(\Gam_{0})) & \mbox{{\rm if }} h(\Gam)<\ome
\\
o(\Gam_{0}) & \mbox{{\rm if }}  h(\Gam)=\ome
\end{array}
\right.
\]
{\rm where $stk(o(\Gam))=c_{1}$ for a list $c_{1}$.}

\item
{\rm $J$ is a rule $(D_{0})$:}
\[
o(\Gam)=D_{0}(c_{0}\oplus o(\Gam_{0})) 
\]
{\rm where $stk(o(\Gam))=c_{0}$ for a list $c_{0}$.}

\eenu
{\rm Finally let
$o(\calP)=o(\Gam_{end})$ for the end-sequent $\Gam_{end}$ of $\calP$.}
\edf
Note that by {\bf (h6)}, there are rules $(D_{1})$ consecutively.
\[
\infer[(D_{1})]{\Gam_{0}}
{
\infer{\Gam_{1}}
 {
  \infer*{}
 {
 \infer[(D_{1})]{\Gam_{n-1}}{\Gam_{n}}
 }
}
}
\]
with $h(\Gam_{0})<\ome$.
Then $o(\Gam_{1})=\cdots=o(\Gam_{n})$ and $o(\Gam_{0})=D_{1}(c_{1}\oplus o(\Gam_{1}))$
for a list $c_{1}$.
We write $stk(o(\Gam_{i+1}))=c_{1}$ and $\Gam_{i+1}; c_{1}\oplus o(\Gam_{i+1})$ for any $i<n$.
Likewise for the upper sequent $\Gam_{1}$ of a $(D_{0})$, we write $\Gam_{1}; c_{0}\oplus o(\Gam_{1})$.

\blem\label{lem:depth}{\rm (Tautology lemma)}\\
For any formula $A(x)$,
there exist ordinal assignments $o$ such that 
\\
$o(\Gam,\lnot A(t),A(t))= \dg(A(x))$
for any $\Gam$ and any term $t$.
\elem
\bprf 
By induction on $ \dg(A)$.
To get $o(\Gam,\lnot A(t),A(t))= \dg(A(x))$, 
use $(\lor)_{b},(\exi)_{b}$, etc. for $b=1$.
Note that in the inference rules for introducing unbounded existential quantifiers $\exi x(N(x)\land\cdots)$
for the predicate $N$, the instance terms are variables.
\eprf

\bdf
{\rm 
For formulas $B$ (possibly with variables), $|B|$ denotes the total number of occurrences of symbols
$0,+,\cdot,\ome,\oplus,D_{0},D_{1},D_{2}, F,\times,\otimes$, $f_{A}\in\calf_{\mu}$ 
and symbols in the language $\calL(\Natural,\in)$.
%Let $|x|=1$ for the variables $x$.
%$|S(t)|=|t|+1$, $|J(s,t)|=|s|+|t|+1$ and $|0^{N}|=|\emptyset|=1$.

For a proof $\mathcal{P}$ with an o.a. $o$,
$|(\mathcal{P},o)|$ denotes the maximum of $|o(\mathcal{P})|$ and
$|A|$ for formulas $A$ occurring in $\mathcal{P}$.
Also $\Natural(\calP)$ denotes the maximum of natural numbers $k$ such that the $k$-th numeral $\bar{k}$
occurs in $\calP$, and 
$S(\calP)$ denotes the maximum of cardinality of the sets $v(s)$ such that the closed set term $s$
occurs in $\calP$.
}
\edf

\bdf\label{df:proofoa}
{\rm A quadruple $(\calP,o,\calf_{0},n)$ of a proof figure $\calP$, an $\calf_{0}$-o.a. (ordinal assignment) 
$o:\Gam\mapsto o(\Gam)\in H(\calf_{0})$, a finite set $\calf_{0}\subset\calf_{\mu}$ and an integer
is an \textit{$(\calf_{0},n)$-proof with o.a.} (ordinal assignment)
if the following conditions are met.}

\bdes
\item[(p0)] 
{\rm $\calP$ is height regulated, and either a $\Sig_{1}$-formula $\exi x\, A(x;t_{1},\ldots,t_{n})$ 
or a $\Pi_{1}$-formula
$\fal x\, \lnot A(x;t_{1},\ldots,t_{n})$ occurs in $\mathcal{P}$, then
the function symbol $f_{A}$ is in the set $\calf_{0}$.
Moreover
$a=o(\mathcal{P})$ is well-behaved, $|(\calP,o)|\leq 2^{2^{n}}$, $\Natural(\calP)\leq 1+n$, and $S(\calP)\leq n$.
}

\item[(p1)]
{\rm For any inference $(ind)$ for induction schema occurring in $\calP$
the induction sequent receives a finite ordinal $a_{1}<\ome$, and
the others $a_{0},a_{2}$ receive the finite ordinal $\dg(A(y))$ for the induction formula $A(y)$, 
cf.\,Definition \ref{df:ordinalassignment}.
}

\item[(p2)]
 \bdes
  
 \item[(p2.1)]
{\rm Let $t$ be a closed term occurring above a $(D_{i})$ with the local stock $c_{i}$.
Then $t<_{i}c_{i}$.}
%\in \calh_{c_{i}}(D_{i}(c_{i}\oplus ))$ where $v(D_{1}(c\oplus))=\Psi_{\rho_{0}}(v(c))$ and
%$v(D_{0}(c\oplus))=\Psi_{\ome_{1}}(v(c))$.}

 \item[(p2.2)]
{\rm Let $J$ be one of rules $(D_{i})_{\alp}$ occurring in $\calP$, 
and $\Gam$ the upper sequent of $J$.
Then $\alp\geq D_{i}(c\oplus o(\Gam))$, 
where $c=stk(o(\Gam))$.}
 \edes
 
\item[(p3)]
{\rm The final part of $\calP$ consists in a $(D_{0})$ followed by a series of paddings,
$(p)_{b_{i}}=(pad)_{b_{i}}$ such that $b_{i}\in H_{0}(\calf_{0})$.}
\[
\calP=
\left.
\begin{array}{c}  
\deduce
{
 \infer[(p)_{b_n}]{\Lam}
 {
  \infer*{\Lam}
  {
   \infer[(p)_{b_0}]{\Lam}
   {
    \infer[(D_{0})]{\Lam}
     {
      \infer*{\Lam}{}
      }
     }
    }
  }
 }
 {}
\end{array}
\right.
\]

{\rm Also there is no $(D_{0})$ above the final $(D_{0})$, i.e., the final is the unique rule $(D_{0})$ in $\calP$,
which is a bottleneck of $\mathcal{P}$.}

\edes
\edf

From {\bf (p3)} and Proposition \ref{prp:H0}.\ref{prp:H0.1}
we see that $dom(o(\Gam_{end}))\in\{0,1,\Natural\}$ for the end-sequent $\Gam_{end}$ of $\mathcal{P}$.

\blem\label{lem:inversion}{\rm (Inversion)}\\
Let $\calP$ be a proof of $\Gam, \fal x\, A(x)$, and $t$ a closed term.
Let $o$ be an o.a. for sequents in $\calP$. Then there exists a proof $\calP'$ of $\Gam,A(t)$ and an o.a. $o'$ such that
$o'(\Gam, A(t))=o(\Gam,\fal x\, A(x))$.

The same holds for proofs ending with $\Gam,A_{0}\land A_{1}$ for conjunctive formulas $A_{0}\land A_{1}$.
\elem
\bprf
For inversion of a universal formula $\fal x\, A(x)$, substitute $t$ for $x$ in $\calP$ to get a proof $\calP'$ of $\Gam,A(t)$.
Consider an $(ind)_{<}$, and let $s$ be the induction term in which $x$ occurs.
Then let
$mj(s[x:=t]):=\rho_{0}=mj(s)$ in the o.a. $o'$ even if $s[x:=t]$ is a closed term.
The same is applied to rules $(ind)_{\Natural}$ and $(ind)_{\in}$.

For inversion of a conjunction $A_{0}\land A_{1}$, replace $A_{0}\land A_{1}$ by $A_{i}$.
Each inference rule $(\land)$ introducing a main formula $A_{0}\land A_{1}$ 
\[
\infer[(\land)]{\Gam,A_{0}\land A_{1}; a_{0}+a_{1}}
{
\Gam,A_{0}; a_{0}
&
\Gam,A_{1}; a_{1}
}
\]
is replaced by one of rules $(pad)_{a_{1-i}}, {}_{a_{1-i}}(pad)$.
\[
\infer[J]{\Gam,A_{i}; a_{0}+a_{1}}
{
\Gam,A_{i}; a_{i}
}
\]
where $J$ is a $(pad)_{a_{1}}$ if $i=0$, and ${}_{a_{0}}(pad)$ otherwise.
\eprf

\blem\label{lem:delta0elim}{\rm (False literal elimination)}\\
Let $A$ be a false closed literal, and
$\calP$ a proof of $\Gam,A$.
Let $o$ be an o.a. for sequents in $\calP$. Then there exists a proof $\calP'$ of $\Gam$ and an o.a. $o'$ such that
$o'(\Gam)=o(\Gam,A)$.
\elem
\bprf
Eliminate the ancestors $A$ of $A$ to get a proof $\calP'$ of $\Gam$.
Consider a $(P\Sig_{1})$.
\[
\infer[(P\Sig_{1})]{\Gam, (\lnot P(t_{0},t_{1}), s\not<t_{0},)  \vphi^{t_{1}}[t_{0},s];a}{\Gam,\vphi[\ome_{1},s];a}
\]
If one of literals $\lnot P(t_{0},t_{1}), s\not<t_{0}$ is a false ancestor of $A$, then eliminate it from the lower sequent.
The same is applied to rules $(P_{\rho_{0}}\Sig_{1}), (b\exi), (ind)_{<},(ind)_{\Natural},(ind)_{\in}$.
\eprf

\subsection{Initial hydras}\label{subsec:embedding}

\blem\label{lem:embed}
Suppose that $T_{1}(\Natural)$ proves a $\Sig_{2}^{*}$-formula $N(x)\to A_{0}(x)$, where
no variable other than $x$ occurs in $A_{0}(x)$.
Then there exist a finite set $\calf_{0}\subset\calf_{\mu}$ and an $\calf_{0}$-o.a. $o$, and for each 
sufficiently large $n\in\Natural$,
there exists
 a proof $\calP_{n}$ of the sequent $\{A_{0}(\bar{n})\}$
such that $(\calP_{n},o,\calf_{0},n)$ is an $(\calf_{0},n)$-proof with o.a.,
and $o(\calP_{n})=o(\calP_{m})$ is an initial hydra.
\elem

Suppose that $T_{1}(\Natural)$ proves a $\Sig_{2}^{*}$-formula $N(x)\to A_{0}(x)$.
We show that there exists a proof $\calP(x)$ of the sequent $\{\lnot N(x),A_{0}(x)\}$ and an $\calf_{0}$-o.a. $o$
such that $(\calP(n),o,\calf_{0},n)$ is an $(\calf_{0},n)$-proof with o.a. for some finite $\calf_{0}$ and 
$o(\calP(x))$ is an initial hydra, where $\calP(n)$ is essentially obtained from $\calP(x)$ by substituting
the numeral $\bar{n}$ for the variable $x$.

Let $\mathcal{Q}_{0}$ be a proof figure of the sequent $\{\lnot N(x),A_{0}(x)\}$ from 
axioms in $T_{1}(\Natural)$.
%axioms (\ref{eq:Z1ord}), (\ref{eq:Z2ord}), (\ref{eq:Z4ord}), (\ref{eq:Z5ord}) and (\ref{eq:Z6ord}),
%and three kinds of induction axiom schemata.

In what follows $\calf_{0}$ denotes the set of function symbols $f_{A}(y_{1},\ldots,y_{n})$ for 
$\Sig_{1}$-formulas $\exi x\, A(x;t_{1},\ldots,t_{n})$ and $\Pi_{1}$-formulas
$\fal x\, \lnot A(x;t_{1},\ldots,t_{n})$ occurring in $\mathcal{Q}_{0}$.

Each leaf in $\mathcal{Q}_{0}$ is either a logical one $(taut)$ or 
one of axioms in $T_{1}(\Natural)$.
% other than $\Pi_{1}$-Collection.
%(\ref{eq:Z1ord}), (\ref{eq:Z2ord}), (\ref{eq:Z4ord}), (\ref{eq:Z5ord}) and (\ref{eq:Z6ord}),
%and three kinds of induction axiom schemata.
Inference rules in $\mathcal{Q}_{0}$ are logical ones, $(\lor),(\land),(\exi)^{N},(\exi),(\fal)$ 
and $(cut)$.

Let us depict pieces of proofs of each leaf in $\mathcal{Q}_{0}$ except $(taut)$'s together with possible ordinal assignments in $\mathcal{I}$ of Definition \ref{df:WB}.

Leaves for axioms (\ref{eq:Z1ord}), (\ref{eq:Z2ord}), (\ref{eq:Z4ord}) and (\ref{eq:Z5ord}) are derived from inference rules 
$(P\Sig_{1})$, $(P\exi)$, $(P_{\rho_{0}}\Sig_{1})$ and $(P_{\rho_{0}}\exi)$, resp.
\[
\infer[(\lor)_{1},(\fal)]{\fal x,y,a(\lnot P(x,y)\lor a\not<x\lor\lnot\vphi[\ome_{1},a]\lor\vphi^{y}[x,a]); 8}
{
 \infer[(P\Sig_{1})]{\lnot P(x,y), a\not<x, \lnot\vphi[\ome_{1},a], \vphi^{y}[x,a]; 2 }
 {
  \infer*{\lnot P(x,y), a\not<x, \lnot\vphi[\ome_{1},a], \vphi[\ome_{1},a];2 }{}
 }
}
\]
with 6 times $(\lor)_{1}$, and $\dg(\vphi)=2$.
\[
\infer[(\lor)_{1},(b\fal)]{\fal a<\ome_{1}\exi x, y<\ome_{1}[a<x\land P(x,y)]; D_{1}(0)}
{
 \infer[(P\exi)]{a\not<\ome_{1},\exi x,y<\ome_{1}[a<x\land P(x,y)]; D_{1}(0)}{}
}
\]
where the formula $\fal a<\ome_{1}(\exi x, y<\ome_{1}[a<x\land P(x,y)])$ is not a $\Del_{0}$-formula.
%cf.\,Definition \ref{df:ordinalassignment}.\ref{df:ordinalassignmentbfal} for rules.

\[
\infer[(\lor)_{1},(\fal)]{\Gam,\fal x,y(\lnot P_{\rho_{0}}(x) \lor y\not<x\lor  \lnot\vphi[y]\lor \vphi^{x}[y]);8}
{
 \infer[(P_{\rho_{0}}\Sig_{1})]{\Gam,\lnot P_{\rho_{0}}(x), y\not<x, \lnot\vphi[y], \vphi^{x}[y]; 2}
 {
 \infer*{\lnot\vphi[y], \vphi[y]; 2}{}
 }
}
\]
with 6 times $(\lor)_{1}$, and $\dg(\vphi)=2$.

\[
\infer[ (\fal)]{\Gam,\fal y\exi x[y<x \land P_{\rho_{0}}(x)]; D_{2}(0)}
{
 \infer[(P_{\rho_{0}}\exi)]{\Gam, \exi x[y<x\land P_{\rho_{0}}(x)]; D_{2}(0)}{}
}
\]
Leaves for transfinite induction schema are replaced as follows.
First consider the schema for ordinals.
{\scriptsize
\[
\hspace{-5mm}
\infer[(\fal),(\lor)_{1}]{\fal \alp(\fal \bet<\alp\, A(\bet)\to A(\alp))\to\fal \alp\, A(\alp); d_{1}\times\rho_{0}+d_{0}+4}
{
\infer[(cut)]{\lnot Prg, A(\gam); d_{1}\times\rho_{0}+d_{0}}
{
 \infer[(b\fal)]{\lnot Prg,\fal \bet<\alp\, A(\bet)}
  {
  \infer[(ind)_{<}]{\bet\not<\alp,\Del; d_{1}\times\rho_{0}}
   {
    \infer[(\land),(\exi)_{1}]{\lnot Prg,\lnot\fal \bet<\gam A(\bet), A(\gam); d_{0}}
    {
    \infer*{\fal \bet< \gam\, A(\bet),\lnot\fal \bet<\gam A(\bet);d}{}
    &
    \infer*{\lnot A(\gam),A(\gam) ; d^{\prime}}{}
    &
    \lnot ON(\gam),ON(\gam); 1
    }
    & 
    \hspace{-2mm}
   \infer*{\Del,A(\bet),\lnot A(\bet);d^{\prime}}{}
    }
   }
  &
  \hspace{-43mm}
  \lnot Prg,\lnot\fal \bet< \gam\, A(\bet),A(\gam); d_{0}
 }
}
\]
}
where $\Del=\lnot Prg\cup\{A(\bet)\}$ with $\lnot Prg=\{\lnot\fal \alp(\fal \bet< \alp\, A(\bet)\to A(\alp))),\lnot ON(\gam)\}$,
$\fal\alp(\cdots):\equiv(\fal \alp(ON(\alp)\to\cdots))$, and 
$d= \dg(\fal \bet<\gam\, A(\bet)), d^{\prime}=\dg(A(\bet))=\max\{d-1,1\}$, $d_{0}=d+d^{\prime}+2$, and
$d_{1}=d_{0}+d^{\prime}+1$.
Also $\rho_{0}=mj(y)$.

Next consider the induction schema for $\Natural$.
%F(0^{N})\land\fal n(F(n)\to F(S(n)))\to\fal n\, F(n)
{\small
\[
\infer[(\lor)_{1},(\fal)]{A(0^{N})\land\fal n(A(n)\to A(Sn))\to\fal n\, A(n);2d+d_{1}\otimes\ome+6}
{
 \infer[(ind)_{\Natural}]{\lnot N(m),\lnot A(0^{N}),\lnot\fal n(A(n)\to A(Sn)),A(m);2d+d_{1}\otimes\ome}
 {
  \infer*{\lnot A(0^{N}),A(0^{N});d}{}
  &
  \hspace{-5mm}
   \infer[(\exi)^{N}_{1}]{\lnot\fal n(A(n)\to A(Sn)),\lnot N(n), \lnot A(n), A(Sn);d_{1}}
   {
    \lnot N(n),N(n);1
     &
    \infer[(\land)]{A(n)\land\lnot A(Sn), \lnot A(n), A(Sn)}
    {
     \infer*{\lnot A(n),A(n);d }{}
     &
     \infer*{\lnot A(Sn),A(Sn); d}{}
     }
   }
  &
  \hspace{-5mm}
  \infer*{\lnot A(m),A(m);d}{}
  }
}
\]
}
where $d=\dg(A(n))$, $d_{1}=2d+2$.

Finally consider the induction schema for sets.
{\small
\[
\infer[(\lor)_{1},(\fal)]{A(\emptyset) \land \fal x,y(Set(x)\land A(x) \to A(J(x,y)))\to \fal x(Set(x)\to A(x));2d+d_{2}\otimes\ome+6}
{
 \infer[(ind)_{\in}]{\lnot Set(y),\lnot A(\emptyset),\lnot\fal x,y(Set(x)\land A(x) \to A(J(x,y))),A(y);2d+d_{2}\otimes\ome}
 {
  \infer*{\lnot A(\emptyset),A(\emptyset);d}{}
  &
 % \hspace{-5mm}
   \infer[(\exi)_{1}]{\lnot\fal x,y(Set(x)\land A(x) \to A(J(x,y))),\lnot Set(x), \lnot A(x), A(J(x,y));d_{2}}
   {
    \infer[(\land)]{Set(x)\land A(x)\land\lnot A(J(x,y)),\lnot Set(x), \lnot A(x), A(J(x,y))}
    {
     \lnot Set(x),Set(x);1
     &
     \infer*{\lnot A(x),A(x);d }{}
     &
     \infer*{\lnot A(J(x,y)),A(J(x,y)); d}{}
     }
   }
  &
  \hspace{-3mm}
  \infer*{\lnot A(y),A(y);d}{}
  }
}
\]
}
where $d=\dg(A(x))$, $d_{2}=2d+3$.
%F(\emptyset) \land \fal x,y(Set(x)\land F(x) \to F(J(x,y)))\to \fal x(Set(x)\to F(x)) $J$ is an $(ind)_{\in}$:

Observe that these pieces enjoy the condition {\bf (p1)}, and there are no nested inference rules for induction schema,
cf.\,{\bf (h4)}.

Leaves for (\ref{eq:Z6ord}) are replaced by
{\small
\[
\infer[(\lor)_{1},(\fal)]{\fal \alp[\fal \bet< \alp\, A(\bet) \to \exi \gam\fal \bet< \alp\,\lnot A^{(\gam)}(\bet)]; 15+D_{2}(0)+4}
{
 \infer[(Rfl)]{\lnot ON(\alp), \lnot \fal \bet< \alp\, A(\bet),\exi \gam\fal \bet< \alp\, A^{(\gam)}(\bet); 15+D_{2}(0)}
 {
  \infer*{ \lnot \fal \bet< \alp\, A(\bet),\fal \bet< \alp\, A(\bet);7}{}
 &
  \infer[(\land),(\exi)_{1}]{\lnot ON(\gam), \alp\not< \gam, \exi \bet< \alp\, \lnot A^{(\gam)}(\bet),\exi \gam\fal \bet< \alp\, A^{(\gam)}(\bet); 8}
  {
  \lnot ON(\gam),ON(\gam);1
  &
   \infer*{\exi \bet< \alp\, \lnot A^{(\gam)}(\bet),\fal \bet< \alp\, A^{(\gam)}(\bet); 6}{}
   }
 }
}
\]
}
where
$7=\dg(\fal \bet< \alp\, A(\bet)), \dg(\fal \bet< \alp\, A^{(\gam)}(\bet))=6$.

Otherwise.
Then there exists a formula $A$ such that 
the formula $\fal\vec{x}\, A$ in the sequent is the universal closure of an axiom in $T_{1}(\Natural)$
not treated so far.
Replace the leaf $\Gam,\fal\vec{x}\, A$ by 
\[
 \infer[(\fal)]{\Gam,\fal\vec{x}\, A}
 {
  \infer[(ax)]{\Gam,A; 1}{}
 }
\]

Next consider inference rules in $\mathcal{Q}_{0}$.
At each $(\lor)$, add $1$, i.e., replace it by $(\lor)_{1}$.
The same for $(\exi)$ introducing an existential formula, and
for $(\exi)^{N}$ with simple instance terms.

Finally consider a $(cut)$:
\[
\infer[(cut)]{\Gam,\Del}
{
\Gam,\lnot A
&
A,\Del
}
\]
Replace it by
\[
\infer[(cut)]{\Gam,\Del;a_{0}+a_{1}}
{
\Gam,\lnot A;a_{0}
&
A,\Del;a_{1}
}
\]

Note that there occurs no inference rules $(D_{i})$ for $i=0,1$ in the constructed $\mathcal{Q}_{0}$.

Let $\mathcal{Q}_{1}$ be the proof of the sequent $\{\lnot N(x),A_{0}(x)\}$ obtained from $\mathcal{Q}_{0}$
as described above with an ordinal $b$ constructed from $1,D_{1}(0),D_{2}(0), n\times\rho_{0}, n\otimes \ome$ and $+$, i.e.,
$b\in\mathcal{I}$ in Definition \ref{df:WB}.

Let $k\geq 6$ be a positive integer such that $k\geq \dg(C)$ for any cut formula $C$ occurring in $\mathcal{Q}_{1}$,
$k\geq\dg(J)$ for any $(ind)\,J$ occurring in $\mathcal{Q}_{1}$,
$3k+11+|b|\leq 2^{2^{k}}$, $\max\{\Natural(\mathcal{Q}_{1}), S(\mathcal{Q}_{1})\}\leq k$
and $|A|\leq 2^{2^{k}}$ for any formula $A$ occurring in $\mathcal{Q}_{1}$.

For each $n\in\Natural$, $\mathcal{Q}_{1}(n)$ denotes a proof of the sequent $\{A_{0}(\bar{n})\}$
obtained from $\mathcal{Q}_{1}$ by substituting the numeral $\bar{n}$ for the variable $x$ and 
eliminating the false literal $\lnot N(\bar{n})$.
Note that $\Natural(\mathcal{Q}_{1}(n))\leq\max\{k,1+n\}$.

\bprp\label{prp:substsize}
\benu
\item\label{prp:substsize.1}
Let $|A|,|t|\leq k$.
Then $|A[x:=t]|\leq k^{2}$ for the result $A[x:=t]$ of substituting the term $t$ for a variable $x$ in the formula $A$.

\item
\label{prp:substsize.2}
Let $\calP'$ be a proof 
obtained from a proof $\calP$ with the restricted rule $(\exi)^{N}$
by substituting a numeral occurring in $\calP$ for a variable.
Then $\Natural(\calP')\leq\Natural(\calP)+1$.

\item
\label{prp:substsize.3}
Let $\calP'$ be a proof 
obtained from a proof $\calP$ with the restricted rule $(\exi)^{S}$
by substituting a closed set term occurring in $\calP$ for a variable.
Then $S(\calP')\leq S(\calP)+1$.
\eenu
\eprp

Add $k$-times $(h)$'s to get a proof $\mathcal{Q}_{2}(n)$:

\[
\mathcal{Q}_{2}(n)=
\left.
\begin{array}{c} 
\infer[(h)]{A_{0}(\bar{n});b_{1}}
{
  \infer*[\mathcal{Q}_{1}(n)]{A_{0}(\bar{n}); b}{}
 }
\end{array}
\right.
\]
where $b_{1}=D_{2}^{(k)}(b)$ with the number $k$ of $(h)$'s.
The conditions {\bf (h3)} and {\bf (h4)} are fulfilled with the proof $\mathcal{Q}_{2}(\bar{n})$.

Next let
\[
\mathcal{P}_{n}=
\left.
\begin{array}{c} 
 \infer[(D_{0})_{\alp_{0}}]{A_{0}(\bar{n});\alp_{0}}
 {
 \infer[(h)]{A_{0}(\bar{n});b_{0}}
 {
  \infer[(D_{1})_{\alp_{1}}]{A_{0}(\bar{n});\alp_{1}}
  {
    \infer*[\mathcal{Q}_{2}(n)]{A_{0}(\bar{n});b_{1}}{}
    }
  }
 }
\end{array}
\right.
\]
where $\alp_{1}=D_{1}(\emptyset\oplus b_{1})$ with the empty stock $\emptyset$,
and another $k$-times $(h)$'s are attached below the $(D_{1})_{\alp_{1}}$.
The conditions {\bf (h2)}, {\bf (h5)} and {\bf (h6)} are fulfilled with the introduced rule $(D_{1})_{\alp_{1}}$.
For {\bf (h5)} note that $k\geq 6=\dg(\fal x<z\, A^{(y)}(x))$ for the formula 
$A^{(y)}(x)\equiv(\exi z<y[P_{\rho_{0}}(z)\land \exi w<y\, B(x)])\, (B\in\Del_{0})$
 in the inference rule $(Rfl)$.
$b_{0}=D_{2}^{(k)}(\alp_{1})$
and $\alp_{0}=D_{0}(c_{0}\oplus b_{0})$ with $c_{0}=D_{2}^{(k+2)}(D_{2}(0)+1)$.
Then {\bf (p2)} is enjoyed for $\mathcal{Q}_{3}$.
For {\bf (p2.1)} note that every closed term $t$ occurring in $\mathcal{Q}_{3}$ is in the closure of 
constants $0^{N},0^{ON},\emptyset,\ome_{1}$ under the function symbols $S, +,\cdot,\lam x.\ome^{x}$ and $J$.
Hence $v(t)$ is in $\calh_{\alp}(\bet)$ for any ordinals $\alp,\bet$.

Thus $\alp_{0}$ is an initial hydra, where the maximum $|(\mathcal{P}_{n},o)|$ of
$|o(\mathcal{P}_{n})|=|\alp_{0}|=3k+11+|b|$ with $b=o(\mathcal{Q}_{1})=o(\mathcal{Q}_{1}(n))$ and the sizes $|A|$ of formulas occurring in
$\mathcal{Q}_{1}(n)$.
Hence $|(\mathcal{P}_{n},o)|\leq (1+n)2^{2^{k}}$ by the choice of the number $k$.

For $n>k$, this $\calP_{n}$ with the o.a. is a proof with o.a. defined in Definition \ref{df:proofoa}.
%Since there occurs no constant other than $0,\ome_{1}$ in $\calP_{0}$, the condition {\bf (p2.1)} holds vacuously.
Obviously we have $a_{0}\in H_{0}(\calf_{0})$.
This shows Lemma \ref{lem:embed}.

\section{Reductions on finite proof figures}\label{sect:consisprf}

\bdf
{\rm Let $A$ be a $\Sig_{2}^{*}$-sentence, and $k$ a natural number.
$k\models A$ iff the result of restricting every unbounded existential $N$-quantifier to $k$ in $A$, i.e.,
restricting $\exi x(N(x)\land\cdots)$ to $\exi x\leq\bar{k}(N(x)\land\cdots)$ holds.
For a finite set $\Gam$ of $\Sig_{2}^{*}$-sentences,
let $k\models\Gam:\Lrarw k\models\bigvee\Gam$.
}
\edf

\blem\label{lem:mainbnd}
Let $(\calP,o,\calf,n)$ be an $(\calf,n)$-proof with an o.a. 
such that
$1+h^{\calf}_{a}(n)\not\models\Gam_{end}$
for $a=o(\mathcal{P})$ and the end-sequent $\Gam_{end}$ of $\Sig_{2}^{*}$-sentences.
Then another $(\calf,n+1)$-proof $(\calP',o',\calf,n+1)$ is constructed such that $a'=o'(\calP')\in (o(\calP))[n]$,
the end-sequent $\Gam_{end}'$ of $\calP'$ is a set of $\Sig_{2}^{*}$-sentences,
and
$1+h^{\calf}_{a'}(n+1)\not\models\Gam_{end}'$.
\elem

Assuming Lemma \ref{lem:mainbnd}, we show Lemma \ref{lem:consis}.
Let $R(x,m)$ be a $\Sig^{*}_{2}$-formula for which
$T_{1}(\Natural)\vdash\fal x[N(x)\to \exi m R(x,m)]$.
By Lemma \ref{lem:embed}
pick a finite set $\calf_{0}\subset\calf_{\mu}$, an $\calf_{0}$-o.a. $o$, and for each 
sufficiently large $n>k$,
 a proof $\calP_{n}$ of the sequent $\{\exi m R(\bar{n},m)\}$
such that $(\calP_{n},o,\calf_{0},n)$ is an $(\calf_{0},n)$-proof with o.a.,
and $a_{0}=o(\calP_{n})=o(\calP_{m})$ is an initial hydra.
Lemma \ref{lem:mainbnd} with the wellfoundedness yields $h^{\calf_{0}}_{a}(n)\models\exi m R(\bar{n},m)$,
i.e.,
$\exi m\leq 1+h^{\calf_{0}}_{a_{0}}(n) R(\bar{n},m)$.
\\

In what follows let $(\calP,o,\calf_{0},n)$ be an $(\calf_{0},n)$-proof with an o.a. $o$ 
such that
$1+h^{\calf}_{a}(n)\not\models\Gam_{end}$ for the end-sequent $\Gam_{end}$ 
of $\Sig_{2}^{*}$-sentences.
We construct another $(\calf_{0},n+1)$-proof $(\calP',o'\calf_{0},n+1)$ such that $a'=o'(\calP')\in (o(\calP))[n]$,
the end-sequent $\Gam_{end}'$ of $\calP'$ is a set of $\Sig_{2}^{*}$-sentences,
and
$1+h^{\calf}_{a'}(n+1)\not\models\Gam_{end}'$.

Note that when a formula $A'$ in $\calP'$ is obtained from a formula $A$ and a term $s$ occurring in $\calP$
by a substitution $A'\equiv A[x:=s]$, then the condition $|A'|\leq (2^{2^{n}})^{2}=2^{2^{n+1}}$ follows 
from $|A|,|s|\leq 2^{2^{n}}$, cf.\,Proposition \ref{prp:substsize}.\ref{prp:substsize.1}.
Also by Proposition \ref{prp:H0}.\ref{prp:H0.1} we have
$|(o(\calP))[n]|\leq  \max\{2^{2^{n}}\cdot 2+2^{2^{n}},2^{2^{n}}\cdot 3,2^{2^{n}}(n+1)\}\leq 2^{2^{n+1}}$ 
if $|o(\calP)|\leq 2^{2^{n}}$.

In each case below the new o.a. $o'$ for the new proof $\calP'$ is defined obviously from the o.a. $o$ and the subscripts $b$ of the displayed inference rules.

\bdf\label{df:mbranchtop}
{\rm The \textit{main branch} of a proof figure $\calP$ is a series $\{\Gam_{i}\}_{i\leq m}$ of occurrences of sequents in $\calP$ such that:}
\benu
\item
{\rm $\Gam_{0}$ is the end-sequent of $\calP$.}
\item
{\rm For each $i<m$, $\Gam_{i+1}$ is the \textit{rightmost} upper sequent of a rule $J_{i}$ with its lower sequent $\Gam_{i}$, and
$J_{i}$ is one of the rules $(cut), (h), (pad)_{0}, {}_{b}(pad)$,
and $(P\Sig_{1}), (P_{\rho_{0}}\Sig_{1}), (D_{i})\, (i=0,1)$.}

\item
{\rm $\Gam_{m}$ is either an axiom or the lower sequent of one of rules 
$(\lor),(\land),(\exi), (b\exi)$,
$(\fal),(b\fal),(ind),(Rfl)$, and $(pad)_{b}$ with $b\neq 0$.}

\eenu
{\rm $\Gam_{m}$ is said to be the \textit{top} of the main branch of $\calP$.}
\edf

Let $\Phi$ denote the top of the main branch of the proof $\calP$ with the o.a. $o$. 
Observe that we can assume that $\Phi$ contains no free variable.

\subsection{top=padding}\label{subsec:toppad}

In this subsection we consider the cases when the top $\Phi$ is a lower sequent of one of rules $(p)_{b}=(pad)_{b}$ with 
$b\neq 0$ or one of rules
$(p)_{b}=(\lor)_{b},(\land)_{b},(\exi)_{b}, (b\exi)_{b}$ with $b>1$.
\\

\noindent
{\bf Case 1}. $dom(b)=\Natural$.
Then $dom(o(\calP))=dom(b)$.
\[
\infer[(p)_{b}]{\Phi; a_{0}+b}
{
 \infer*{\cdots; a_{0}}{}
 }
\]
is replaced by
\[
\infer[(p)_{b[n]}]{\Phi; a_{0}+b[n]}
{
 \infer*{\cdots; a_{0}}{}
 }
\]
The condition {\bf (p2.2)} is fulfilled with the replacement $\mathcal{P}^{\prime}$ by Lemma \ref{prp:H0.2}.
\\

\noindent
{\bf Case 2}. $dom(b)= H_{i}(\calf_{0})$ for $i=0,1$: 
$b=b[D_{i+1}(0)]$ and $dom(o(\calP))=\Natural$. 
Consider the uppermost $(D_{v})\, (v\leq i)$ on the main branch at which $D_{v}$ is applied for hydra.
Such a $(D_{v})$ exists by {\bf (p3)}.
\[
\infer[(D_{v})]{\Gam; D_{v}(c\oplus a[D_{i+1}(0)])}
{
\infer*{\cdots}
 {
  \infer[(p)_{b}]{\Phi;a_{0}+b[D_{i+1}(0)]}
  {
   \infer*{\cdots; a_{0}}{}
   }
  }
 }
\]
where $c=stk(o(\Gam))$.

We have $(D_{v}(c\oplus a))[n]=\ell\cdot v+r$ with $\ell=D_{v}(c\oplus a[1])$ and $r=D_{v}((c+D_{2}(a[1])+1)\oplus a[1])$.
For $v=1$,
\[
 \infer[{}_{\ell}(pad)]{\Gam; (D_{1}(c\oplus a))[n]}
 {
  \infer[(D_{1})]{\Gam; r}
  {
  \infer*{\cdots}
   {
    \infer[(p)_{b[1]}]{\Phi;a_{0}+b[1]}
    {
     \infer*{\cdots; a_{0}}{}
     }
    }
   }
  }
\]
For $v=0$,
\[
  \infer[(D_{0})]{\Gam; r}
  {
  \infer*{\cdots}
   {
    \infer[(p)_{b[1]}]{\Phi;a_{0}+b[1]}
    {
     \infer*{\cdots; a_{0}}{}
     }
    }
   }
\]
In the new $\calP'$, the stock is enlarged to $stk(o(\Gam))=c+D_{2}(a[1])+1$.
The condition {\bf (p2)} is fulfilled with $\calP'$ since $t<_{v}c\Rarw t<_{v}c+D_{2}(a[1])+1$
and $(D_{v}(c\oplus a))[n]<D_{v}(c\oplus a)$.
\\

\noindent
{\bf Case 3}. $b=b_{0}+1$ and $dom(o(\calP))=1$: Then $(p)_{b}=(pad)_{b}$ is below the bottleneck $(D_{0})$ by {\bf (p3)}, and $dom(b)=1=dom(o(\calP))$.

\[
\infer*{\cdots; a+b_{0}+1}
{
 \infer[(p)_{b}]{\Phi; a_{0}+b_{0}+1}
 {
  \infer*{\cdots; a_{0}}{}
  }
 }
\]
$\calP$ is replaced by
\[
\infer*{\cdots; a+b_{0}}
{
 \infer[(p)_{b_{0}}]{\Phi; a_{0}+b_{0}}
 {
  \infer*{\cdots; a_{0}}{}
  }
 }
\]
{\bf Case 4}. $b=b_{0}+1$, $dom(o(\calP))=\Natural$, and $b_{0}\neq 0$ when $(p)_{b}\neq (pad)_{b}$:
Consider the uppermost $(D_{v})\, (v=0,1,2)$ on the main branch at which $D_{v}$ is applied for hydra.
Let $\Gam$ be the lower sequent of the $(D_{v})$, i.e.,
the uppermost sequent below the top $\Phi$ such that $h(\Gam)<h(\Phi)$ or 
$\Gam$ is the lower sequent of a $(D_{i})$ for $i=0,1$.
We have $(D_v(d\oplus(a+b_0+1)))[n]=c\cdot 2$ for $c=D_v(d\oplus (a+b_0))$.
\[\calP=
\left.
\begin{array}{c} 
\infer*{}
{
 \infer{\Gam;\, D_v(d\oplus(a+b_0+1))}
 {
  \infer*{\cdots;\, a+b_{0}+1}
  {
   \infer[(p)_{b}]{\Phi;\, a_{0}+b_0+1}
   {
    \infer*{\cdots;\, a_{0}}{}
   }
  }
 }
}
\end{array}
\right.
\]
Replace the $(p)_{b_{0}+1}$ by $(p)_{b_0}$ and insert a new $(pad)_{c}$ immediately below $\Gam$. 
Note that $c\in H_{0}(\calf_{0})$ when $v=0$, cf.\,{\bf (p3)}.
\[\calP':= 
\left.
\begin{array}{c}  
\deduce
{\hskip0.0cm
 \infer*{}
 {
  \infer[(pad)_{c}]{\Gam;\, c\cdot 2}
  {
   \infer{\Gam;\, D_v(d\oplus(a+b_0))}
   {}
  }
 }
}
{\hskip-1.0cm
  \infer*{\cdots;\, a+b_{0}}
  {
   \infer[(p)_{b_{0}}]{\Phi;\, a_{0}+b_0}
   {
    \infer*{\cdots;\, a_{0}}{}
   }
  }
}
\end{array}
\right.
\]
{\bf Case 5}. 
$dom(b)\in\{multi_{t,2}(\calf_{0}), multi_{t,1}(c_{1};\calf_{0}) :t\in Tm(\calf_{0})\cup\{D_{2}(0)\}, c_{1}\subset H(\calf_{0})\}$.

\[\calP=
\left.
\begin{array}{c} 
\infer*{}
{
 \infer{\Gam;\, D_{0}(c_{0}\oplus a)}
 {
  \infer*{\cdots;\, a}
  {
   \infer[(p)_{b}]{\Phi;\, a_{0}+b}
   {
    \infer*{\cdots;\, a_{0}}{}
   }
  }
 }
}
\end{array}
\right.
\]
where 
$dom(a)\in\{multi_{t,2}(\calf_{0}), multi_{t,1}(c_{1};\calf_{0}) :t\in Tm(\calf_{0}), c_{1}\subset H(\calf_{0})\}$.

Let $multi_{n}$ denote the set in Definition \ref{df:dom}.
Let $s\in multi_{n}$ be a term such that
$D_{0}(c_{0}\oplus a[s])\in (D_{0}(c_{0}\oplus a))[n]$.
Let 
\[\calP'=
\left.
\begin{array}{c} 
\infer*{}
{
 \infer{\Gam;\, D_{0}(c_{0}\oplus a[s])}
 {
  \infer*{\cdots;\, a[s]}
  {
   \infer[(p)_{b[s]}]{\Phi;\, a_{0}+b[s]}
   {
    \infer*{\cdots;\, a_{0}}{}
   }
  }
 }
}
\end{array}
\right.
\]
The condition {\bf (p2.2}) is fulfilled with $\mathcal{P}^{\prime}$ by Lemma \ref{prp:H0.2}.

\subsection{top=axiom}\label{subsec:topaxiom}
In this subsection we consider the cases when the top $\Phi$ is an axiom.
\\

\noindent
{\bf Case 1}.
The top $\Phi=A,\Del_{0}$ is either an $(ax)$ or a $(taut)$.
Then $\Phi$ contains a true $\Del_{0}$-formula $A$ or a literal $A=(\lnot)P(t_{0},t_{1}), (\lnot)P_{\rho_{0}}(t)$.
In each case $ \dg(A)=1$.
\\
{\bf Case 1.1}.
There exists a $(D_{v})$ between the top and the vanishing cut:
Consider the uppermost such $(D_{v})$ at which $D_{v}$ is applied for hydra.
We have $v=1,2$ by {\bf (p3)}.
\[
\infer*{\Gam_{end}:a_{1}}
{
\infer[(cut)]{\Gam,\Del;a+b}
{
 \infer*{\Gam,\lnot A;a}{}
 &
 \infer*{A,\Del;b}
 {
  \infer[(D_{v})]{A,\Del_{1};D_{v}(c_{v}\oplus(b_{0}+1))}
  {
   \infer*{A,\cdots;b_{0}+1}{A,\Del_{0};1}
   }
  }
 }
}
\]
When $A$ is a $\Del_{0}$-formula, let $\mathcal{P}^{\prime}$ be the following with the false $\Del_{0}$-formula $\lnot A$ down to the end-sequent $\Gam_{end}$.
\[
\infer*{\Gam_{end},\lnot A: a_{1}[n]}
{
\infer[(pad)_{b[n]}]{\Gam,\Del,\lnot A;a+b[n]}
{
\infer*{\Gam,\lnot A;a}{}
}
}
\]
where $b[n]=b[D_{v}(c_{v}\oplus b_{0})\cdot 2]$ for $b=b[D_{v}(c_{v}\oplus(b_{0}+1))]$.

Otherwise $A$ is a $P$-literal.
Eliminate the false literal $\lnot A$ by Lemma \ref{lem:delta0elim} to get the following $\mathcal{P}'$.
\[
\infer*{\Gam_{end}: a_{1}[n]}
{
\infer[(pad)_{b[n]}]{\Gam,\Del;a+b[n]}
{
\infer*{\Gam;a}{}
}
}
\]
{\bf Case 1.2}. Otherwise:
Consider the uppermost $(D_{v})\, (v=0,1,2)$ below the vanishing cut.
Such a $(D_{v})$ exists by {\bf (p3)}.
\[\calP= 
\left.
\begin{array}{c}  
\infer*{\Gam_{end}:a_{1}}
{
\infer[(D_{v})]{\cdots;D_{v}(c_{v}\oplus(a_{0}+1))}
{
 \infer*{\cdots;a_{0}+1}
 {
  \infer{\Gam,\Del;\, a+b+1}
  {
   \infer*{\Gam,\lnot A;\, a}{}
  &
   \infer*{A,\Del;\, b+1}{A,\Del_{0};1}
   }
 }
}
}
\end{array}
\right.
\]
where $(D_{v}(c_{v}\oplus(a_{0}+1)))[n]=D_{v}(c_{v}\oplus a_{0})\cdot 2$.
Let $\mathcal{P}^{\prime}$ be the following.
\[\calP':=
\left.
\begin{array}{c}   
\infer*{\Gam_{end}(,\lnot A): a_{1}[n]}
{
\infer[(pad)_{D_{v}(c_{v}\oplus a_{0})}]{\cdots; D_{v}(c_{v}\oplus a_{0})\cdot 2}
{
\infer[(D_{v})]{\cdots;D_{v}(c_{v}\oplus a_{0})}
{
 \infer*{\cdots;a_{0}}
{
 \infer[(pad)_{b}]{\Gam,\Del(,\lnot A);a+b}
 {
  \infer*{\Gam(,\lnot A);\, a}{}
  }
 }
}
}
}
\end{array}
\right.
\]
where $\lnot A$ is absent when $A$ is not a $\Del_{0}$-formula as in {\bf Case 1.1}.
\\

\noindent
{\bf Case 2}. The top is an axiom $(P_{\rho_{0}}\exi)$.

Let $C\equiv(\exi x[t<x\land P_{\rho_{0}}(x)])$.
Consider the uppermost and the lowest $(D_{1})$'s below the $(cut)$ whose cut formula is $C$.
We see that such a $(D_{1})$ exists below the cut from {\bf (h2)}.
\[
\infer[(D_{1})_{\bet}]{\cdots;D_{1}(c_{1}\oplus a)}
{
 \infer*{\cdots}
 {
\infer[(D_{1})_{\alp}]{\Gam^{\prime}; c_{1}\oplus a}
{
 \infer*{\Gam; a[D_{2}(0)]}
 {
  \infer[(cut)]{\Gam_{1},\Del_{1};b_{1}+a_{1}[D_{2}(0)]}
  {
   \infer*{\Del_{1},\lnot C;b_{1}}{}
   &
   \infer*{C,\Gam_{1};a_{1}[D_{2}(0)]}
    {
    \infer[(P_{\rho_{0}}\exi)]{\Gam_{0},\exi x[t<x\land P_{\rho_{0}}(x)]; D_{2}(0)}{}
    }
  }
 }
}
}
}
\]
where 
there is no $(D_{1})$ above the $(cut)$ by {\bf (h2)}.
$(D_{1}(c_{1}\oplus a))[n]=\ell+r$ for
$\ell=D_{1}(c_{1}\oplus a[1])$
and $r=D_{1}(c_{2}\oplus a[1])$ with $c_{2}=c_{1}+D_{2}(a[1])+1$.

We have $t\in\calh_{c_{1}}(D_{1}(c_{1}\oplus a[1]))\cap D_{2}(0)=D_{1}(c_{1}\oplus a[1])=\ell$ by {\bf (p2.1)}.
By inversions for the $A$-formula $\lnot C$ and eliminating false literals $t\not<\ell, \lnot P_{\rho_{0}}(\ell)$ we obtain
the following $\mathcal{P}^{\prime}$, cf.\,Lemma \ref{lem:inversion}.

\[
\infer[{}_{\ell}(pad)]{\Gam^{\prime};\ell+r}
{
\infer[(D_{1})_{\bet}]{\cdots;D_{1}( c_{2}\oplus a[1])(=r)}
{
\infer*{\cdots}
{
 \infer[(D_{1})_{\alp}]{\Gam^{\prime}; c_{2}\oplus a[1]}
 {
  \infer*{\Gam; a[1]}
  {
   \infer[(pad)_{a_{1}[1]}]{\Gam_{1},\Del_{1};b_{1}+a_{1}[1]}
   {
    \infer*[x:=\ell]{\Del_{1};b_{1}}{}
   }
  }
 }
}
}
}
\]
Let us check the condition {\bf (p2.1)} for the $(D_{1})_{\alp}$ in $\mathcal{P}^{\prime}$.
Any term occurring in $\calP'$ is in the closure of $\ell$ and terms occurring in $\calP$
under function symbols in $\mathcal{L}(\Natural,\in)$.
Hence it suffices to show $\ell=D_{1}(c_{1}\oplus a[1])<_{1}c_{2}$,
which follows from $v(c_{1}\oplus a[1])=v(c_{1})\#\ome^{v(a[1])}<v(c_{1})\#\ome^{v(a[1])}+1=v(c_{2})$ and
$\{c_{1}\}\cup E_{1}(a)<_{1}c_{1}$ with $E_{1}(a[1])\subset E_{1}(a)$, which implies $c_{1},a[1]<_{1}c_{2}$.

The condition {\bf (p2.2)} is fulfilled with $\mathcal{P}^{\prime}$ by Lemma \ref{prp:H0.2}, and
{\bf (p0)} by $|\ell|\leq |D_{1}(c_{1}\oplus a)|\leq |o(\calP)|\leq 2^{2^{n}}$.
\\

\noindent
{\bf Case 3}. The top is an axiom $(P\exi)$.

First let $t\not<\ome_{1}$ and $\calP$ be the following.
\[
\infer[(D_{0})_{\alp}]{\Lam;D_{0}(c_{0}\oplus a)}
{
 \infer*{\Lam;a[D_{1}(0)]}
 {
  \infer[(cut)]{\Gam_{1},\Del_{1};b_{0}+a_{0}}
  {
   \infer*{\Del_{1},t<\ome_{1};b_{0}}{}
   &
   \infer*{t\not<\ome_{1},\Gam_{1};a_{0}[D_{1}(0)]}
    {
    \infer[(P\exi)]{\Gam_{0}, t\not<\ome_{1}, \exi x, y< \ome_{1}[t<x\land P(x,y)];D_{1}(0)}{}
    }
  }
 }
}
\]
where there occurs a rule $(D_{0})$ below the $(cut)$, and  
there occurs no $(D_{0})$ above the $(cut)$ by {\bf (p3)}.
We have $(D_{0}(c_{0}\oplus a))[n]=D_{0}(c_{2}\oplus a[1])$ for $c_{2}=c_{0}+D_{2}(a[1])+1$.
Eliminate the false $t<\ome_{1}$ to get the following for the enlarged stock 
$stk(\Lam)=c_{2}=c_{0}+D_{2}(a[1])+1$.
\[
\infer[(D_{0})_{\alp}]{\Lam;D_{0}(c_{2}\oplus a[1])}
{
 \infer*{\Lam;a[1]}
 {
  \infer[(pad)_{a_{0}[1]}]{\Gam_{1},\Del_{1};b_{0}+a_{0}[1]}
  {
   \infer*{\Del_{1};b_{0}}{}
  }
 }
}
\]

Next let $t$ be a closed term such that $v(t)<\ome_{1}$, and
$C\equiv(\exi x, y< \ome_{1}[t<x\land P(x,y)])$.
\[
\infer[(D_{0})_{\alp}]{\Lam;D_{0}(c_{0}\oplus a)}
{
 \infer*{\Lam; a[D_{1}(0)]}
 {
  \infer[(cut)]{\Gam_{1},\Del_{1};b_{1}+a_{1}[D_{1}(0)]}
  {
   \infer*{\Del_{1},\lnot C;b_{1}}{}
   &
   \infer*{C,\Gam_{1};a_{1}[D_{1}(0)]}
    {
    \infer[(P\exi)]{\Gam_{0}(, t\not<\ome_{1}), \exi x, y< \ome_{1}[t<x\land P(x,y)]; D_{1}(0)}{}
    }
  }
 }
}
\]
where there occurs a rule $(D_{0})$ below the $(cut)$, and  
there occurs no $(D_{0})$ above the $(cut)$ by {\bf (p3)}.

Let $(D_{0}(c_{0}\oplus a))[n]=r$ for
$\alp\geq D_{0}(c_{0}\oplus a)>\ell=D_{0}(c_{0}\oplus a[1])$ and 
$r=D_{0}(c_{2}\oplus a[1])$ with $c_{2}=c_{0}+D_{2}(a[1])+1$.
Then $\ell>v(t)$ 
by $t\in\calh_{c_{0}}(D_{0}(c_{0}\oplus a[1]))\cap D_{1}(0)=\ell$, {\bf (p2.1)}.
Let $s=F(c_{0}\oplus a[1])$, i.e., $v(s)=F_{\ell\cup\{\ome_{1}\}}(\rho_{0})$.
We have $|s|=|\ell|\leq|o(\calP)|$ for {\bf (p0)}.

By inversions for the $A$-formula $\lnot C$ and eliminating false literals $\ell\not< \ome_{1}, s\not<\ome_{1}, t\not<\ell, \lnot P(\ell,s)$ we obtain the following, cf.\,Lemmas \ref{lem:inversion} and \ref{lem:delta0elim}.

\[
\infer[(D_{0})_{\alp}]{\Lam;D_{0}(c_{2}\oplus a[1])}
{
 \infer*{\Lam; a[1]}
 {
  \infer[(pad)_{a_{1}[1]}]{\Gam_{1},\Del_{1};b_{1}+a_{1}[1]}
  {
   \infer*[x:=\ell, y:=s]{\Del_{1};b_{1}}{}
   }
  }
 }
\]
Let us check the condition {\bf (p2.1)} for the $(D_{0})_{\alp}$ in $\mathcal{P}^{\prime}$.
Any term occurring in $\calP'$ is in the closure of $\ell, s$ and terms occurring in $\calP$
under function symbols in $\mathcal{L}(\Natural,\in)$.
$s<_{0}c_{2}$ follows from 
$\ell=D_{0}(c_{0}\oplus a[1])<_{0}c_{2}$,
which in turn follows from
$c_{0},a[1]<_{0}c_{0}$
since $D_{0}(c_{0}\oplus a)$ is well-behaved.

The condition {\bf (p2.2)} is fulfilled with $\mathcal{P}^{\prime}$ by Lemma \ref{prp:H0.2}.

\subsection{top=rule}\label{subsec:toprule}

In this subsection we consider the cases when the top $\Phi$ is a lower sequent of 
one of explicit rules $(\lor)_{1},(\land),(\exi)_{1},(b\exi)_{1},(\fal),(b\fal)$ or $(Rfl)$ or one of rules for
induction schema.
% or one of rules $(\lor)_{b},(\exi)_{b}$ with $dom(b)=1$.
\\

\noindent 
{\bf Case 1}.
The top is the lower sequent of an explicit logical rule $J$.
Since the end-sequent consists solely of closed formulas, the main formula of $J$ is also closed.
By virtue of subsection \ref{subsec:toppad} we can assume that $b=1$ for the added hydra $b$ at $J$.
\\

\noindent
{\bf Case 1.1}.
$J$ is a $(\fal)$:
Since the end-sequent consists solely in $\Sig_{2}^{*}$-sentences,
$J$ is a rule introducing unbounded universal quantifier on ordinals.
Consider the uppermost rule $(D_{v})\, (v=0,1,2)$ below $J$ where $D_{v}$ is applied for hydra.
Let $\mathcal{P}$ be the following.
\[
 \infer[(D_{v})]{\Gam,\fal \alp\,\lnot A(\alp);D_{v}(c\oplus (b+1))}
 {
  \infer*{\Gam,\fal \alp\,\lnot A(\alp);b+1}
  {
    \infer[(\fal)\,J]{\Gam_{0},\fal \alp\, \lnot A(\alp); a_{0}+1}
   {
    \infer*{\Gam_{0},\lnot ON(\alp), \lnot A(\alp); a_{0}}{}
  }
 }
}
\]
where $\fal \alp\, \lnot A(\alp)$ is a \textit{false} closed $\Pi_{1}$-formula by the assumption.
Note that the predicate $P_{\rho_{0}}$ does not occur in any $\Del_{0}$-formula, and hence
any $(D_{1})$ does not change the descendants of the formula $\fal \alp\, \lnot A(\alp)$.
$\lnot A(s)$ is a false $\Del_{0}$-formula with the closed term $s=\mu y.A(y)$.
Let $\mathcal{P}^{\prime}$ be the following
with the false $\lnot A(s)$.
\[
\infer[(pad)_{D_{v}(c\oplus b)}]{\Gam,\fal \alp\,\lnot A(\alp),\lnot A(s);D_{v}(c\oplus(b+1))[n]}
{
\infer[(D_{v})]{\Gam,\fal \alp\,\lnot A(\alp),\lnot A(s);D_{v}(c\oplus b)}
{
\infer*{\Gam,\fal \alp\,\lnot A(\alp),\lnot A(s);b}
{
 \infer{\Gam_{0},\fal \alp\, \lnot A(\alp),\lnot A(s);a_{0}}
 {
  \infer*[\alp:=s]{\Gam_{0},\lnot A(s);a_{0}}{}
  }
 }
}
}
\]
where $D_{v}(c\oplus(b+1))[n]=D_{v}(c\oplus b)\cdot 2$,
the closed term $s$ is substituted for the eigenvariable $y$, cf.\,Lemma \ref{lem:inversion},
and the false literal $\lnot ON(s)$ is eliminated by Lemma \ref{lem:delta0elim}.
Note that there is no rule $(D_{1})$ above the rule $(\fal)$
since no free variable occurs below $(D_{1})$ by {\bf (h1)}.

Let us check the condition {\bf (p2.1)} for a rule $(D_{i})$ with its stock $d$ in $\mathcal{P}^{\prime}$.
Let $A(y)\equiv A(y;t_{1},\ldots,t_{k})$. Then
$f_{ A}(t_{1},\ldots,t_{k})=\mu y.\, A(y)\in\calh_{d}(D_{i}(d\oplus ))$
since $\{t_{1},\ldots,t_{k}\}\subset\calh_{d}(D_{i}(d\oplus ))$ by {\bf (p2.1)} in $\mathcal{P}$.
Moreover we have $|s|=|f_{A}(t_{1},\ldots,t_{k})|\leq|A(\alp;t_{1},\ldots,t_{k})|\leq2^{2^{n}}$ for {\bf (p0)}.

The case when $J$ is a $(b\fal)$ with a $\Del_{0}$-main formula is similar.
\\

\noindent
{\bf Case 1.2}.
$J$ is a $(b\fal)$ introducing a bounded universal quantifier for integers: 
\[
\infer*{\Gam,\fal x<t\,A(x);a}
{
 \infer[(b\fal)\,J]{\Gam_{0},\fal x<t\,A(x);a_{0}}
 {
  \infer*{\Gam_{0},x\not<t,A(x);a_{0}}{}
  }
 }
\]
Assume that there is no rule $(D_{v})\,(v=1,2)$ affecting the hydra below $J$.
Since $\fal x<t\,A(x)$ is false, $v(t)>0$ and $t$ is a numeral.
Pick a $k$ such that $k<v(t)$ and $\lnot A(\bar{k})$, where $|\bar{k}|=2k-1<2v(t)-1=|t|$.
Substitute $\bar{k}$ for the variable $x$, and eliminate the false literal $\bar{k}\not<t$,
we obtain
\[
\calP'=
\infer*{\Gam,\fal x<t\,A(x),A(\bar{k});a}
{
  \infer*[x:=\bar{k}]{\Gam_{0},A(\bar{k});a_{0}}{}
 }
\]
\\

\noindent
{\bf Case 1.3}.
$J$ is a $(b\fal)$ introducing a bounded universal quantifier for sets: 
\[
\infer*{\Gam,\fal x\in t\,A(x);a}
{
 \infer[(b\fal)\,J]{\Gam_{0},\fal x\in t\, A(x);a_{0}}
 {
  \infer*{\Gam_{0},x\not\in t, A(x);a_{0}}{}
  }
 }
\]
Assume that there is no rule $(D_{v})\,(v=1,2)$ affecting the hydra below $J$.
Since $\fal x\in t\, A(x)$ is false, $t$ is a closed set term with $v(t)\neq\emptyset$.
Pick a closed term $s$ such that $v(s)\in v(t)$ and $\lnot A(s)$ with $|s|<|t|$.
Substitute $s$ for the variable $x$, and eliminate the false literal $s\not\in t$,
we obtain
\[
\calP'=
\infer*{\Gam,\fal x\in t\, A(x),A(s);a}
{
  \infer*[x:=s]{\Gam_{0},A(s);a_{0}}{}
 }
\]
{\bf Case 1.4}.
$J$ is an $(\exi)_{1}$: 
\[
\infer*{\Gam,\exi y<t\,A^{\prime}(y);a}
{
 \infer[(\exi)\,J]{\Gam_{0},\exi y\,A(y);a_{0}+1}
 {
  \Gam_{0},A(s);a_{0}
  }
 }
\]
where $a=o(\mathcal{P})$, $s$ is a closed term, and
there is a rule $J_{0}$ affecting on a descendant $\exi y\, A(y)$ of the main formula.
$J_{0}$ is one of the rules $(P\Sig_{1})$ and $(P_{\rho_{0}}\Sig_{1})$ 
since the predicate $P_{\rho_{0}}$ does not occur in the end-sequent.
%$\Del_{0}$-formula $\exi y<t\,A^{\prime}(y)$.

If there is a rule $(D_{v})\,(v=1,2)$ between $J$ and $J_{0}$, then 
insert a $(\exi)_{d}$ below the $(D_{v})$, where $d=D_{v}(c\oplus b_{0})$ with $o(\Del)=D_{v}(c\oplus(b_{0}+1))$
for the lower sequent $\Del$ of the $(D_{v})$, cf.\,{\bf Case 4.2} below.
Assume that there is no such rule $(D_{v})$.
\\

\noindent
{\bf Case 1.4.1}.
The rule is a $(P\Sig_{1})$:
Then $(\exi y\,A(y))\equiv(\vphi[\ome_{1},s_{0}])$ and
$(\exi y<t\,A^{\prime}(y))\equiv(\vphi^{t}[t_{0},s_{0}])$ for some closed terms $s_{0},t_{0}$.
\[
\infer*{\Gam,\exi y<t\,A^{\prime}(y);a}
{
\infer[(P\Sig_{1})]{\Gam_{1},(\lnot P(t_{0},t),s_{0}\not<t_{0},) \vphi^{t}[t_{0},s_{0}];b+1}
{
 \infer*{\Gam_{1},\vphi[\ome_{1},s_{0}]; b+1}
 {
 \infer{\Gam_{0},\exi y\,A(y);a_{0}+1}
 {
  \Gam_{0},A(s);a_{0}
  }
 }
}
}
\]
If one of $\lnot P(t_{0},t)$ and $s_{0}\not<t_{0}$ is true, then eliminate 
one of the false literals $P(t_{0},t)$ and $s_{0}<t_{0}$
as in {\bf Case 3} of subsection \ref{subsec:topaxiom}.

Suppose that both $P(t_{0},t)$ and $s_{0}<t_{0}$ are true.
Then $\vphi[\ome_{1},s_{0}]$, i.e., $\exi y\,A(y)$ is false since $\vphi^{t}[t_{0},s_{0}]$ is false.
Hence the closed $\Del_{0}$-formula $A(s)$ is false, too.
Let $\mathcal{P}^{\prime}$ be the following with a $(pad)_{d}$ below the $(P\Sig_{1})$.
Then $o(\mathcal{P}^{\prime})=a[n]$.
Specifically there is a $(D_{v})$ below $(P\Sig_{1})$ at which $D_{v}$ is applied first to hydras.
Then its lower sequent receives $D_{v}(c\oplus (b_{1}+1))$ for a $b_{1}$ in $\mathcal{P}$,
and let $d=D_{v}(c\oplus b_{1})$.
\[
\infer*{\Gam,\exi y<t\,A^{\prime}(y),A(s); a[n]}
{
\infer[(P\Sig_{1})]{\Gam_{1},(\lnot P(t_{0},t),s_{0}\not<t_{0},) \vphi^{t}[t_{0},s_{0}],A(s);b}
{
 \infer*{\Gam_{1},\vphi[\ome_{1},s_{0}],A(s); b}
 {
 \infer{\Gam_{0},\exi y\,A(y),A(s); a_{0}}
 {
  \infer*{\Gam_{0},A(s); a_{0}}{}
  }
 }
}
}
\]
{\bf Case 1.4.2}.
The rule is a $(P_{\rho_{0}}\Sig_{1})$.

Then $(\exi y\,A(y))\equiv(\vphi[s_{0}])$ and
$(\exi y<t\,A^{\prime}(y))\equiv(\exi y<t\,A(y))\equiv(\vphi^{t}[s_{0}])$ for a closed term $s_{0}$.
\[
\infer*{\Gam,\exi y<t\,A(y); a}
{
\infer[(P_{\rho_{0}}\Sig_{1})]{\Gam_{1},(\lnot P_{\rho_{0}}(t),s_{0}\not<t,) \vphi^{t}[s_{0}]; b+1}
{
 \infer*{\Gam_{1},\vphi[s_{0}]; b}
 {
 \infer{\Gam_{0},\exi y\,A(y); a_{0}+1}
 {
  \infer*{\Gam_{0},A(s); a_{0}}{}
  }
 }
}
}
\]
If one of $\lnot P_{\rho_{0}}(t)$ and $s_{0}\not<t$ is true, then eliminate 
one of the false literals $P_{\rho_{0}}(t)$ and $s_{0}<t_{0}$
as in {\bf Case 3} of subsection \ref{subsec:topaxiom}.

Suppose that both $P_{\rho_{0}}(t)$ and $s_{0}<t$ are true.
Then $\vphi[s_{0}]$, i.e., $\exi y\,A(y)$ is false since $\vphi^{t}[s_{0}]$ is false.
Hence the closed $\Del_{0}$-formula $A(s)$ is false, too.
Let $\mathcal{P}^{\prime}$ be the following with a $(pad)_{d}$.
\[
\infer*{\Gam,\exi y<t\,A(y),A(s);a[n]}
{
\infer[(P_{\rho_{0}}\Sig_{1})]{\Gam_{1},(\lnot P_{\rho_{0}}(t),s_{0}\not<t,) \vphi^{t}[s_{0}],A(s); b}
{
 \infer*{\Gam_{1},\vphi[s_{0}],A(s); b}
 {
 \infer{\Gam_{0},\exi y\,A(y),A(s); a_{0}}
 {
  \infer*{\Gam_{0},A(s); a_{0}}{}
  }
 }
}
}
\]
Other cases $(\lor),(\land), (b\exi)$ are similar.
\\

\noindent
{\bf Case 1.5}.
$J$ is a $(\exi)^{N}_{1}$ introducing an unbounded existential quantifier for integers: 
\[
\infer*{\Gam,\exi x(N(x)\land A(x));a}
{
 \infer[(\exi)^{N}_{1}\,J]{\Gam_{0},\exi x( N(x)\land A(x));a_{0}+a_{1}+1}
 {
 \Gam_{0},N(s); a_{0}
 &
  \Gam_{0},A(s);a_{1}
  }
 }
\]
where $s$ is an $N$-simple and closed term.
Assume that there is no rule $(D_{v})\,(v=1,2)$ affecting the hydra below $J$.
If $s$ is not a numeral, then eliminate the false literal $N(s)$.
Suppose that $s$ is a numeral $\bar{k}$.
By the assumption $\exi x\leq 1+h^{\calf}_{a}(n)(N(x)\land A(x))$ does not hold, and
$k\leq\Natural(\calP)\leq 1+n\leq 1+h^{\calf}_{a}(n)$.
Hence $A(s)$ is false.
Let $\calP'$ be the following.
\[
\infer*{\Gam,\exi x(N(x)\land A(x)),A(s);a[n]}
{
 \infer[{}_{a_{0}}(pad)]{\Gam_{0},\exi x( N(x)\land A(x)),A(s);a_{0}+a_{1}}
 {
  \Gam_{0},A(s);a_{1}
  }
 }
\]
{\bf Case 1.6}.
$J$ is a $(\exi)^{S}_{1}$ introducing an unbounded existential quantifier on sets: 
\[
\infer*{\Gam,\exi x(Set(x)\land A(x));a}
{
 \infer[(\exi)^{S}_{1}\,J]{\Gam_{0},\exi x( Set(x)\land A(x));a_{0}+a_{1}+1}
 {
 \Gam_{0},Set(s); a_{0}
 &
  \Gam_{0},A(s);a_{1}
  }
 }
\]
where $s$ is an $S$-simple and closed term.
Assume that there is no rule $(D_{v})\,(v=1,2)$ affecting the hydra below $J$.
If $s$ is not a set term, then eliminate the false literal $Set(s)$.
Suppose that $s$ is a set term.
By the assumption $\exi x(Set(x)\land A(x))$ does not hold, $A(s)$ is false.
Let $\calP'$ be the following.
\[
\infer*{\Gam,\exi x(Set(x)\land A(x)),A(s);a[n]}
{
 \infer[{}_{a_{0}}(pad)]{\Gam_{0},\exi x( Set(x)\land A(x)),A(s);a_{0}+a_{1}}
 {
  \Gam_{0},A(s);a_{1}
  }
 }
\]
{\bf Case 2}.
The top is the lower sequent of a $(Rfl)$:
Let $A(x)\equiv(\exi z\exi w[z\in P_{\rho_{0}}\land B(x)])\, (B\in\Del_{0})$,
$A^{(y)}(x)\equiv(\exi z<y\exi w< y[z\in P_{\rho_{0}}\land B])$.
\[
\infer[J]{\Del; D_{1}(c_{1}\oplus a[D_{2}(0)])}
{
 \infer*{\Del_{1}}
 {
  \infer[J_{1}]{\Del_{2}^{\prime}; a}
  {
   \infer*{\Del_{2}; a}
   {
    \infer[(Rfl)]{\Gam; a_{0}+a_{1}+D_{2}(0)}
     {
     \Gam, \fal x< t\, A(x); a_{0}
     &
     t\not< y, \exi x< t \lnot A^{(y)}(x),\Gam; a_{1}
     }
    }
  }
 }
}
\]
where 
$J_{1}$ is the uppermost $(D_{1})_{\alp_{1}}$ and $J$ is the lowermost $(D_{1})_{\alp}$ below the $(Rfl)$.
Such a $(D_{1})$ exists by {\bf (h5)}.

We have 
$t\in\calh_{c_{1}}(D_{1}(c_{1}\oplus a[1] ))$ by {\bf (p2.1)},
and hence $v(t)<\ell=D_{1}(c_{1}\oplus a[1])$.
$(D_{1}(c_{1}\oplus a))[n]=\ell+r$ for $r=D_{1}(c_{2}\oplus a[1])$
with $c_{2}=c_{1}+D_{2}(a[1])+1$.
Let $\mathcal{P}'$ be the following.
\[
\infer[(cut)]{\Del; \ell+r}
{
 \infer{\Del,\fal x< t A^{(\ell)}(x); D_{1}(c_{1}\oplus a[1])}
 {
  \infer*{\Del_{1},\fal x< t A^{(\ell)}(x); a[1]}
  {
   \infer{\Del_{2}^{\prime},\fal x< t A^{(\ell)}(x); a[1]}
   {
    \infer[(D_{1})_{\ell}]{\Del_{2},\fal x< t A^{(\ell)}(x); a[1]}
    {
     \infer*{\Del_{2},\fal x< t A(x);a[1]}
     {
      \infer[(pad)_{a_{1}+1}]{\Gam,\fal x< t A(x); a_{0}+a_{1}+1}
      {\Gam,\fal x< t A(x); a_{0}}
      }
      }
     }
   }
  }
&
 \infer{\exi x< t\lnot A^{(\ell}(x),\Del; D_{1}(c_{2}\oplus a[1])}
 {
  \infer*{\exi x< t\lnot A^{(\ell)}(x),\Del_{1}; a[1]}
  {
   \infer[(D_{1})_{\alp_{1}}]{\exi x< t\lnot A^{(\ell)}(x),\Del_{2}^{\prime}; a[1]}
   {
    \infer*{\exi x< t\lnot A^{(\ell)}(x),\Del_{2};a[1]}
      {
      \infer[{}_{a_{0}}(pad)_{1}]{\exi x< t\lnot A^{(\ell}(x),\Gam; a_{0}+a_{1}+1}
      {
       \infer*[y:=\ell]{\exi x< t\lnot A^{(\ell)}(x),\Gam; a_{1}}{}
       }
     }
    }
  }
 }
}
\]
In $\mathcal{P}$,  $h(\Del)\geq \dg(\exi x< t \lnot A^{(y)}(x))= \dg(\fal x< t A^{(b)}(x))$ by {\bf (h5)}.
Thus the introduced $(cut)$ in $\mathcal{P}'$ enjoys {\bf (h3)}.
There is no $(D_{1})$ above the $(Rfl)$ by {\bf (h2)}.
In the left part of the $(cut)$, a new $(D_{1})_{\ell}$ arises with its stock $c_{1}$ and $\ell=D_{1}(c_{1}\oplus a[1])$, 
cf.\,{\bf (p2.2)}.
In the upper sequent of the right rule $(D_{1})_{\alp_{1}}$, 
a bounded sentence $\exi x< t\lnot A^{(\ell)}(x)$ is added, cf.\,the definition of the rule $(D_{1})$.
For the condition {\bf (p2.1)} of the right rule $(D_{1})_{\alp_{1}}$ we have
$\ell=D_{1}(c_{1}\oplus a[1])<_{1}c_{2}$
by $c_{1}, a[1]<_{1}c_{1}$ and $v(c_{1}\oplus a[1])=v(c_{1})\#\ome^{v(a[1]}=v(c_{1}+D_{2}(a[1]))<v(c_{2})$.
Moreover we have $|1|=|D_{0}(0)|=|D_{2}(0)|$ and
$|\ell|=|D_{1}(c_{1}\oplus a[1])|=|D_{1}(c_{1}\oplus a[D_{2}(0)]|\leq|o(\calP)|$ for {\bf (p0)}.
\\

\noindent
{\bf Case 3}.
The top is the lower sequent of an $(ind)$.
\\
{\bf Case 3.1} The top is the lower sequent of an $(ind)_{<}$.
\[
\infer[(D_{0})]{\cdots; D_{0}(c_{0}\oplus b_{0})}
{
\infer*{\cdots;b_{0}}
{
\infer[(D_{1})]{\cdots;c_{1}\oplus b_{1}}
{
\infer*{\cdots;b_{1}}
{
\infer[(ind)_{<}]{(s\not< t,) \Gam;a}{\Gam,\lnot\fal x< y A(x), A(y);a_{1} & \Gam,\lnot A(s);a_{2}}
}
}
}
}
\]
where $(D_{1})$ is the uppermost one.
Such a $(D_{1})$ exists by {\bf (h4)}.
There is no $(D_{1})$ above the $(ind)$ by {\bf (h2)}.
By {\bf (p1)} we have $ \dg(A(y))=a_{2}$, $a_{1}<\ome$ and
$a=(a_{1}+a_{2}+1)\times mj(t)$,
where 
$\rho_{0}\geq mj(t)\geq v(t)$ for the closed term $t$.
%$t\in Tm(\calf_{0})$.
Also
$|s|\leq 2^{2^{n}}$.
% for $s\in Tm(\calf_{0})$.
\\
{\bf Case 3.1.1}. $s\not< t$:
Then the true literal $s\not< t$ remains in the lower sequent.
Eliminate the false literal $s<t$, and insert a $(pad)_{d}$ to have 
$o(\mathcal{P}^{\prime})=o(\mathcal{P})$.
Then we are in {\bf Case 5} of subsection \ref{subsec:toppad}.
\\

\noindent
{\bf Case 3.1.2}. $s< t$:
Then $t$ is a closed ordinal term, i.e., $t\in Tm(\calf_{0})$ with $v(t)>0$, 
and $s$ is either $s\in Tm(\calf_{0})$ or not a well formed term.
Let 
$multi_{n}=\{s_{0}\in Tm(\calf_{0}): v(s_{0})<v(mj(t)), |s_{0}|\leq 2^{2^{n}}, \bigwedge_{i=0,1}s_{0}\in\calh_{c_{i}}(D_{i}(c_{i}\oplus))\}$. 
Then $0^{ON}\in multi_{n}$.
Let
$s^{\prime}$ be the term such that 
if $s\in Tm(\calf_{0})$, then $s^{\prime}\equiv s$, and $s^{\prime}\equiv 0^{ON}$ otherwise.
Then $s^{\prime}\in multi_{n}$ and $s= s^{\prime}< mj(t)$ holds {\bf (p2.1)}.
Let $mj(s)=s^{\prime}$.

Assuming $\lnot A(s)$ is an $\exi$-formula, 
let $P'$ be the following: 
{\footnotesize
\[\hspace*{-7mm}
\infer[(D_{0})]{\cdots; D_{0}(c_{0}\oplus b_{0}')}
{
\infer*{\cdots; b_{0}'}
{
\infer[(D_{1})]{\cdots; c_{1}\oplus b_{1}'}
{ 
\infer*{\cdots;b_{1}'}
{
 \infer[(cut)]{(s\not< t,) \Gam;a'}
 {
  \infer[(b\fal)]{\Gam,\fal x< s A(x)}
   {
    \infer[(ind)_{<}]{\Gam,y\not< s,A(y);(a_{1}+a_{2}+1)\times mj(s)}
     {
      \infer*{\Gam,\lnot\fal x< y A(x),A(y);a_{1}}{}
      &
       \infer*[P(A)]{\Gam,\lnot A(y),A(y);a_{2}}{}
      }
    }
   &
    \infer[(cut)]{\lnot\fal x< s\, A(x),\Gam}
     {
      \infer*[y:=s]{\Gam,\lnot\fal x< s\, A(x),A(s);a_{1}}{}
      &
      \infer*{\lnot A(s),\Gam;a_{2}}{}
     }
  }
 }
}
}
}
\hspace*{-10mm}
P'
\]
}
where $P(A)$ denotes a proof of $\Gam,\lnot A(y),A(y)$ which is canonically constructed from logical inferences, 
cf.\,Tautology lemma \ref{lem:depth}.
For the part of the substitution $[y:=s]$, cf.\,Lemma \ref{lem:inversion}.

We have $h_{0}(\Gam)\geq \dg(\fal x< a A(x))\geq \dg(A(a))$ by {\bf (h4)}, and hence {\bf (h3)} holds for the introduced $(cut)$'s.
Also $a'=(a_{1}+a_{2}+1)\times mj(s)+a_{1}+a_{2}=a[mj(s)]$, and $D_{0}(c_{0}\oplus b_{0}')\in (D_{0}(c_{0}\oplus b_{0}))[n]$.
Since no essentially new term is created here, {\bf (p2.1)} is fulfilled with $\mathcal{P}^{\prime}$.

If $\lnot A(s)$ is not an $\exi$-formula, then upper sequents of the upper cut should be interchanged.
Note that $a_{2}+a_{1}=a_{1}+a_{2}$ for $a_{1},a_{2}<\ome$:
\[
\infer{(s\not<t, ) \Gam}
{
 \infer*{\Gam,\fal x< s A(x)}{}
&
 \infer{\lnot \fal x< s A(x),\Gam}
 {
  \infer*{\lnot A(s),\Gam;a_{2}}{}
 &
  \infer*{\Gam,\lnot\fal x< s A(x),A(s);a_{1}}{}
 }
}
\msten P'
\]
{\bf Case 3.2} The top is the lower sequent of an $(ind)_{\Natural}$.
\[
\infer[(ind)_{\Natural}]{(\lnot N(s),) \Gam;a}
 {
  \Gam,A(\bar{0});a_{0} 
  & 
  \Gam,\lnot N(x),\lnot A(x),A(Sx);a_{1}
  &
  \Gam,\lnot A(s);a_{2}
  }
\]
where $a_{0}=a_{2}=\dg(A)\leq h_{0}((\lnot N(s),),\Gam)$ and $a=a_{0}+a_{2}+a_{1}\otimes \ome$ with $a_{1}>0$.
Then $a[n]=a_{0}+a_{2}+a_{1}\cdot(n+1)\geq 1$.

If $\lnot N(s)$ is a true literal,
then let $\calP'$ be the following:
\[
\infer[(pad)]{\lnot N(s), \Gam;a[n]}
 {
  \lnot N(s),\Gam; 1
  }
\]
Let $s\equiv\bar{k}$ be a numeral.
Then $k\leq 1+n$. 
Let $\calP'$ be the following.
\[
\infer[(pad)_{a_{1}\cdot(1+n-k)}]{\Gam,a[n]}
{
\infer[(cut)]{\Gam;a_{0}+a_{2}+a_{1}\cdot k}
 {
  \Gam,A(\bar{0});a_{0} 
  & 
  \infer*[x:=\bar{m}]{\Gam,\lnot A(\bar{m}),A(\ovl{m+1});a_{1}\,(m<k)}{}
 &
 \Gam,\lnot A(s);a_{2}
 }
}
\]
{\bf Case 3.3} The top is the lower sequent of an $(ind)_{\in}$.
\[
\infer[(ind)_{\in}]{(Set(s),) \Gam;a}
 {
 \Gam,A(\emptyset);a_{0}
 &
 \Gam,\lnot Set(x),\lnot A(x), A(J(x,y));a_{1} 
 & 
 \Gam,\lnot A(s);a_{2}
 }
\]
where $a_{0}=a_{2}=\dg(A)\leq h_{0}((\lnot Set(s),),\Gam)$ and $a=a_{0}+a_{2}+a_{1}\otimes \ome$ with $a_{1}>0$.
Then $a[n]=a_{0}+a_{2}+a_{1}\cdot(1+n)$.

If $\lnot Set(s)$ is a true literal,
then let $\calP'$ be the following:
\[
\infer[(pad)]{\lnot Set(s), \Gam;a[n]}
 {
  \lnot Set(s),\Gam; 1
  }
\]
Let $v(s)=\{v(t_{0}),\ldots,v(t_{k-1})\}$ with subterms $t_{i}$ of $s\equiv s_{k}$ such that
$s_{0}\equiv\emptyset$ and $s_{i+1}\equiv J(s_{i},t_{i})$ for $i<k$.
Then $v(s_{i})=\{v(t_{j}):j<i\}$ and $k\leq n$.
Let $\calP'$ be the following.
\[
\infer[(pad)_{a_{1}\cdot(1+n-k)}]{\Gam,a[n]}
{
\infer[(cut)]{\Gam;a_{0}+a_{2}+a_{1}\cdot k}
 {
  \Gam,A(\emptyset);a_{0} 
  & 
  \infer*[x:=s_{i},y:=t_{i}]{\Gam,\lnot A(s_{i}),A(s_{i+1});a_{1}\,(i<k)}{}
 &
 \Gam,\lnot A(s);a_{2}
 }
}
\]
{\bf Case 4}.
The top $\Phi$ is the lower sequent of one of logical inferences $(\lor)_{b}, (\exi)_{b}, (b\exi)_{b}$.
By virtue of subsection \ref{subsec:toppad} we can assume $b=1$.
Consider the cases when the logical inference is one of $ (\exi)_{1},(b\exi)_{1}$, which is denoted $(\exi)_{1}$.
The case $(\lor)_{1}$ is similar.
Let the main formula of the logical inference be a formula $\exi x< t \, A(x)$ with a minor formula $A(s)$,
where $t$ denotes either a term or $\rho_{0}$, $(\exi x<\rho_{0}\,A(x)):\equiv(\exi x\, A(x))$.
Let $J$ denote the $(cut)$ at which the descendant $\exi x< t'\, A'(x)$ of $\exi x<t\, A$ vanishes.
\\

\noindent
{\bf Case 4.1}.
$\exi x< t'\, A'(x)$ is a $\Del_{0}$-formula:
Let $\mathcal{P}$ be the following.
\[
 \infer*{\Gam_{end};c}
 {
  \infer[(cut)]{\Gam,\Del;a+b}
  {
   \infer*{\Gam,\lnot \exi x< t^{\prime} \, A^{\prime}(x);a}{}
   &
   \infer*{\exi x< t^{\prime} \, A^{\prime}(x),\Del;b}
   {
    \infer[(\exi)_{1}]{\exi x< t\, A(x),\Del_{0};b_{0}+1}{\exi x< t\, A(x),A(s),\Del_{0};b_{0}}
   }
  }
 }
\]

One of $\lnot \exi x<t^{\prime}\, A^{\prime}(x),\exi x<t^{\prime}\, A^{\prime}(x)$ is false.
When $\exi x<t^{\prime}\, A^{\prime}(x)$ is false, let the false $\Del_{0}$-formula $\exi x<t^{\prime}\, A^{\prime}(x)$
go down to the end-sequent.
\[
 \infer*{\Gam_{end},\exi x<t^{\prime}\, A^{\prime}(x);c}
 {
  \infer[{}_{a}(pad)]{\Gam,\Del,\exi x<t^{\prime}\, A^{\prime}(x);a+b}
  {
   \infer*{\exi x<t^{\prime}\, A^{\prime}(x),\Del;b}
      {
    \infer[(\exi)_{1}]{\exi x< t\, A(x),\Del_{0};b_{0}+1}{\exi x< t\, A(x),A(s),\Del_{0};b_{0}}
   }
  }
 }
\]
This is in {\bf Case 1.2} of this subsection.

When $\exi x<t^{\prime}\, A^{\prime}(x)$ is true, we are in {\bf Case 1.1} of this subsection.
\[
 \infer*{\Gam_{end},\lnot \exi x<t^{\prime}\, A^{\prime}(x);c}
 {
  \infer[(pad)_{b}]{\Gam,\Del,\lnot\exi x<t^{\prime}\, A^{\prime}(x);a+b}
  {
  \infer*{\Gam,\lnot\exi x<t^{\prime}\, A^{\prime}(x);a}{}
  }
  }
\]

In what follows assume that $\exi x< t'\, A'(x)$ is not a $\Del_{0}$-formula
\\

\noindent
{\bf Case 4.2}.
The case when there exists a $(D_{v})\,(v=1,2)$ between $\Phi$ and $J$ at which $D_{v}$ is applied to hydras.
Consider the uppermost such $(D_{v})$.
\[
  \infer[(D_{v})]{\exi x< t'\, A'(x),\Del'; D_{v}(c_{v}\oplus(a+1))}
  {
   \infer*{\cdots;a+1}
   {
    \infer[(\exi)_{1}]{(s\not<t,)\exi x< t\, A(x),\Del_{0};a_{0}+1}{\exi x< t\, A(x),A(s),\Del_{0};a_{0}}
   }
  }
\]
where $\exi x<t^{\prime}\,A^{\prime}(x)$ may differ from $\exi x<t\,A(x)$ due to a rule $(D_{1})$
with $t=\rho_{0}$
when either $(\exi x<t\,A(x))\equiv(\exi x\exi w[x\in P_{\rho_{0}}\land B(x,w)])\, (B\in\Del_{0})$,
or
$(\exi x<t\,A(x))\equiv(\exi x[s\in P_{\rho_{0}}\land B(s,x)])$.
The case when a rule $(P\Sig_{1}), (P_{\rho_{0}}\Sig_{1})$ change a descendant of the main formula is excluded
since we are assuming that $\exi x< t'\, A'(x)$ is not a $\Del_{0}$-formula.

Then lower the $(\exi)$ below the $(D_{v})$.
\[
  \infer[(\exi)_{D_{v}(c_{v}\oplus a)}]{(s\not<t',)\exi x< t'\, A'(x),\Del'; D_{v}(c_{v}\oplus a)\cdot 2}
   {\infer[(D_{v})]{\exi x< t'\, A'(x),A'(s),\Del'; D_{v}(c_{v}\oplus a)}
   {
    \infer*{\cdots;a}
    {
     \exi x< t\, A(x),A(s),\Del_{0};a_{0}
    }
   }
  }
\]
We have to verify that this is a legitimate proof.
Assume that $v=1$ and 
there exists a rule $(D_{1})_{\alp}$ affecting $(\exi x<t\,A(x))\equiv(\exi x\exi w[x\in P_{\rho_{0}}\land B(x,w)])$, and
$(\exi x<t^{\prime}A^{\prime}(x))\equiv (\exi x<\alp\exi w<\alp[x\in P_{\rho_{0}}\land B(x,w)])$ with $t=\rho_{0}$.
We have $s\in\calh_{c_{1}}(D_{1}(c_{1}\oplus ))$ by {\bf (p2.1)} and $D_{1}(c_{1}\oplus(a+1))\leq\alp$ by {\bf (p2.2)}.
Hence $v(s)<D_{1}(c_{1}\oplus(a+1))\leq\alp=t^{\prime}$.
Also note that the new $(\exi)_{D_{v}(c_{v}\oplus a)}$ does not divide a series of rules $(D_{1})$
since it is inserted below the lowest, cf.\,{\bf (h2)}.
\\

\noindent
{\bf Case 4.3}.
By virtue of {\bf Case 4.2} we can assume that there is no $(D_{v})$ between $\Phi$ and
the vanishing cut at which $D_{v}$ is applied to hydras for $v=1,2$.
Then the descendants of the main formula $\exi x<t\,A(x)$ does not change up to the cut formula $\exi x<t\,A(x)$.
Note that there is no $(D_{1})$ nor $(D_{0})$ above the $(cut)\, J$ by {\bf (h2)}, and there is a $(D_{0})$ below the vanishing cut by {\bf (p3)}.
Consider the uppermost $(D_{v})\,(v=0,1,2)$ at which $D_{v}$ is applied to hydras.
\[
\infer[(D_{v})]{\Lam;D_{v}(d_{v}\oplus(c+1))}
{
 \infer*{\cdots;c+1}
 {
  \infer[(cut)]{\Gam,\Del;a+b+1}
  {
   \infer*{\Gam,\lnot \exi x< t \, A(x);a}{}
   &
   \infer*{\exi x< t \, A(x),\Del;b+1}
   {
    \infer[(\exi)_{1}]{\exi x< t\, A(x),\Del_{0};b_{0}+1}{\exi x< t\, A(x),A(s),\Del_{0};b_{0}}
   }
  }
 }
}
\]
Since $\exi x<t\,A(x)$ is not a $\Del_{0}$-formula,
$ \dg(\exi x< t \, A(x))>0$,
and there exists an $(h)$ below the vanishing cut by {\bf (h3)}.
This means that the rule $(D_{v})$ is an $(h)=(D_{2})$ and $v=2$.
Hence $D_{v}(d_{v}\oplus(c+1))=D_{2}(c+1)$.

Since $h_{0}(\Gam,\Del)\geq \dg(\exi x< t \, A(x))> \dg(A(s))$, we have 
$h_{0}(\Lam)=h_{0}(\Gam,\Del)-1\geq \dg(A(s))$ for {\bf (h3)}.
Assuming that $\lnot A(s)$ is an $E$-formula, let $\calP'$ be the following for $D_{2}(c+1)[n]=D_{2}(c)\cdot 2$, cf.\,Lemma \ref{lem:inversion}.

\[
\infer[(cut)]{\Lam; D_{2}(c)\cdot 2}
{
 \infer[(h)]{\Lam,A(s); D_{2}(c)}
 {
  \infer*{\cdots;c}
  {
   \infer{\Gam,\Del,A(s);a+b}
   {
    \infer*{\Gam,\lnot \exi x< t \, A(x);a}{}
    &
    \infer*{\exi x< t\, A(x),A(s)\Del;b}
    {
     \exi x< t\, A(x),A(s),\Del_{0};b_{0}
    }
   }
  }
 }
&
\infer[(h)]{\lnot A(s),\Lam;D_{2}(c)}
{
 \infer*{\cdots;c}
 {
  \infer[(pad)_{b}]{\lnot A(s),\Gam,\Del;a+b}
  {
   \infer*[x:=s]{\lnot A(s),\Gam;a}{}
   }
  }
 }
}
\]
Note that there may occur a $(D_{1})$ above the left part of the $(cut)$ in $\mathcal{P}$.
Let $(D_{1})_{\alp}$ be a rule occurring above the left upper sequent of the $(cut)$ such that
its lower sequent contains an ancestor $\lnot\exi x<t\,A(x)$ of the left cut formula.
We have to verify the condition {\bf (p2.1)} for the $(D_{1})$ in $\mathcal{P}^{\prime}$.
Let $c_{1}$ be the local stock of the $(D_{1})_{\alp}$.
Then $t<_{1}c_{1}$, where
$t<D_{2}(0)=\rho_{0}$ since 
an implicit formula $\lnot\exi x<t\,A^{\prime}(x)$ is in the upper sequent of the $(D_{1})_{\alp}$
where either $A^{\prime}\equiv A$ or $(A^{\prime})^{(\alp)}\equiv A$,
and there occurs no unbounded universal quantifier in an implicit formula in an upper sequent of a rule $(D_{1})$
by the definition of the rule.
Hence $s<t<D_{1}(c_{1}\oplus d)$ for any $d$.
Thus $s\in \calh_{c_{1}}(D_{1}(c_{1}\oplus d))$, i.e., $s<_{1}c_{1}$. This shows {\bf (p2.1)}.
\\

The case when the top $\Phi$ is the lower sequent of one of logical inferences $(\exi)^{N}_{1},(\exi)^{S}_{1}$
is similar.

This completes a proof of Lemma \ref{lem:mainbnd}, and hence of Lemma \ref{lem:consis} and
Theorem \ref{th:hydra}.\ref{th:hydra1}.

\end{document}